\documentclass[twoside,american,english]{article}
\usepackage[T1]{fontenc}
\usepackage[latin9]{luainputenc}
\usepackage{geometry}
\usepackage{float}
\usepackage{url}
\usepackage{amsmath}
\usepackage{amsthm}
\usepackage{amssymb}
\usepackage{graphicx}

\makeatletter

\providecommand{\tabularnewline}{\\}

\numberwithin{equation}{section}
\numberwithin{figure}{section}
\newcommand{\lyxaddress}[1]{
\par {\raggedright #1
\vspace{1.4em}
\noindent\par}
}
\theoremstyle{plain}
\newtheorem{thm}{\protect\theoremname}
\theoremstyle{definition}
\newtheorem{example}[thm]{\protect\examplename}
\theoremstyle{plain}
\newtheorem{cor}[thm]{\protect\corollaryname}

\makeatother

\usepackage{babel}
\addto\captionsamerican{\renewcommand{\corollaryname}{Corollary}}
\addto\captionsamerican{\renewcommand{\examplename}{Example}}
\addto\captionsamerican{\renewcommand{\theoremname}{Theorem}}
\addto\captionsenglish{\renewcommand{\corollaryname}{Corollary}}
\addto\captionsenglish{\renewcommand{\examplename}{Example}}
\addto\captionsenglish{\renewcommand{\theoremname}{Theorem}}
\providecommand{\corollaryname}{Corollary}
\providecommand{\examplename}{Example}
\providecommand{\theoremname}{Theorem}

\begin{document}

\title{Computational Results for the Higgs Boson Equation in the de Sitter
Spacetime}

\author{Andras~Balogh; Jacob Banda; Karen Yagdjian}
\maketitle

\lyxaddress{\begin{center}
The University of Texas Rio Grande Valley, Edinburg TX, 78539, United
States
\par\end{center}}
\begin{abstract}
High performance computations are presented for the Higgs Boson Equation
in the de Sitter Spacetime using explicit fourth order Runge-Kutta
scheme on the temporal discretization and fourth order finite difference
discretization in space. In addition to the fully three space dimensional
equation its one space dimensional radial solutions are also examined.
The numerical code for the three space dimensional equation has been
programmed in CUDA Fortran and was performed on NVIDIA Tesla K40c
GPU Accelerator. The radial form of the equation was simulated in
MATLAB. The numerical results demonstrate the existing theoretical
result that under certain conditions bubbles form in the scalar field.
We also demonstrate the known blow-up phenomena for the solutions
of the semilinear Klein-Gordon equation with imaginary mass. Our numerical
studies suggest several previously not known properties of the solution
for which theoretical proofs do not exist yet: 1. smooth solution
exists for all time if the initial conditions are compactly supported
and smooth; 2. under some conditions no bubbles form; 3. solutions
converge to step functions related to unforced, damped Duffing equations.
\end{abstract}

\section{\label{sec:Introduction}Introduction}

There are several open mathematical questions about the Higgs boson
in the de Sitter spacetime. We are interested in the feature of this
issue related to the partial differential equations theory and, especially,
to the problem of the global in time existence of the solution. In
fact, the equation is semilinear since it contains the Higgs potential
and has time dependent coefficient. Because of the lack of mathematically
rigorous proof of the existence of global in time solution we turn
to numerical investigations which can shed a light on that issue,
and also to indicate the creation of so-called bubbles. In order to
achieve that aim, in this article we perform high-performance numerical
computations using Graphical Processing Units for examining the behavior
of solutions to the Higgs boson equation in the de Sitter spacetime. 

The Klein-Gordon equation with the Higgs potential (the Higgs boson
equation) in the de Sitter spacetime can be written as
\begin{align}
\phi_{tt}-e^{-2t}\Delta\phi+3\phi_{t} & =\mu^{2}\phi-\lambda\phi^{3}\label{eq:HBE}\\
\phi\left(\vec{x},0\right) & =\varphi_{0}\left(\vec{x}\right),\quad\vec{x}\in\mathbb{R}^{3},\\
\phi_{t}\left(\vec{x},0\right) & =\varphi_{1}\left(\vec{x}\right),\quad\vec{x}\in\mathbb{R}^{3}.\label{eq:HBEinit}
\end{align}
Here $\Delta$ is the Laplace operator in $\vec{x}\in\mathbb{R}^{3}$,
$t>0$ is the time variable, the parameters are $\lambda>0$ and $\mu>0$,
while $\phi=\phi\left(\vec{x},t\right)$ is a real-valued function.
From now on we consider solution $\phi=\phi\left(\vec{x},t\right)$
to the Klein-Gordon equation which is at least continuous. It is of
considerable interest for particle physics and inflationary cosmology
to study the so-called bubbles \cite{coleman_1985,Linde,Voronov}.
In the quantum field theory a bubble is defined as a simply connected
domain surrounded by a wall such that the field approaches one of
the vacuums outside of a bubble (see, e.g., \cite{coleman_1985}).
It is mathematically reasonable to define the bubble as a maximal
connected set of spatial points $\vec{x}\in\mathbb{R}^{3}$ at which
solution to the Cauchy problem \eqref{eq:HBE}\textendash \eqref{eq:HBEinit}
changes sign. 

First we discuss the existence of global in time solution of equations
\eqref{eq:HBE}\textendash \eqref{eq:HBEinit} and of the more general
equation 
\begin{equation}
\phi_{tt}-e^{-2t}\Delta\phi+n\phi_{t}=\mu^{2}\phi-F\left(\phi\right),\label{eq:HBEF}
\end{equation}
where $F\left(\phi\right)=\pm\phi^{p}$ or $F\left(\phi\right)=\pm|\phi|^{p-1}\phi$,
or some more general function, and $\vec{x}\in\mathbb{R}^{n}$. For
this equation the local existence of solution in different spaces
of functions is well investigated for appropriate values of $p$ and
$n$. The estimate for the lifespan of the solution is given as follows.
The main parameter that controls estimates and solvability is the
principal square root $M:=(\mu^{2}+n^{2}/4)^{\frac{1}{2}}$ . The
following result is from Theorem 0.1 in \cite{KY-ArXiv2017}. 

Let $H^{s}\left(\mathbb{R}^{n}\right)=W^{s,2}\left(\mathbb{R}^{n}\right)$
be the usual Sobolev space \cite{Adams}. The function $F$ is said
to be Lipschitz continuous with exponent $\alpha\geq0$ in the space
$H^{s}\left(\mathbb{R}^{n}\right)$ if there is a constant $C\geq0$
such that 
\begin{equation}
\left\Vert F\left(\phi_{1}\left(x\right)\right)-F\left(\phi_{2}\left(x\right)\right)\right\Vert {}_{H^{s}\left(\mathbb{R}^{n}\right)}\leq C\left\Vert \phi_{1}-\phi_{2}\right\Vert {}_{H^{s}\left(\mathbb{R}^{n}\right)}\left(\left\Vert \phi_{1}\right\Vert {}_{H^{s}\left(\mathbb{R}^{n}\right)}^{\alpha}+\left\Vert \phi_{2}\right\Vert {}_{H^{s}\left(\mathbb{R}^{n}\right)}^{\alpha}\right)\label{eq:calM}
\end{equation}
for all $\phi_{1},\phi_{2}\in H^{s}\left(\mathbb{R}^{n}\right)$.
Assume that the nonlinear term $F\left(u\right)$ is Lipschitz continuous
in the space $H^{s}\left(\mathbb{R}^{n}\right)$, $s>n/2\geq1$, $F\left(0\right)=0$,
$\alpha>0$ and $M\in\mathbb{C}$. According to $\mbox{{\rm (iii)}}$
of Theorem 0.1 from \cite{KY-ArXiv2017} we have the following statement:
\textit{If $\mu>0$, then the lifespan $T_{ls}$ of the solution can
be estimated from below as follows: 
\begin{eqnarray*}
T_{ls} & \geq & -\frac{1}{M-\frac{n}{2}}\ln\left(\left\Vert \varphi_{0}\right\Vert {}_{H^{s}\left(\mathbb{R}^{n}\right)}+\left\Vert \varphi_{1}\right\Vert {}_{H^{s}\left(\mathbb{R}^{n}\right)}\right)-C\left(m,n,\alpha\right)
\end{eqnarray*}
with some constant $C\left(m,n,\alpha\right)$ for sufficiently small
$\left\Vert \varphi_{0}\right\Vert {}_{H^{s}\left(\mathbb{R}^{n}\right)}$
and $\left\Vert \varphi_{1}\right\Vert _{H^{s}\left(\mathbb{R}^{n}\right)}$.
} In particular, this covers the cases 
\[
F\left(\phi\right)=\pm\left|\phi\right|^{\alpha}\phi\quad\mbox{{\rm and}}\quad F\left(\phi\right)=\pm\left|\phi\right|^{\alpha+1}\quad\mbox{{\rm and}}\quad F\left(\phi\right)=\lambda\phi^{3}.
\]

For the function $u=e^{\frac{n}{2}t}\phi$, the equation \eqref{eq:HBEF}
leads to 
\begin{equation}
u_{tt}-e^{-2t}\bigtriangleup u-M^{2}u=-e^{\frac{n}{2}t}F\left(e^{-\frac{n}{2}t}u\right),\label{eq:K_G_Higgs}
\end{equation}
with $M\geq n/2$. The last equation can be regarded as the Klein-Gordon
equation whose squared physical mass $m^{2}$ is negative, $m^{2}=-M^{2}<0$.
The quantum fields whose squared physical masses are negative (imaginary
mass) represent tachyons (see, e.g., \cite{B-F-K-L}). According to
\cite{B-F-K-L} the free tachyons in the Minkowski spacetime have
to be rejected on stability grounds since the localized disturbances
of the Klein-Gordon equation with imaginary mass spread with at most
the speed of light, but grow exponentially. 

Epstein and Moschella in \cite{Epstein-Moschella} give a complete
study of a family of scalar tachyonic quantum fields which are linear
Klein-Gordon quantum fields on the de Sitter manifold whose squared
masses are negative 
\begin{equation}
\phi_{tt}-e^{-2t}\Delta\phi+n\phi_{t}+m^{2}\phi=0\,,
\end{equation}
and take an infinite set of discrete values $m^{2}=-k(k+n)$, $k=0,1,2,\ldots$.
The nonexistence of a global in time solution of the semilinear Klein-Gordon
massive tachyonic (self-interacting quantum fields) equation in the
de~Sitter spacetime, that is a finite lifespan, is proved in \cite{yagdjian_DCDS}.
In fact, Theorem 1.1 in \cite{yagdjian_DCDS} states that if $c\not=0$,
$\alpha>0$, and $m\not=0$, then for every positive numbers $\varepsilon$
and $s$ there exist functions $\varphi_{0}$, $\,\varphi_{1}\in C_{0}^{\infty}\left(\mathbb{R}^{n}\right)$
such that $\left\Vert \varphi_{0}\right\Vert _{H^{s}\left(\mathbb{R}^{n}\right)}+\left\Vert \varphi_{1}\right\Vert _{H^{s}\left(\mathbb{R}^{n}\right)}\leq\varepsilon$
but the solution $\phi=\phi\left(x,t\right)$ of the semilinear equation
\begin{equation}
\phi_{tt}-e^{-2t}\Delta\phi+n\phi_{t}-m^{2}\phi=c\left|\phi\right|^{1+\alpha}\,,
\end{equation}
with the initial values $\phi\left(\vec{x},0\right)=\varphi_{0}\left(\vec{x}\right)$,
$\phi_{t}\left(\vec{x},0\right)=\varphi_{1}\left(\vec{x}\right)$,
blows up in finite time. This would also imply the blowup of the \textit{sign-preserving}
solutions (under some additional conditions) of the equation
\begin{equation}
\phi_{tt}-e^{-2t}\Delta\phi+n\phi_{t}-m^{2}\phi=-|\phi|^{\alpha}\phi\,.
\end{equation}
Thus the issue of the existence of a global in time of solution of
equation \eqref{eq:HBE} is still an open problem. 

There have been numerous numerical approaches for solving the various
types of nonlinear Klein-Gordon equations and other nonlinear wave
equations. Most of the numerical results are for one space dimension,
including a cubic B-spline collocation method presented in \cite{RASHIDINIA20101866};
the Adomian decomposition method for solitary waves in \cite{KAYA2004341};
and the method of lines in one space dimensions used in \cite{griffiths2010traveling}
for the mth-order Klein\textendash Gordon equation. For computational
analysis of nonlinear hyperbolic equations it is important to preserve
not just the dissipation of energy (see \cite{MACIASDIAZ2010552})
but to minimize the dispersion error as well (see \cite{BOGEY2004194}).
Explicit methods, even though they are only conditionally stable,
have a tendency to have smaller dispersion error compared to the usually
unconditionally stable implicit methods. Comparison of several explicit
finite difference methods have been presented for the one spatial
variable case in \cite{JIMENEZ199061}. A differential transform method
with variational iteration method with Adomian\textquoteright s polynomials
was presented very recently for the Higgs boson equation in de Sitter
Spacetime in \cite{Yazici}. The Adomian decomposition method was
also used for the Klein-Gordon equation with quadratic nonlinearity
in \cite{BASAK2009718} with a general method presented for the three
space dimensional case and an example with a known solution presented
for a one space dimensional case. Even in higher spatial dimensions
if one uses radial basis functions the resulting problem becomes one-dimensional
(see \cite{DEHGHAN2009400} and \cite{Donninger}). With a simple
unit cube for the computation domain we chose a ``quick and dirty''
finite difference discretization with fourth order of accuracy for
the spatial component along with a matching fourth-order Runge-Kutta
method for the time variable. The resulting explicit numerical code
(stencil code) is very well suited for high performance computations
using Graphical Processing Units. In Section \ref{sec:The-Numerical-Method}
we describe the numerical approach for the $\left(3+1\right)$\textendash dimensional
general case and for $\left(1+1\right)$\textendash dimensional radial
solutions. In Section \ref{sec:Computational-Examples} the main computational
results are presented via various computational examples in order
to test the numerical code as well as to examine the properties of
the Higgs boson equation. Finally in Section \ref{sec:Concluding-Remarks}
conclusions are given in the form of conjectures.

\section{\label{sec:The-Numerical-Method}The Numerical Method}

Our numerical approach uses a fourth order finite difference method
in the three-dimensional space in combination with an explicit fourth
order Runge-Kutta method in time for the discretization and numerical
solution of the Higgs boson equation. In addition to the general case
of three spatial dimensions we also investigate radial solutions in
one spatial dimension, which is much less demanding computationally.
In this section we describe these choices of approach. 

\subsection{General 3D Solutions}

We only consider solutions of the Higgs boson equation \eqref{eq:HBE}
with compact support in space. Since the solution has the finite speed
of propagation ${\displaystyle e^{-t}}$ (see, e.g., \cite{Hormander_1997}),
the total distance travelled by the solution is $\int_{0}^{\infty}e^{-t}dt=1.$
The finite cone of influence along with a rescaling of the spatial
domain enables us to use zero boundary conditions on the unit box
$\Omega=\left[0,1\right]\times\left[0,1\right]\times\left[0,1\right]$
as computational domain. The rescaled Higgs boson equation then has
the form
\begin{align}
\phi_{tt}-\frac{1}{L^{2}}e^{-2t}\Delta\phi+3\phi_{t} & =\mu^{2}\phi-\lambda\phi{}^{3},\quad\vec{x}\in\Omega,\quad t\in\left(0,T\right],\label{eq:higgs-boson-rescalled}\\
\phi\left(\vec{x},0\right) & =\varphi_{0}\left(\vec{x}\right),\quad\vec{x}\in\Omega,\\
\phi_{t}\left(\vec{x},0\right) & =\varphi_{1}\left(\vec{x}\right),\quad\vec{x}\in\Omega,\\
\phi\left(\partial\Omega,t\right) & =0,\quad t\in\left[0,T\right].\label{eq:BC}
\end{align}
Here we set the initial condition to have compact support inside $\Omega$
and choose $L>0$ large enough so that information propagating from
the initial compact support will never reach the boundary of the unit
box. Since the solution will be zero on continuous regions outside
the compact support, we use finite difference method for spatial discretization,
which is a local approach as opposed to global approaches like spectral
methods. A uniform grid is chosen in all three space variables with
grid spacing $\delta x=\delta y=\delta z=1/n$ with $n\in\mathbb{N}$
and with notation 

\begin{eqnarray*}
\phi_{jkl}\left(t\right) & = & \phi\left(x_{j},y_{k},z_{l},t\right)=\phi\left(j\delta x,\,k\delta y,\,l\delta z,\,t\right)
\end{eqnarray*}
for $j,k,l=0\ldots n$ and for $t\in\left[0,T\right]$. The second
partial derivatives in the Laplacian ${\displaystyle \Delta\phi=\frac{\partial^{2}\phi}{\partial x^{2}}+\frac{\partial^{2}\phi}{\partial y^{2}}+\frac{\partial^{2}\phi}{\partial z^{2}}}$
are discretized using the fourth-order central difference scheme (see,
e.g., \cite{finitediff}){\small{}
\begin{align*}
\left.\frac{\partial^{2}\phi}{\partial x^{2}}\right|_{jkl} & =\frac{1}{\delta x^{2}}\left(-\frac{1}{12}\phi_{j-2,kl}+\frac{4}{3}\phi_{j-1,kl}-\frac{5}{2}\phi_{jkl}+\frac{4}{3}\phi_{j+1,kl}-\frac{1}{12}\phi_{j+2,kl}\right)+O\left(dx^{4}\right),\\
\left.\frac{\partial^{2}\phi}{\partial y^{2}}\right|_{jkl}^{i} & =\frac{1}{\delta y^{2}}\left(-\frac{1}{12}\phi_{jk-2,l}+\frac{4}{3}\phi_{jk-1,l}-\frac{5}{2}\phi_{jkl}+\frac{4}{3}\phi_{jk+1,l}-\frac{1}{12}\phi_{jk+2,l}\right)+O\left(dy^{4}\right),\\
\left.\frac{\partial^{2}\phi}{\partial z^{2}}\right|_{jkl}^{i} & =\frac{1}{\delta z^{2}}\left(-\frac{1}{12}\phi_{jkl-2}+\frac{4}{3}\phi_{jkl-1}-\frac{5}{2}\phi_{jkl}+\frac{4}{3}\phi_{jkl+1}-\frac{1}{12}\phi_{jkl+2}\right)+O\left(dz^{4}\right).
\end{align*}
}We also transform equation \eqref{eq:higgs-boson-rescalled} into
a first-order-in-time system of two equations via notations
\begin{align}
\vec{v}\left(t\right) & =\left(\begin{array}{l}
\vec{v}_{1}\\
\vec{v}_{2}
\end{array}\right)\left(t\right)=\left(\begin{array}{l}
\vec{\phi}_{jkl}\\
\frac{\partial}{\partial t}\vec{\phi}_{jkl}
\end{array}\right)\left(t\right)
\end{align}
and
\begin{align}
\vec{f}\left(t,\vec{v}\left(t\right)\right) & =\left(\begin{array}{c}
\vec{v}_{2}\\
-3\vec{v}_{2}+\frac{e^{-2t}}{L^{2}}\Delta\vec{v}_{1}+\mu^{2}\vec{v}_{1}-\lambda\vec{v}_{1}^{3}
\end{array}\right),\label{eq:rhs}
\end{align}
where $\vec{v}_{1}=\vec{\phi}_{jkl}=\left(\phi_{000},\ldots\phi_{nnn}\right)^{T}$.
This way we obtain an evolution system in the form
\begin{align}
\vec{v}'\left(t\right) & =\vec{f}\left(t,\vec{v}\left(t\right)\right)
\end{align}
with initial and boundary conditions
\begin{align*}
\vec{v}\left(\vec{x},0\right) & =\left(\begin{array}{l}
\vec{\varphi}_{0_{jkl}}\\
\vec{\varphi}_{1_{jkl}}
\end{array}\right)\left(\vec{x}\right),\quad\vec{x}\in\Omega,\\
\vec{v}\left(\partial\Omega,t\right) & =\vec{0},\quad t\in\left[0,T\right].
\end{align*}
Note that the second component of the right-hand-side function \eqref{eq:rhs}
has the form
\begin{multline*}
\left(f_{2}\right)_{jkl}=\mu^{2}\left(v_{1}\right)_{jkl}-\lambda\left(v_{1}\right)_{jkl}^{3}-3\left(v_{2}\right)_{jkl}+\frac{e^{-2t}}{L^{2}}\Delta\left(v_{1}\right)_{jkl}\\
=\left(\mu^{2}-\frac{15}{2}\frac{e^{-2t}}{L^{2}\delta x^{2}}\right)\left(v_{1}\right)_{jkl}-\lambda\left(v_{1}\right)_{jkl}^{3}-3\left(v_{2}\right)_{jkl}\\
-\frac{1}{12}\frac{e^{-2t}}{L^{2}\delta x^{2}}\left(\left(v_{1}\right)_{jkl+2}+\left(v_{1}\right)_{jkl-2}+\left(v_{1}\right)_{jkl+2n}+\left(v_{1}\right)_{jkl-2n}+\left(v_{1}\right)_{jkl+2n^{2}}+\left(v_{1}\right)_{jkl-2n^{2}}\right)\\
+\frac{4}{3}\frac{e^{-2t}}{L^{2}\delta x^{2}}\left(\left(v_{1}\right)_{jkl+1}+\left(v_{1}\right)_{jkl-1}+\left(v_{1}\right)_{jkl+n}+\left(v_{1}\right)_{jkl-n}+\left(v_{1}\right)_{jkl+n^{2}}+\left(v_{1}\right)_{jkl-n^{2}}\right).
\end{multline*}
For time discretization we use the classical, explicit fourth-order
Runge-Kutta method (see, e.g., \cite{RKbook})

\begin{align*}
\vec{v}\left(t+\delta t\right) & =\vec{v}\left(t\right)+\frac{\vec{k}_{1}+2\vec{k}_{2}+2\vec{k}_{3}+\vec{k}_{4}}{6}\delta t,
\end{align*}
where
\begin{align*}
\vec{k_{1}} & =\vec{f}\left(t,\vec{v}\left(t\right)\right),\\
\vec{k}_{2} & =\vec{f}\left(t+\frac{\delta t}{2},\vec{v}\left(t\right)+\vec{k}_{1}\frac{\delta t}{2}\right),\\
\vec{k}_{3} & =\vec{f}\left(t+\frac{\delta t}{2},\vec{v}\left(t\right)+\vec{k}_{2}\frac{\delta t}{2}\right),\\
\vec{k}_{4} & =\vec{f}\left(t+\delta t,\vec{v}\left(t\right)+\vec{k}_{3}\delta t\right).
\end{align*}
This conditionally stable explicit numerical scheme enables us to
use Graphical Processing Units for high-performance computing. We
also reused variables for decreasing memory storage by setting
\begin{align*}
\vec{k_{1}} & =\vec{f}\left(t,\vec{v}\left(t\right)\right),\\
\vec{k}_{2} & =\vec{f}\left(t+\frac{\delta t}{2},\vec{v}\left(t\right)+\vec{k}_{1}\frac{\delta t}{2}\right),\\
\vec{k_{1}} & =\vec{k}_{1}+2\vec{k}_{2},\\
\vec{k}_{3} & =\vec{f}\left(t+\frac{\delta t}{2},\vec{v}\left(t\right)+\vec{k}_{2}\frac{\delta t}{2}\right),\\
\vec{k}_{2} & =\vec{f}\left(t+\delta t,\vec{v}\left(t\right)+\vec{k}_{3}\delta t\right),\\
\vec{v}\left(t+\delta t\right) & =\vec{v}\left(t\right)+\frac{\vec{k}_{1}+2\vec{k}_{3}+\vec{k}_{2}}{6}\delta t.
\end{align*}
The numerical code has been programmed using PGI CUDA Fortran Compiler
\cite{PGI} and was performed on NVIDIA Tesla K40c GPU Accelerators.
Using texture data/memory (see, e.g., \cite{CUDA-Fortran-Guide})
we were able to speed up computations by more than $20\%$. The longest
calculation took 27 hours to reach the non-dimensional time $t=50$
in the simulation. The visualization was done using software packages
ParaView \cite{Ayachit:2015:PGP:2789330} and MATLAB \cite{MATLAB}. 

\subsection{Radial Solutions}

Some of our numerical examples are for radial solutions of the Higgs
boson equation. This simplification is justified partly according
to the cosmological principle (see, e.g., \cite{Liddle201507}) that
the universe is homogeneous and isotropic on large scales. The radial
form of equation \eqref{eq:higgs-boson-rescalled} is
\begin{equation}
\phi_{tt}-\frac{e^{-2t}}{L^{2}}\left(\frac{2}{r}\frac{\partial}{\partial r}\phi+\frac{\partial^{2}}{\partial r^{2}}\phi\right)+3\phi_{t}=\mu^{2}\phi-\lambda\phi^{3},\quad r\in\left(0,1\right),\quad t\in\left(0,T\right]\label{eq:radial}
\end{equation}
with boundary conditions 
\[
\phi_{r}\left(0,t\right)=\phi\left(1,t\right)=0,\quad t\in\left[0,T\right]
\]
and with initial condition
\begin{align}
\phi\left(r,0\right) & =\varphi_{0}\left(r\right),\quad r\in\left[0,1\right],\\
\phi_{t}\left(r,0\right) & =\varphi_{1}\left(r\right),\quad r\in\left[0,1\right].
\end{align}
The singularity in equation \eqref{eq:radial} at $r=0$ is treated
using L'Hospital's rule
\begin{equation}
\lim_{r\to0}\frac{\frac{\partial}{\partial r}\phi\left(r,t\right)}{r}=\lim_{r\to0}\frac{\partial^{2}\phi\left(r,t\right)}{\partial r^{2}}=\frac{\partial^{2}\phi\left(0,t\right)}{\partial r^{2}}.
\end{equation}
At the boundary $r=0$ symmetry is used for boundary conditions using
grid points only from the $\left[0,1\right]$ spatial domain:
\begin{align*}
\frac{\partial^{2}\phi\left(0,t\right)}{\partial r^{2}} & \approx\frac{1}{\delta r^{2}}\left(-\frac{1}{12}\phi\left(-2\delta r\right)+\frac{4}{3}\phi\left(-\delta r\right)-\frac{5}{2}\phi\left(0\right)+\frac{4}{3}\phi\left(\delta r\right)-\frac{1}{12}\phi\left(2\delta r\right)\right)\\
 & =\frac{1}{\delta r^{2}}\left(-\frac{5}{2}\phi\left(0\right)+\frac{8}{3}\phi\left(\delta r\right)-\frac{1}{6}\phi\left(2\delta r\right)\right).
\end{align*}
Simulations for the radial equation \eqref{eq:radial} and for the
full three-space-dimensional equation \eqref{eq:higgs-boson-rescalled}
produced equivalent results. Below we provide all results obtained
from numerical simulations for the full three-space-dimensional case
both for radial and for general solutions.

\section{\label{sec:Computational-Examples}Computational Examples}

\begin{figure}[H]
\caption{\label{fig:Line-segments}Line segments used for line plots in the
computational domain.}
\centering{}%
\begin{tabular}{ccc}
(a) Line segment parallel to the $x$-axis &  & (b) Diagonal line segment\tabularnewline
\includegraphics[width=0.5\columnwidth]{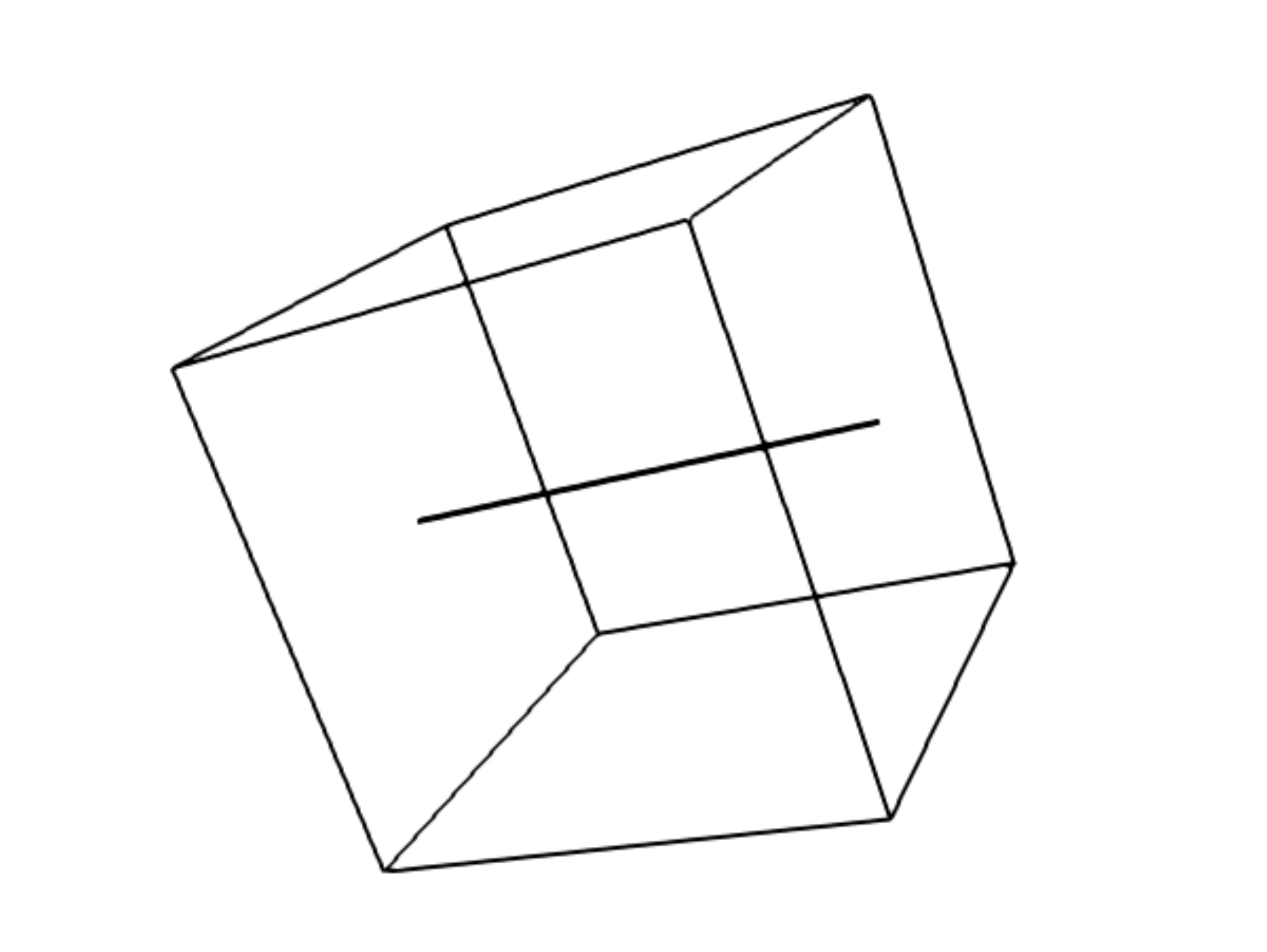} &  & \includegraphics[width=0.5\columnwidth]{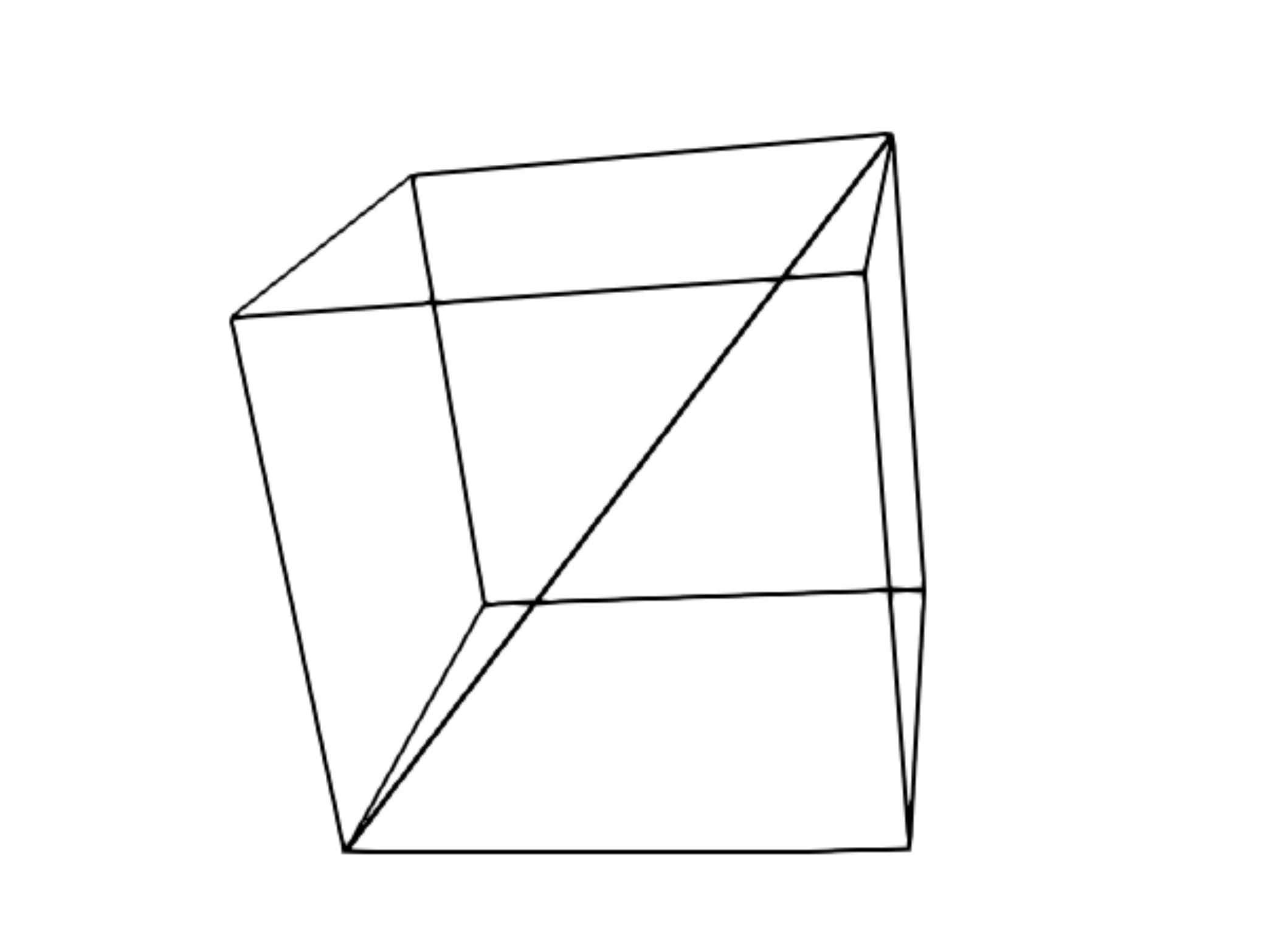}\tabularnewline
\end{tabular}
\end{figure}

\begin{figure}[H]
\caption{\label{fig:Bump-Function}Line plot of a bump function.}
\centering{}\includegraphics[width=0.5\columnwidth]{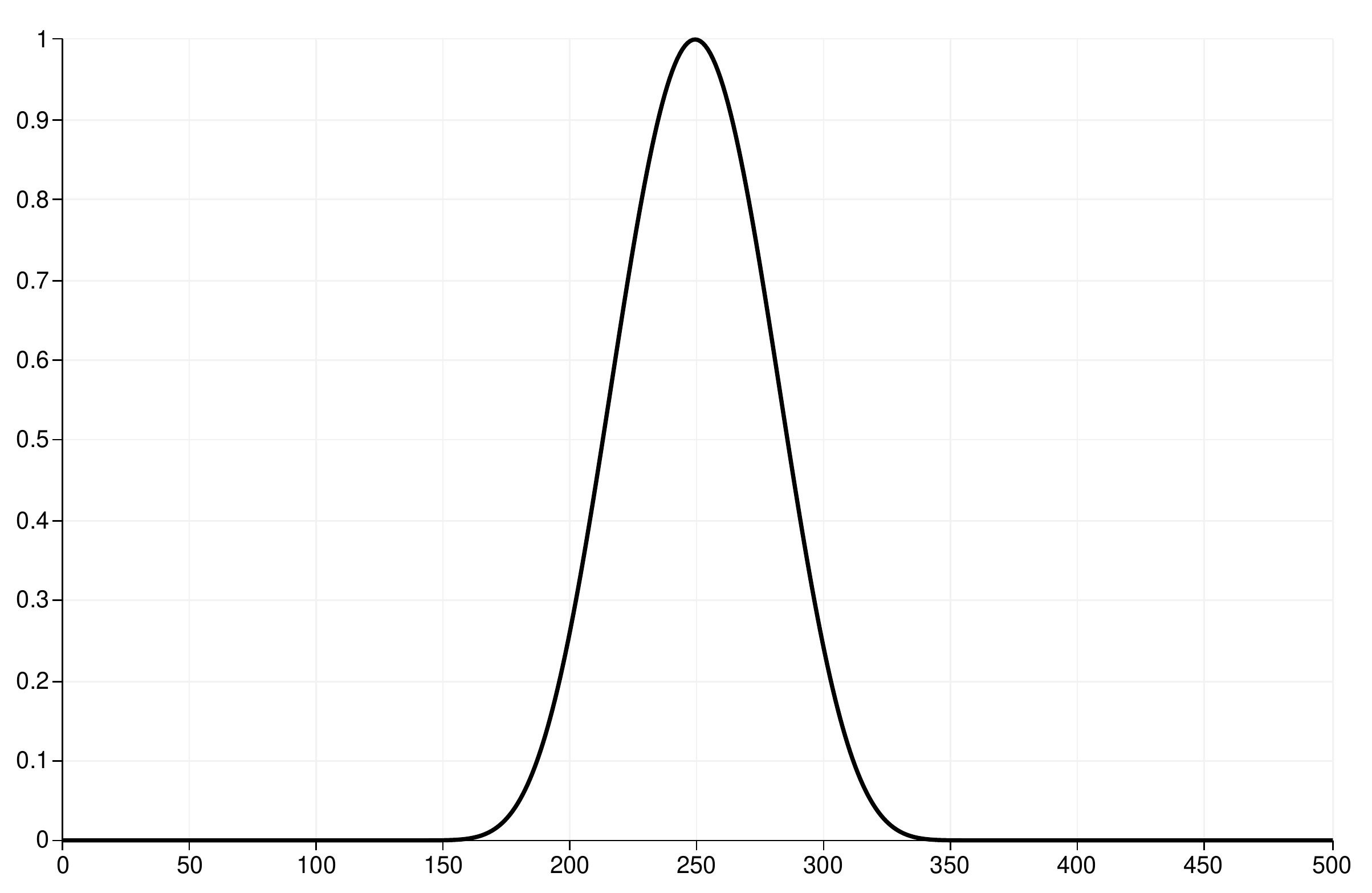}
\end{figure}

\begin{figure}[H]
\caption{\label{fig:Convergence}Convergence of the numerical method.}
\centering{}%
\begin{tabular}{cc}
\multicolumn{2}{c}{Difference between numerical solutions along a line. $w_{n}$ is the
numerical solution for grid size $n$. }\tabularnewline
$t=1$ & $t=2$\tabularnewline
\includegraphics[width=0.5\textwidth]{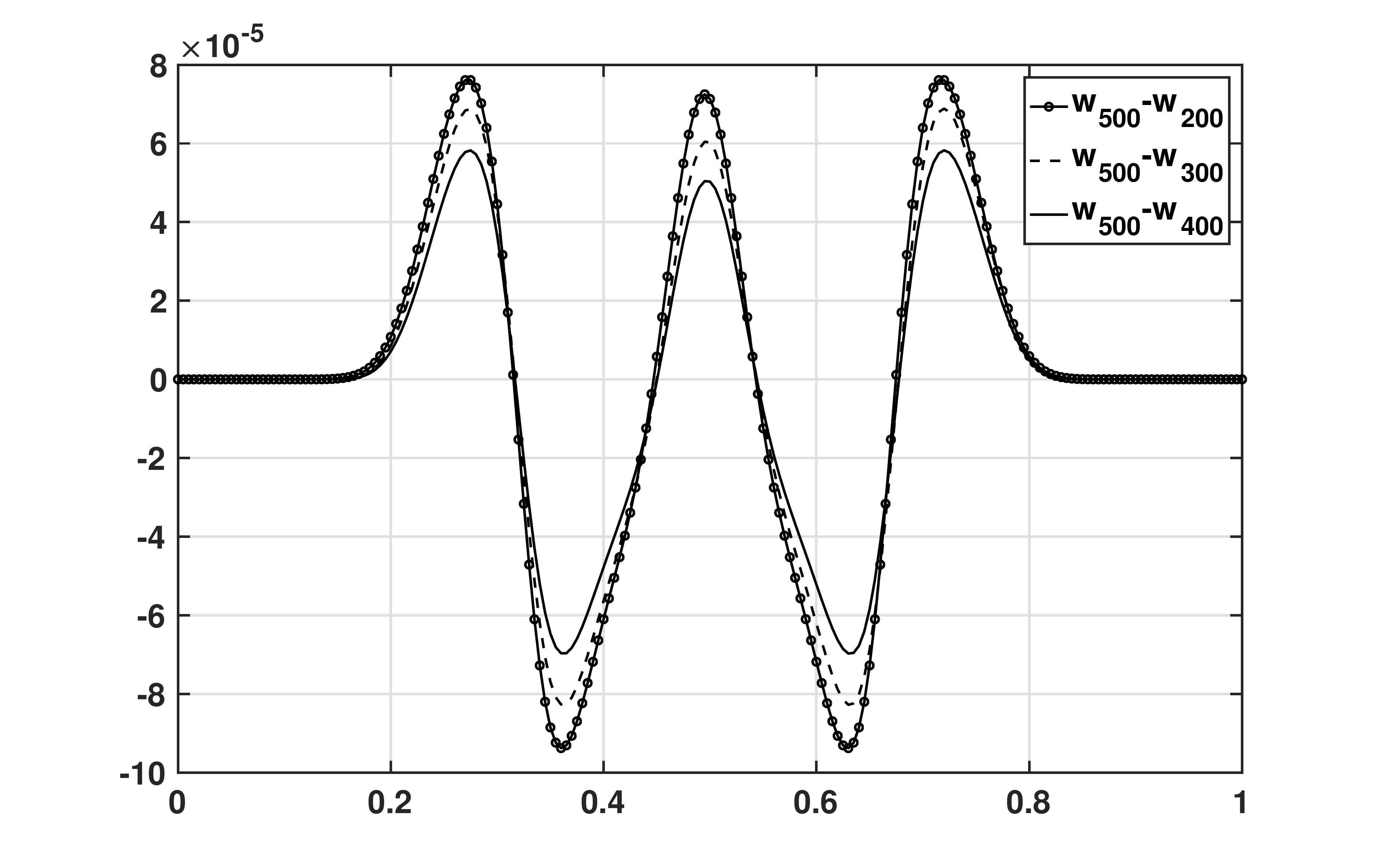} & \includegraphics[width=0.5\textwidth]{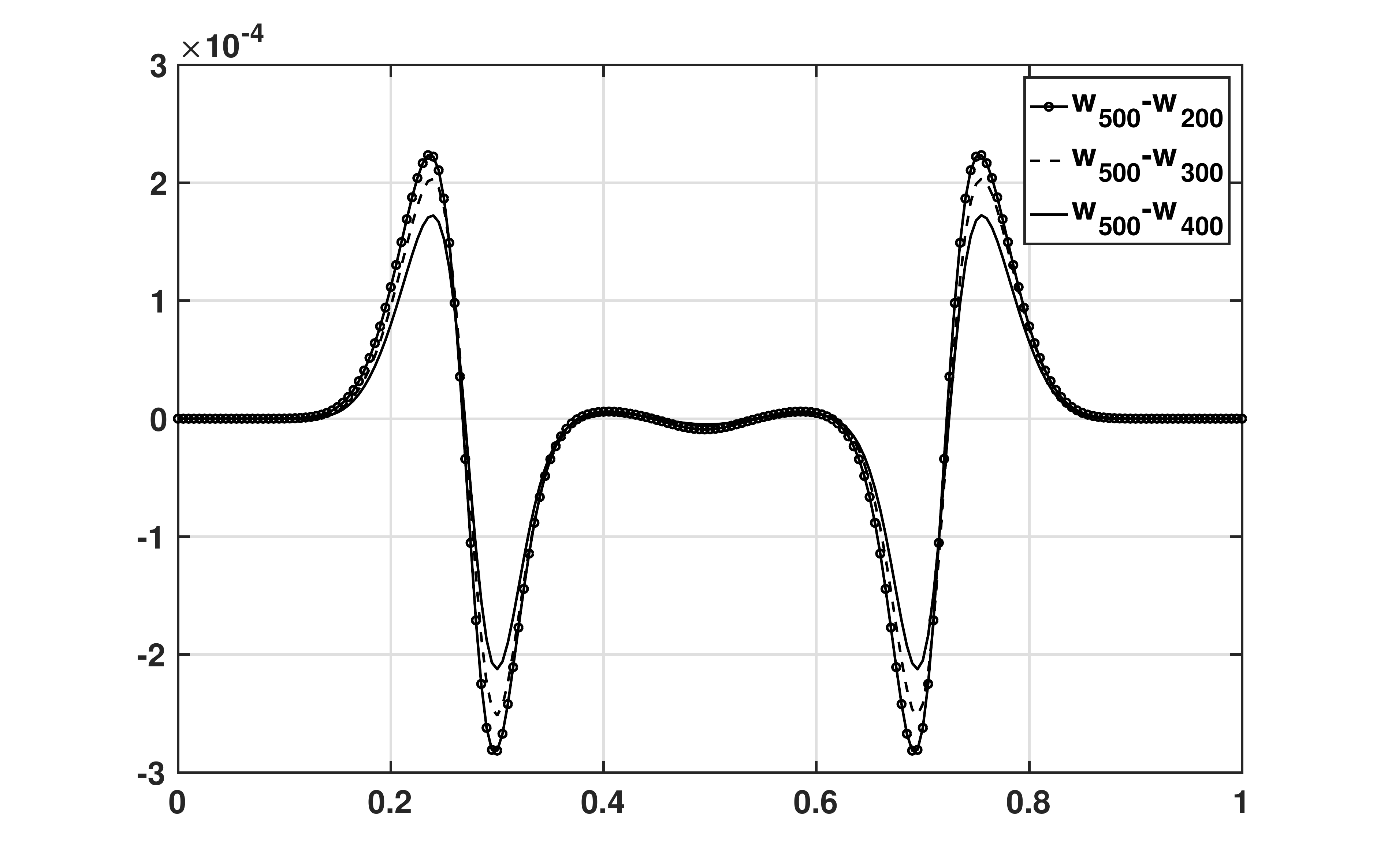}\tabularnewline
 & \tabularnewline
\multicolumn{2}{c}{Numerical solution $w_{500}$ along a line }\tabularnewline
$t=1$ & $t=2$\tabularnewline
\includegraphics[width=0.45\textwidth]{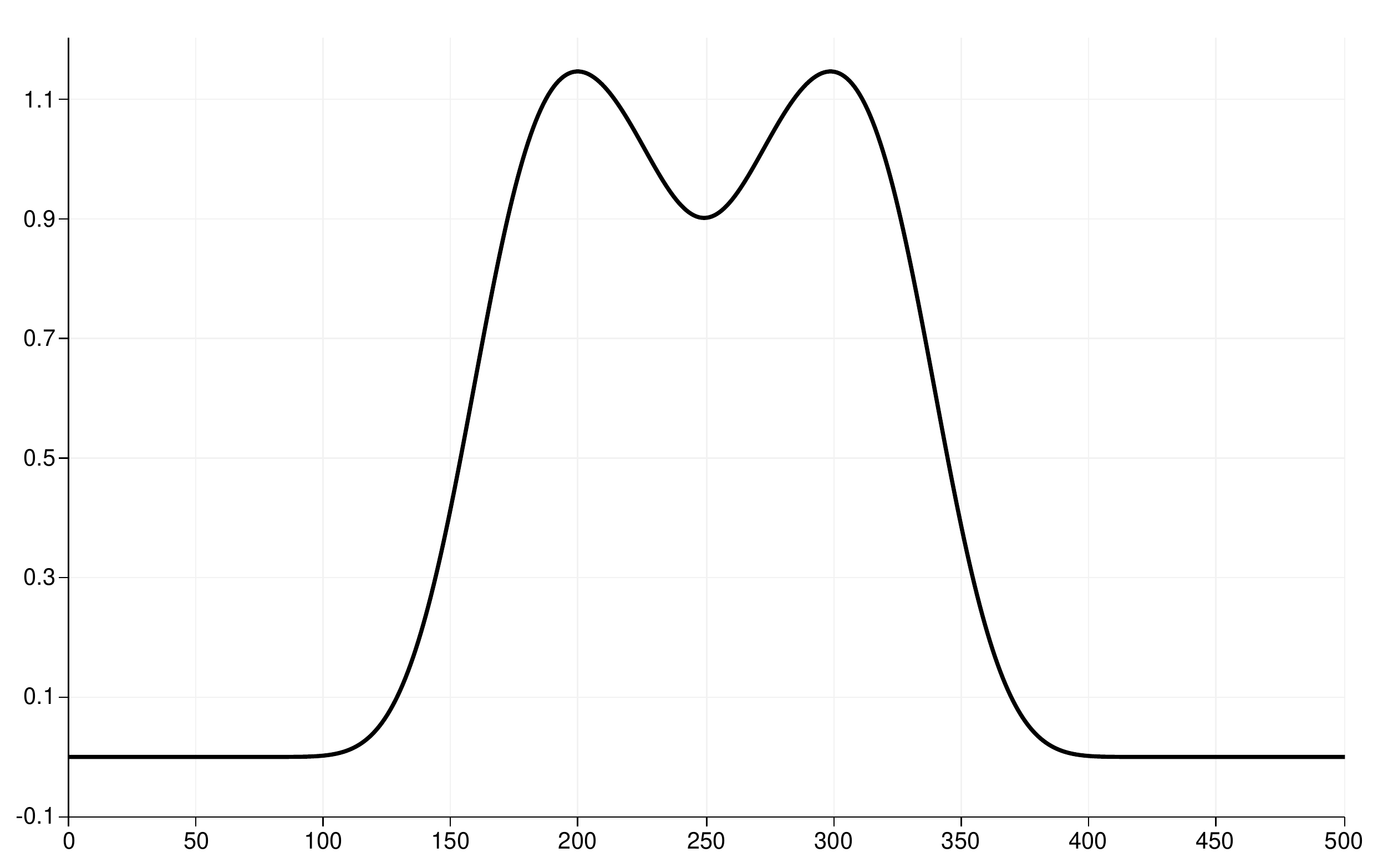} & \includegraphics[width=0.45\textwidth]{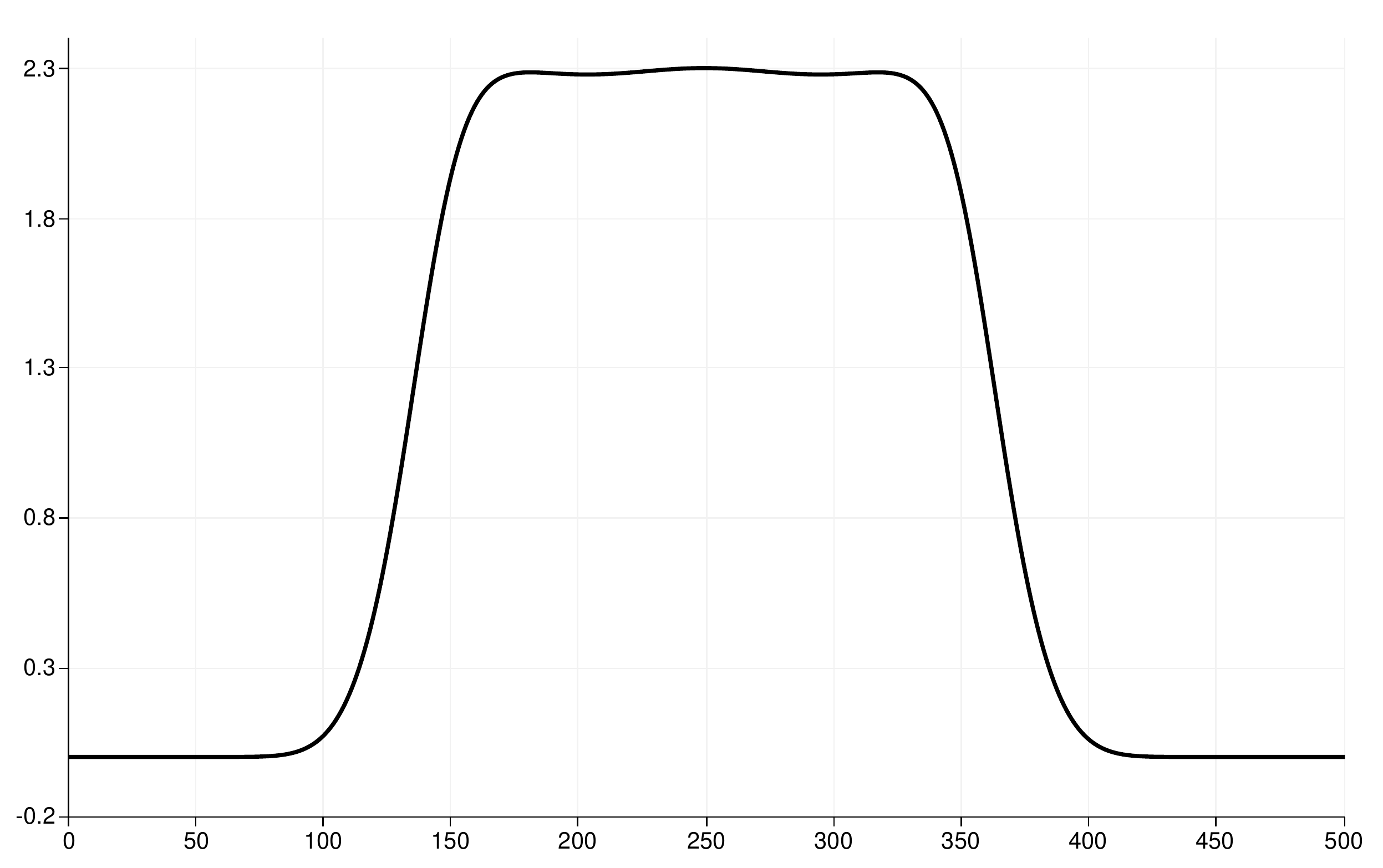}\tabularnewline
\end{tabular}
\end{figure}

\begin{figure}[H]
\caption{\label{fig:Precision}Difference between double and single precision.}
\centering{}%
\begin{tabular}{cc}
$t=1$ & $t=2$\tabularnewline
\includegraphics[width=0.5\columnwidth]{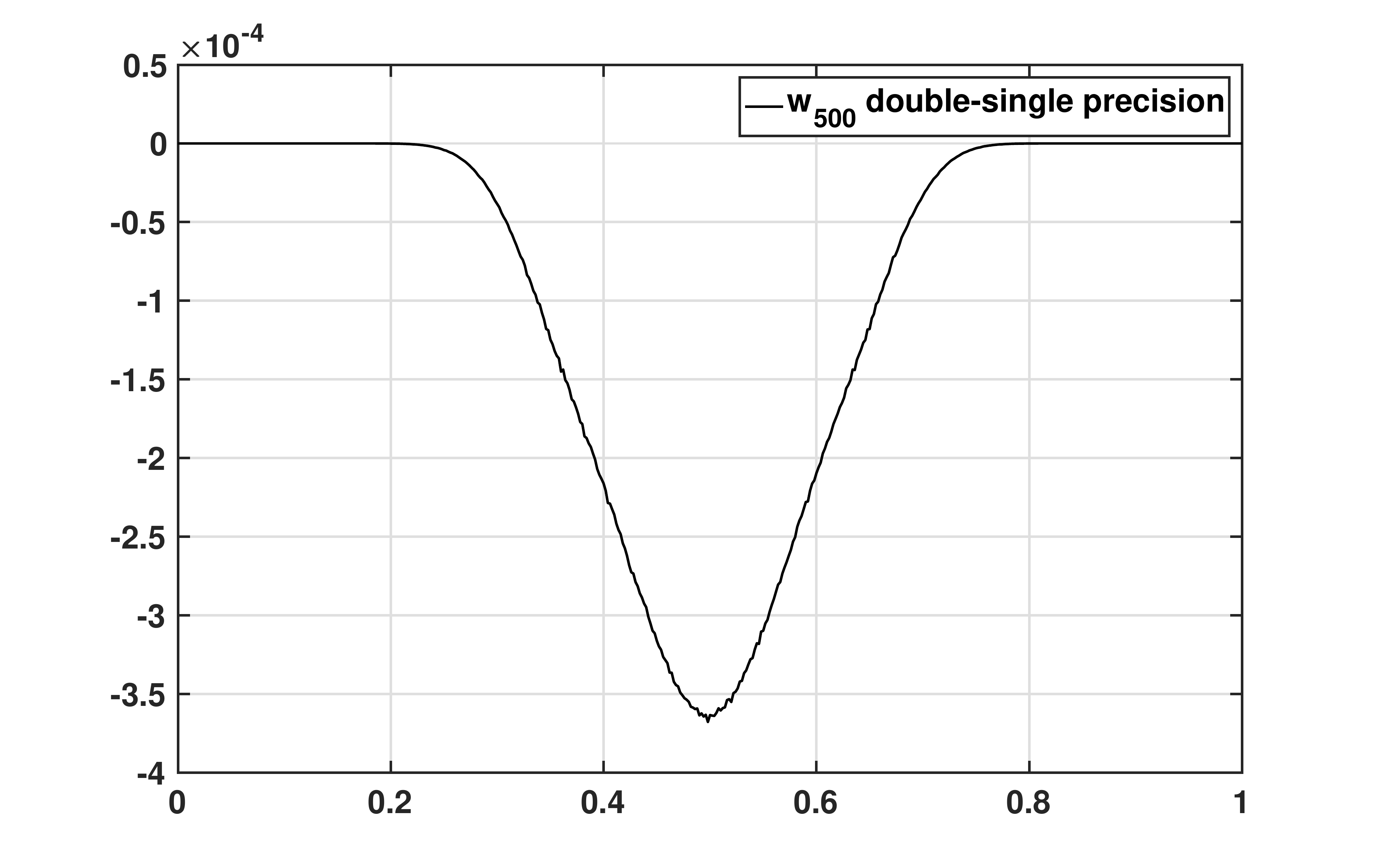} & \includegraphics[width=0.5\columnwidth]{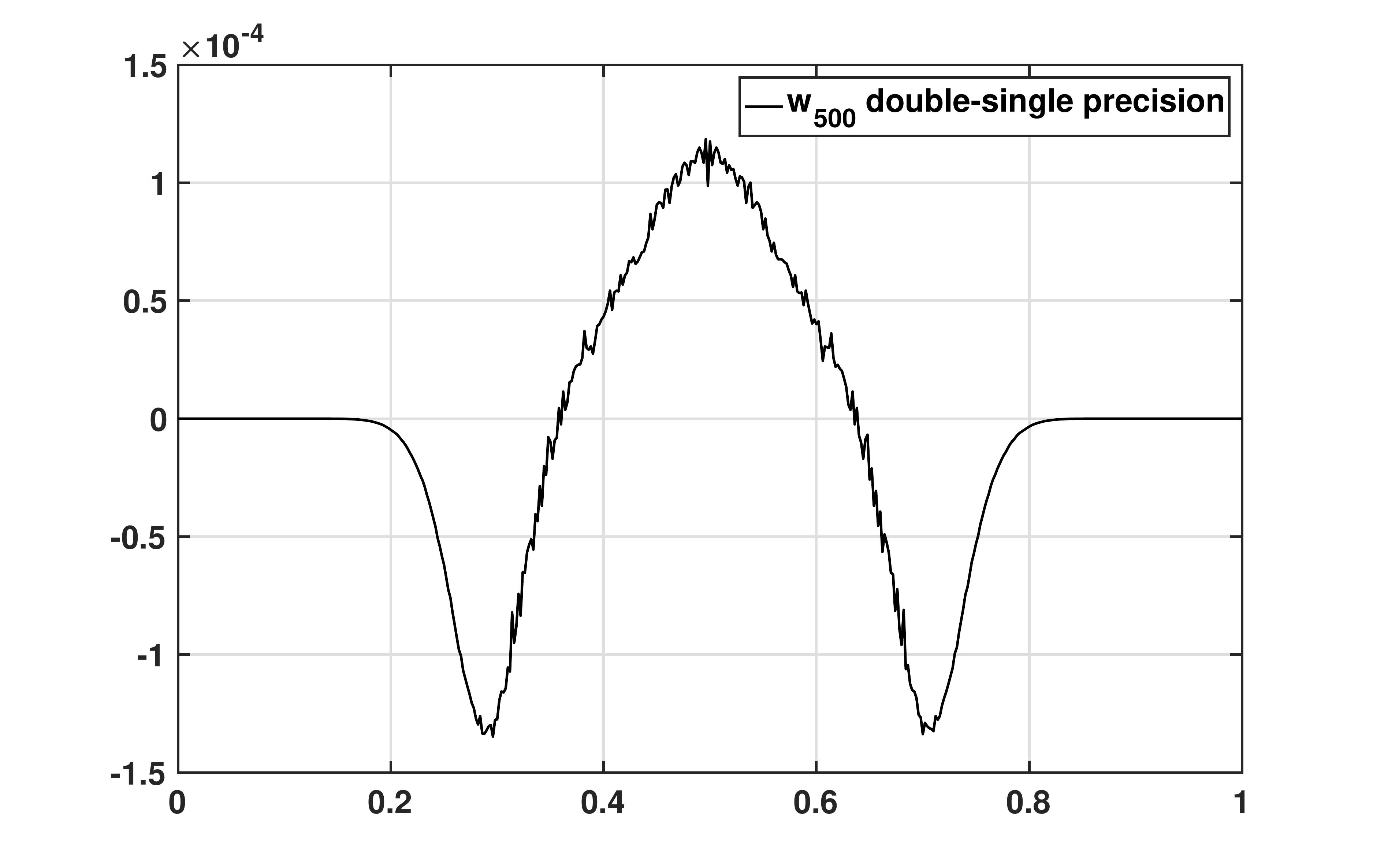}\tabularnewline
\end{tabular}
\end{figure}

\begin{figure}[H]
\caption{\label{fig:Blowup}Blow up of solution to the Klein-Gordon equation
with an imaginary mass. }
\centering{}%
\begin{tabular}{cc}
(a) Line plot of solution at $t=0$  & (b) Line plot of solution at $t=2.0$\tabularnewline
\includegraphics[width=0.48\columnwidth]{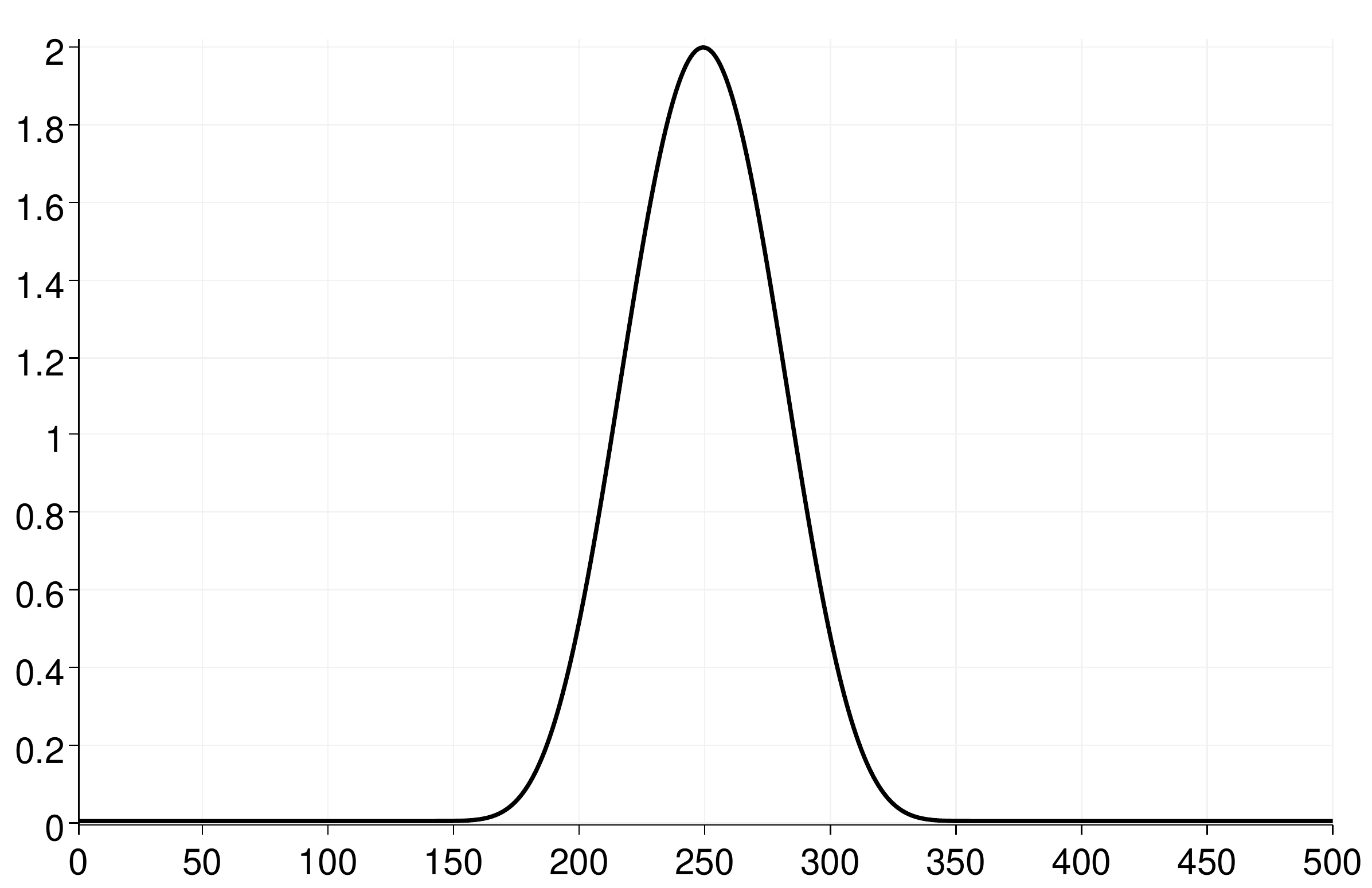} & \includegraphics[width=0.48\columnwidth]{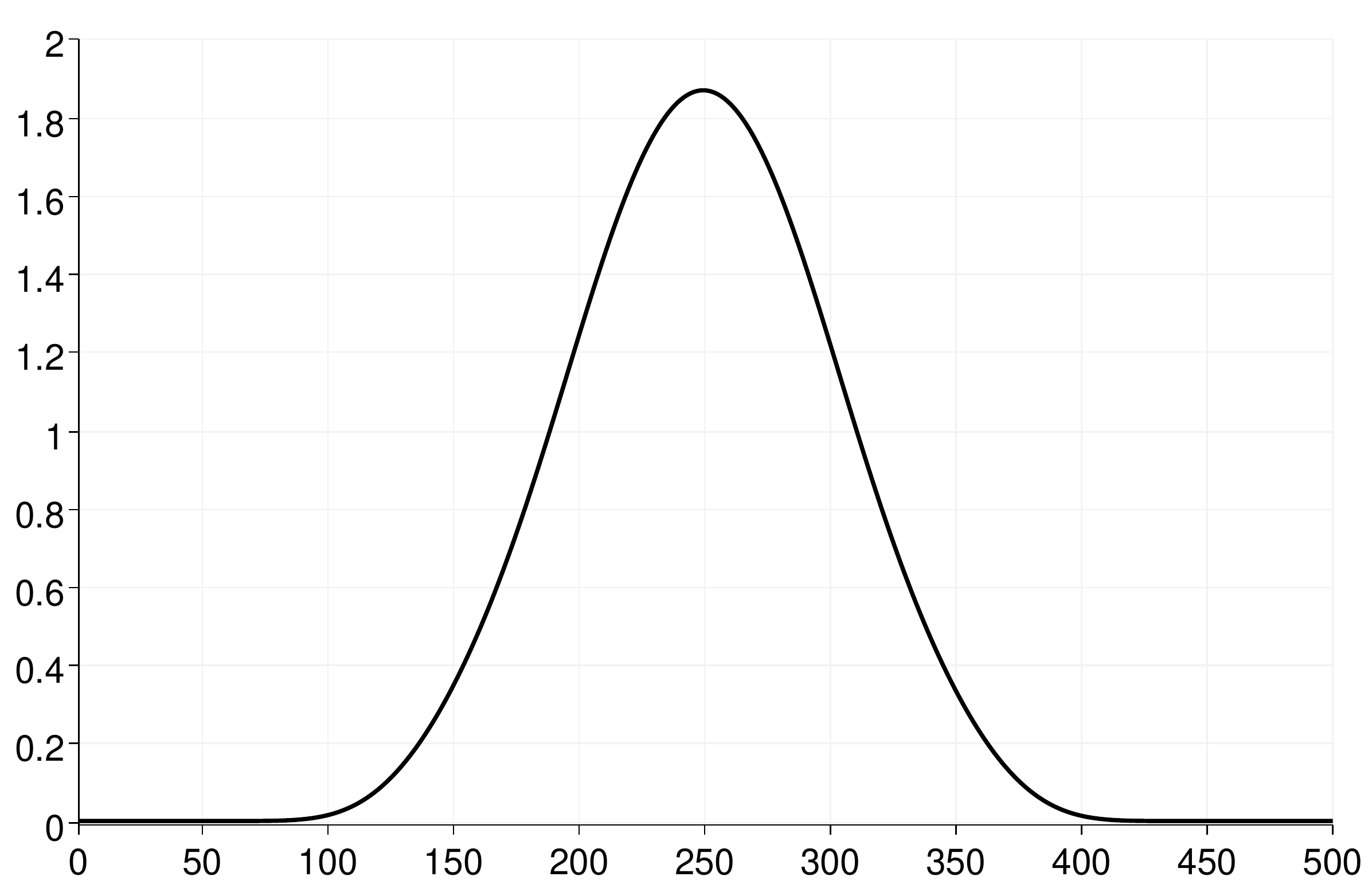}\tabularnewline
 & \tabularnewline
(c) Line plot of solution at $t=2.8$ & (d) Line plot of solution at $t=2.9$\tabularnewline
\includegraphics[width=0.48\columnwidth]{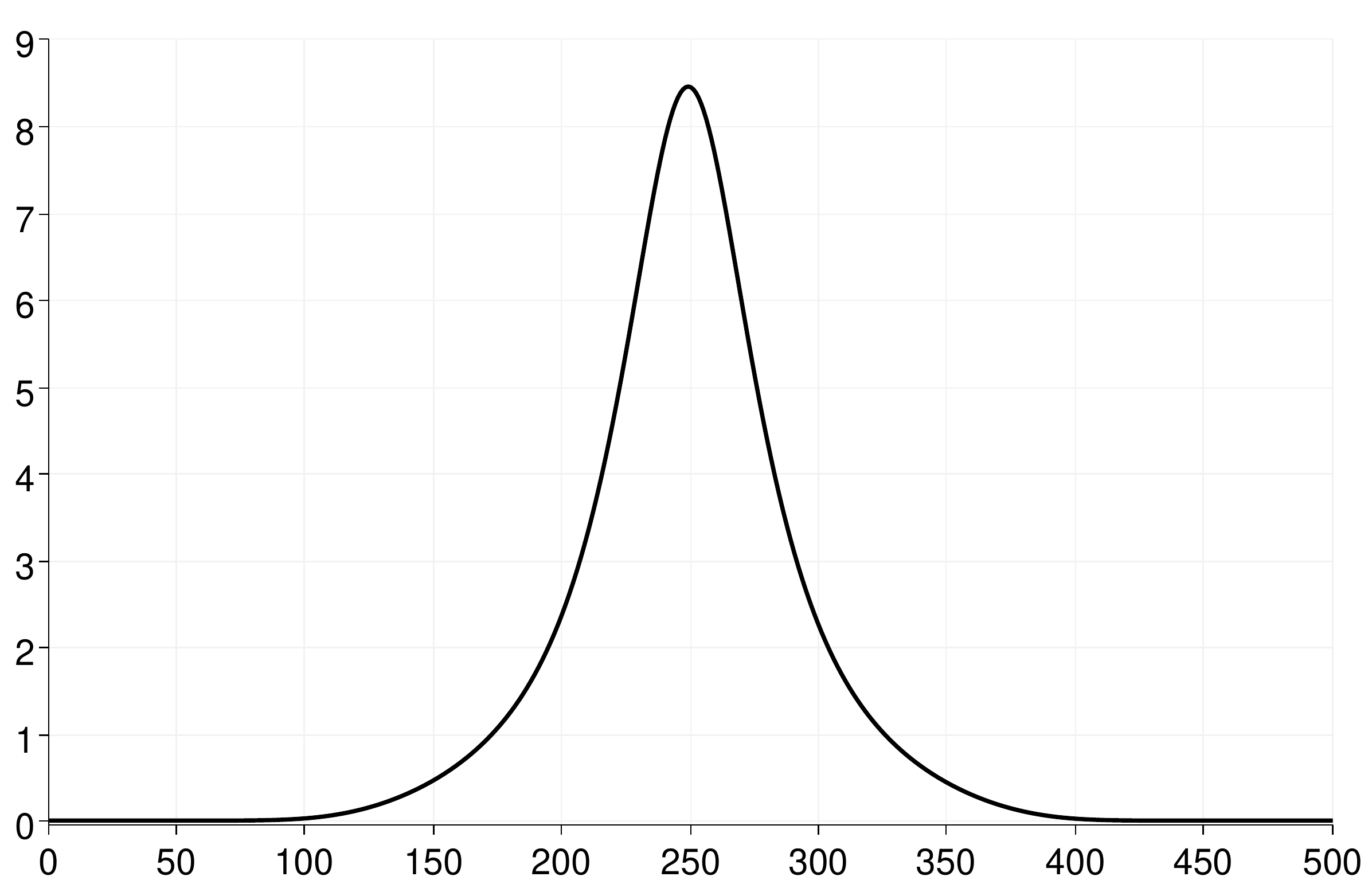} & \includegraphics[width=0.48\columnwidth]{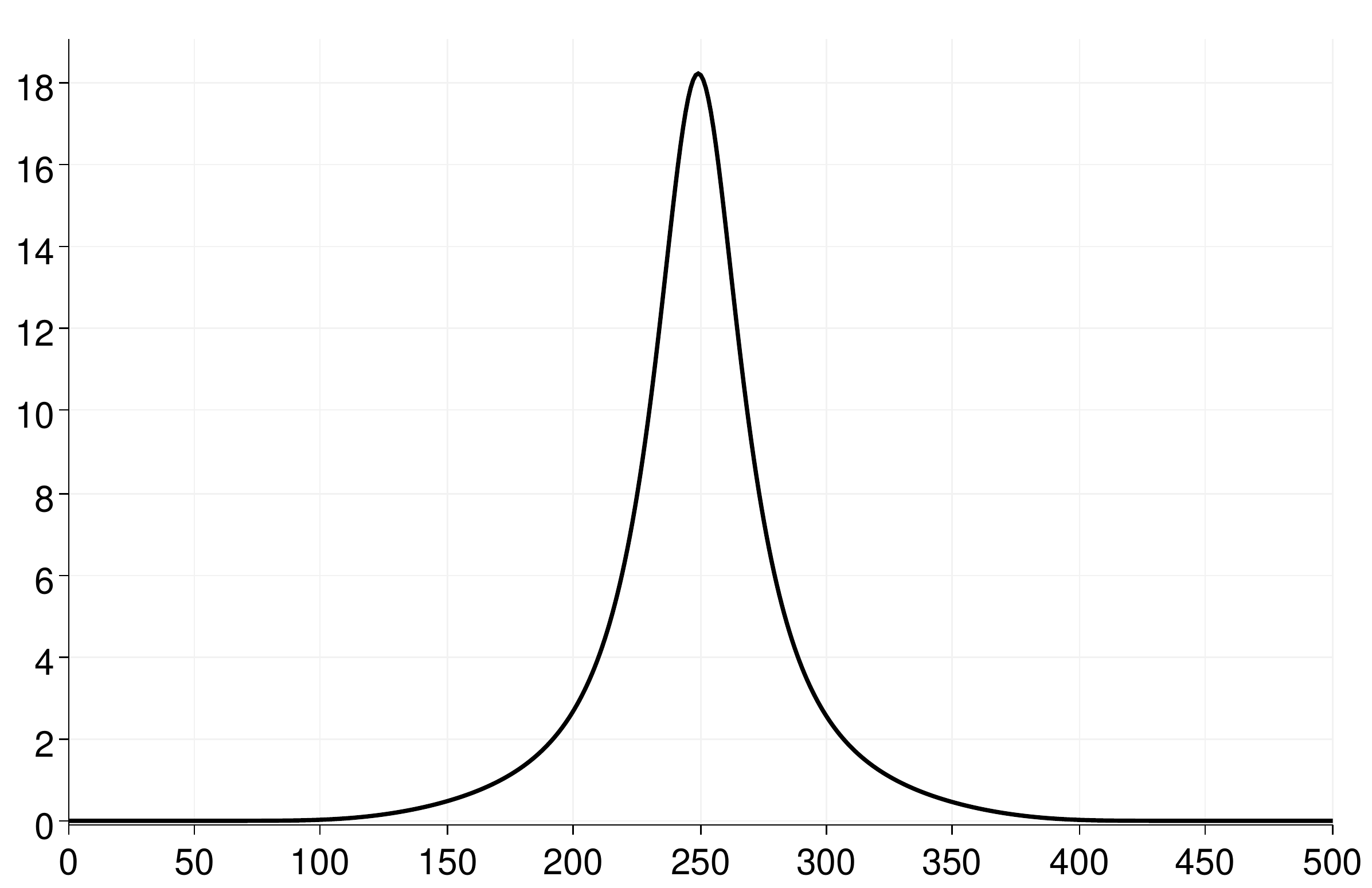}\tabularnewline
 & \tabularnewline
(e) Line plot of solution at $t=2.97$ & (f) Blow up of ${\displaystyle \int_{\Omega}\phi\left(\vec{x},t\right)d\vec{x}}$
around $t=2.98$\tabularnewline
\includegraphics[width=0.48\columnwidth]{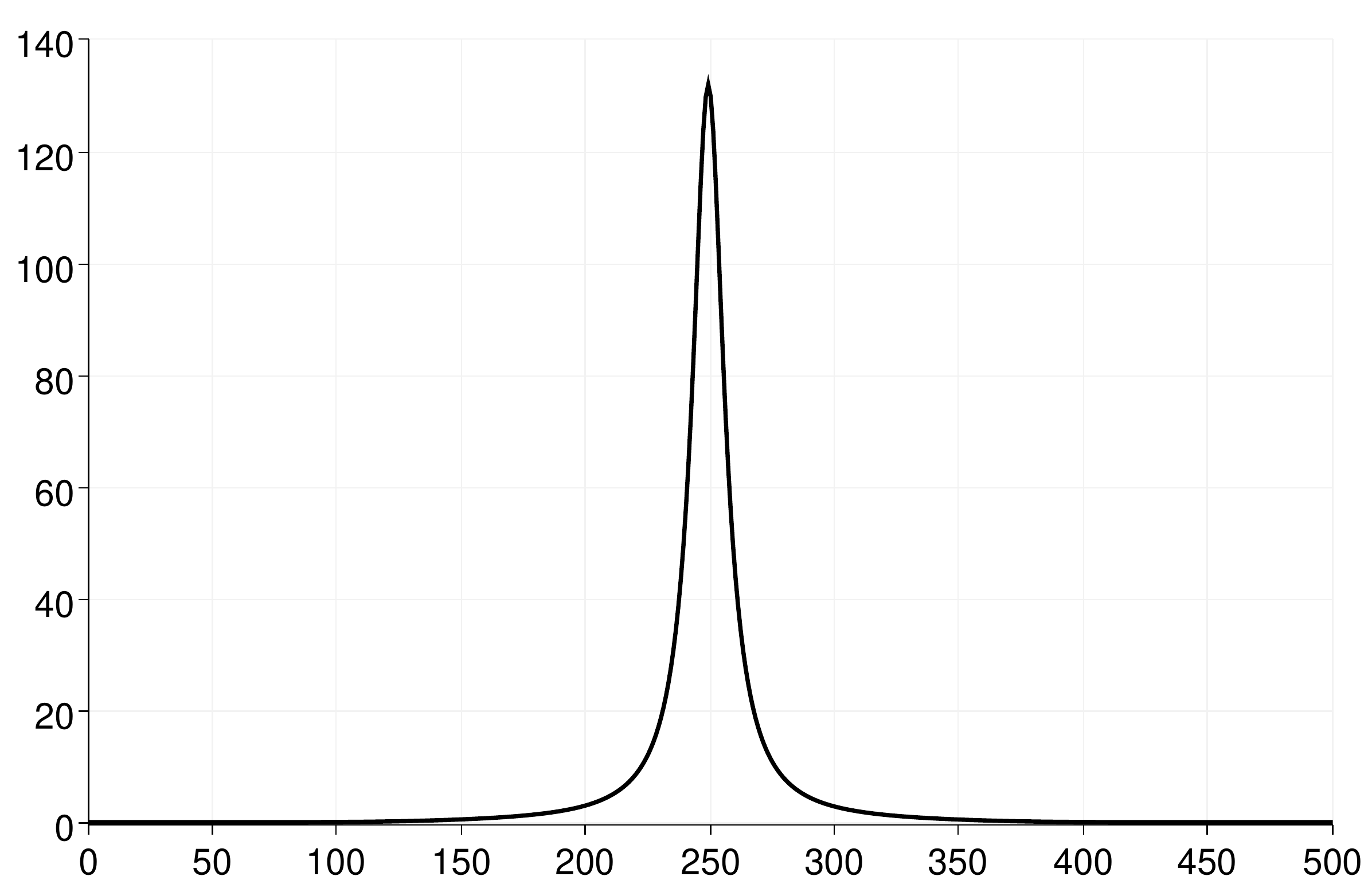} & \includegraphics[width=0.48\columnwidth]{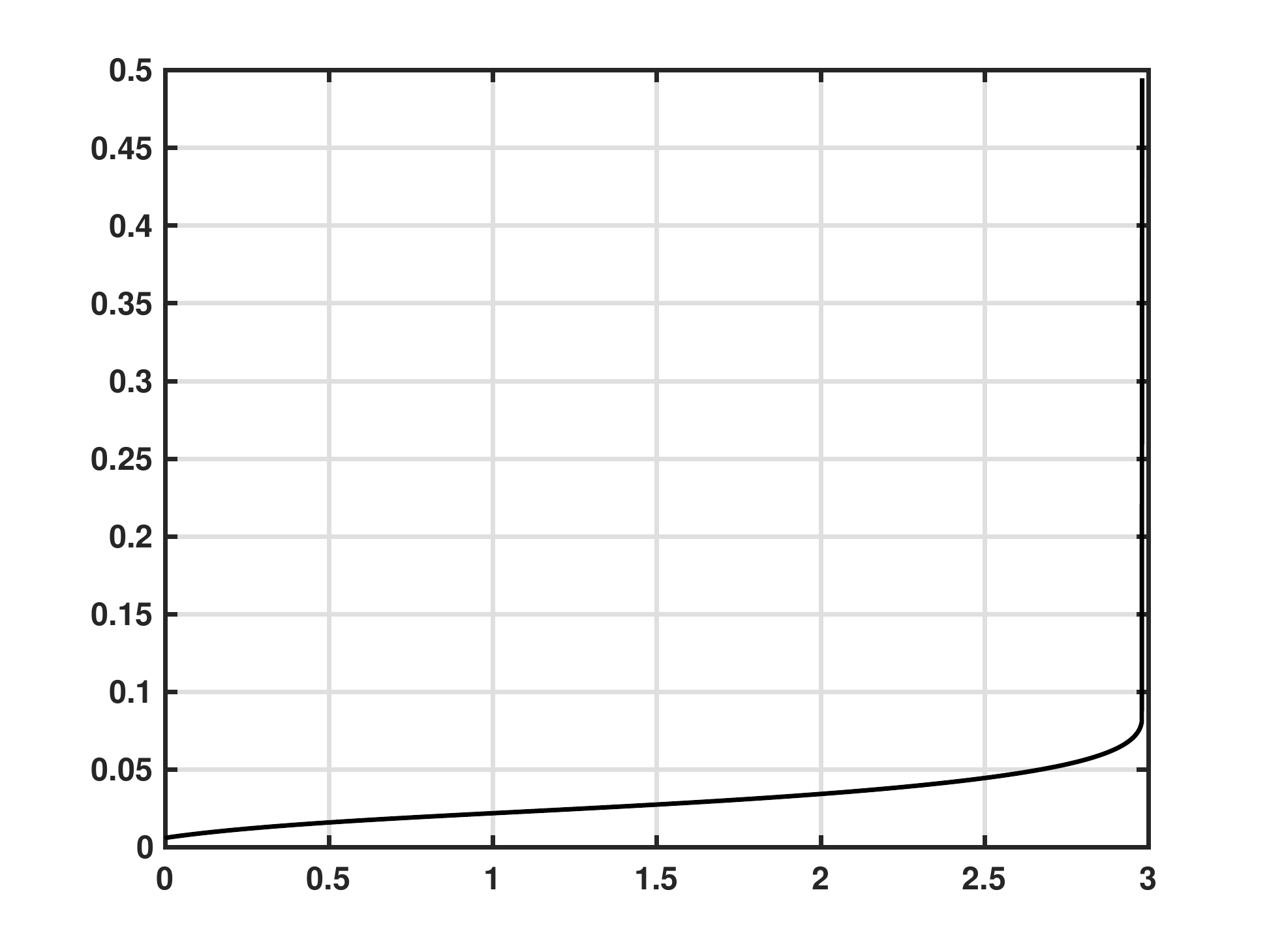}\tabularnewline
\end{tabular}
\end{figure}

\begin{figure}[H]
\caption{\label{fig:Bubble-Line}Bubble formation. }
\centering{}%
\begin{tabular}{cc}
(a) Line plot of solution at $t=0$  & (b) Line plot of solution at $t=0.21$\tabularnewline
\includegraphics[width=0.48\columnwidth]{bubble-lineplot-000} & \includegraphics[width=0.48\columnwidth]{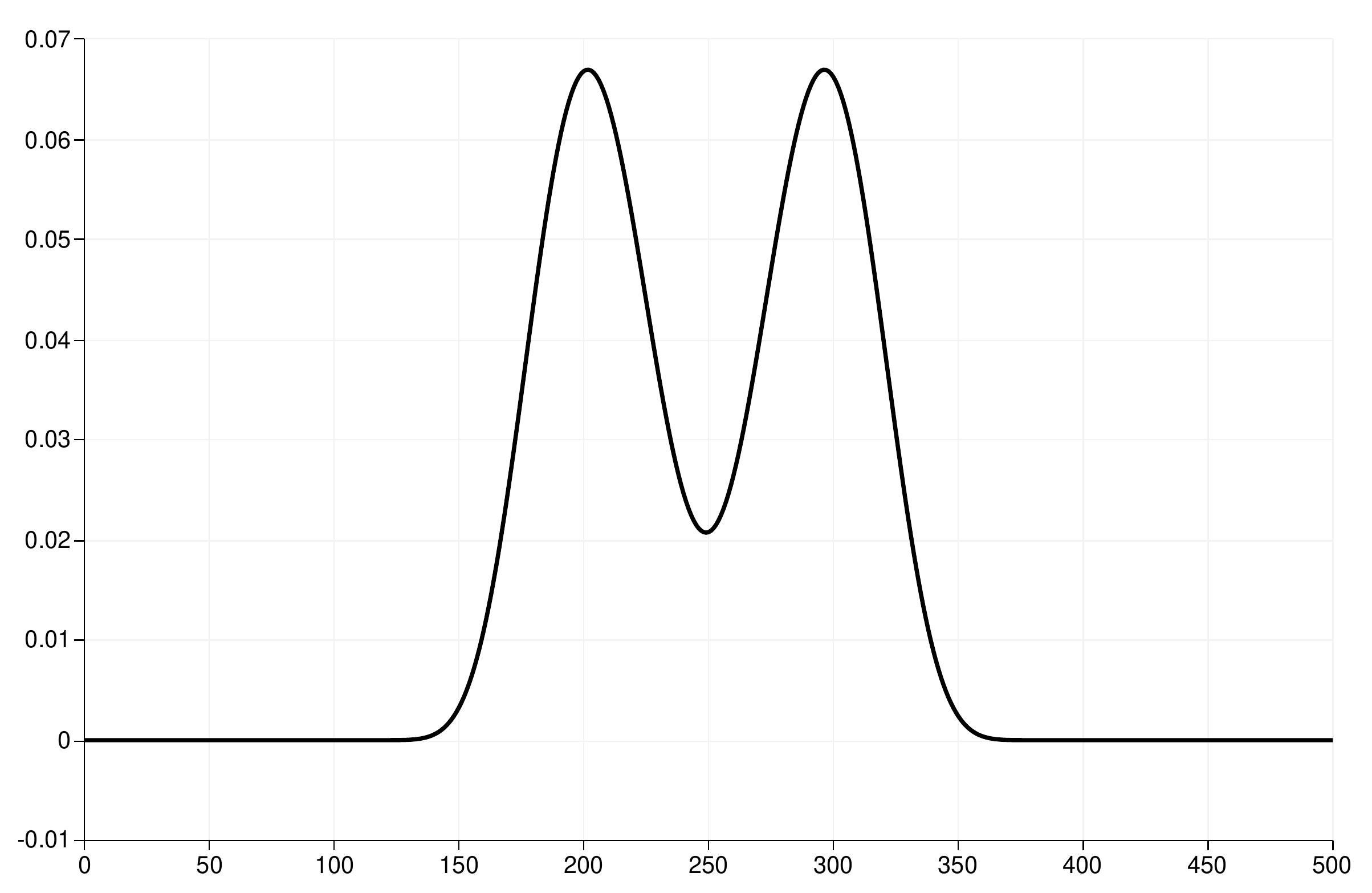}\tabularnewline
 & \tabularnewline
(c) Line plot of solution at $t=0.22$ & (d) Line plot of solution at $t=0.23$\tabularnewline
\includegraphics[width=0.48\columnwidth]{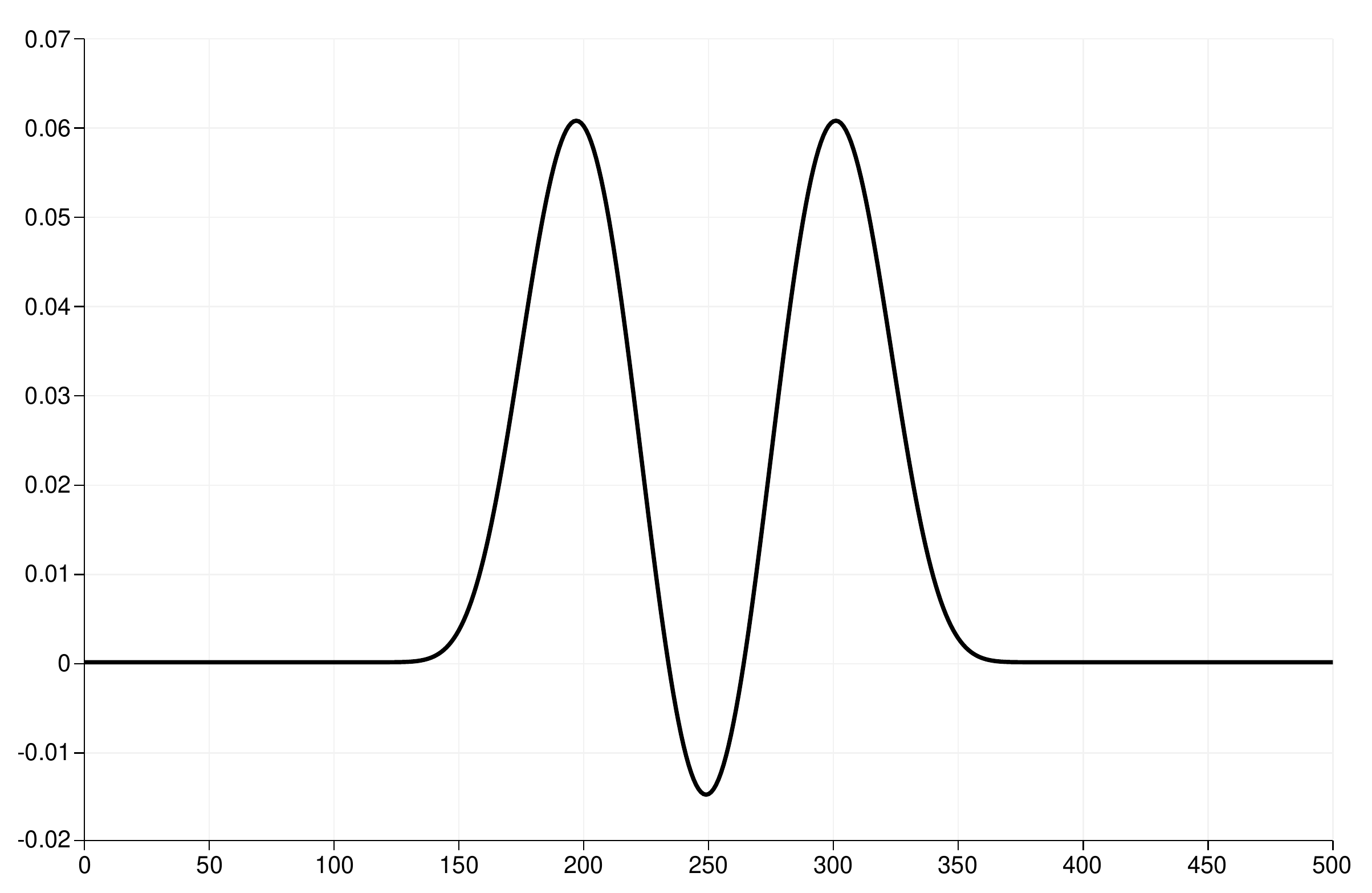} & \includegraphics[width=0.48\columnwidth]{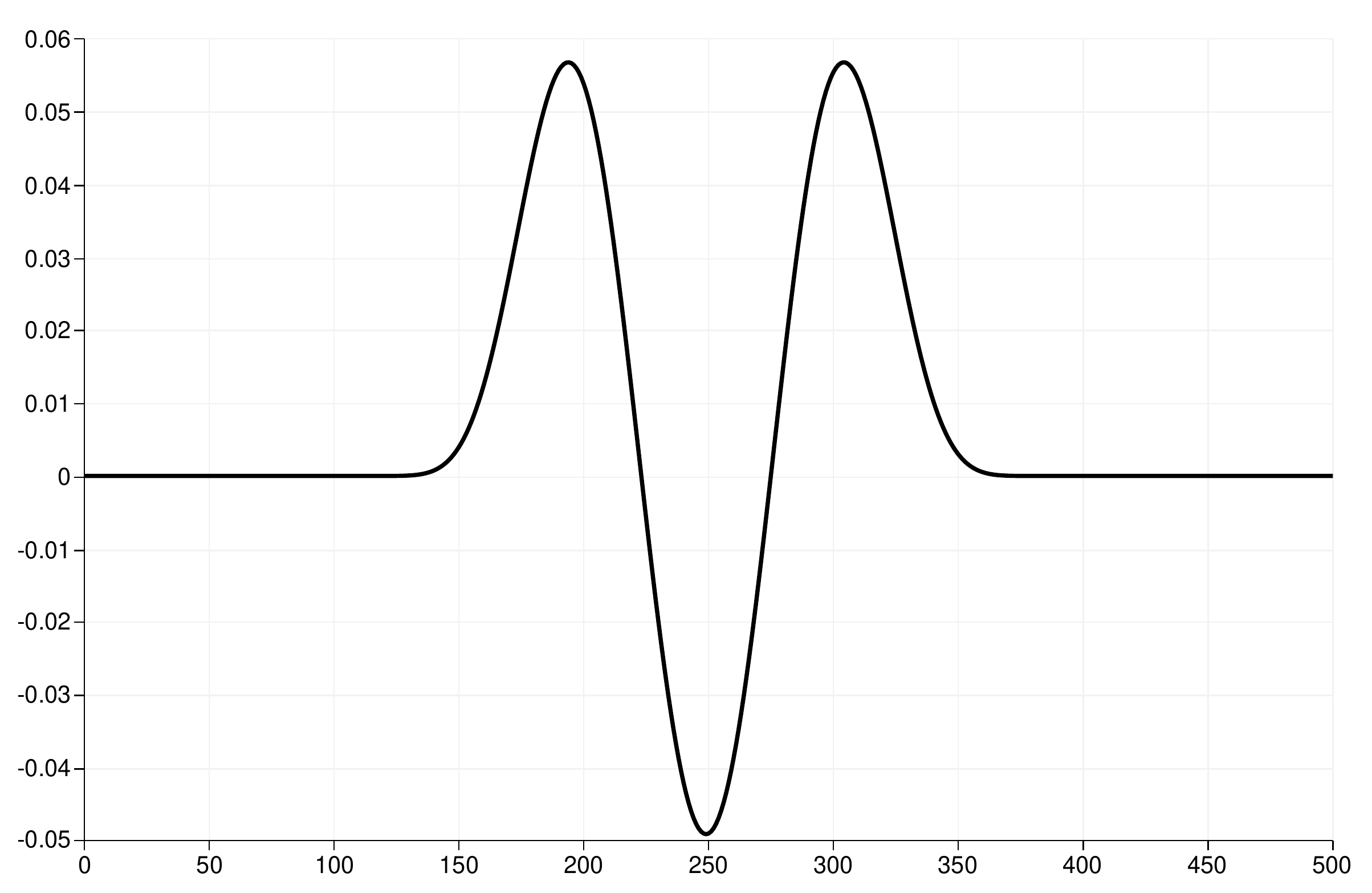}\tabularnewline
 & \tabularnewline
(e) Line plot of solution at $t=0.4$ & (f) 3D bubble at $t=0.4$\tabularnewline
\includegraphics[width=0.48\columnwidth]{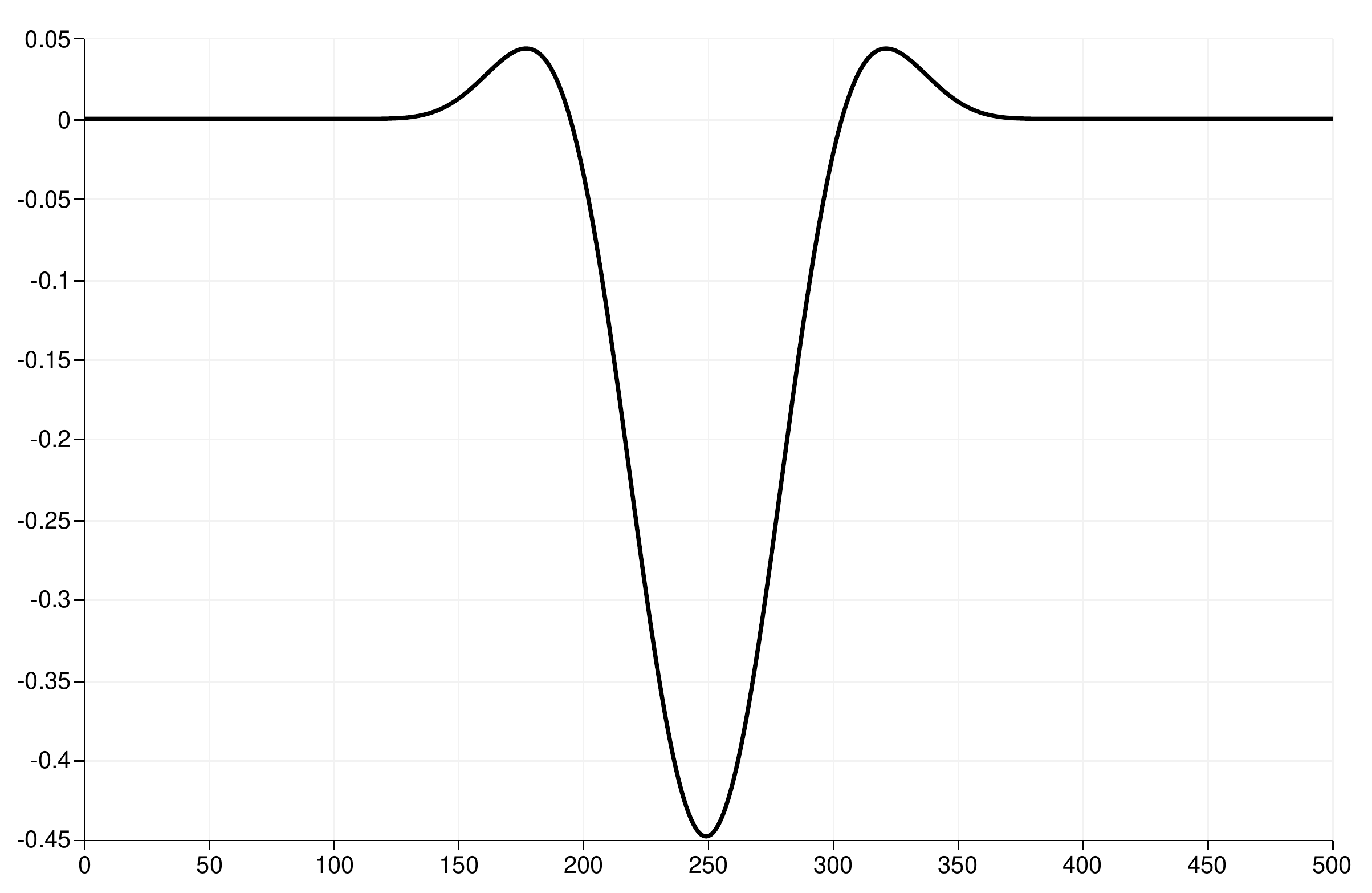} & \includegraphics[width=0.48\columnwidth]{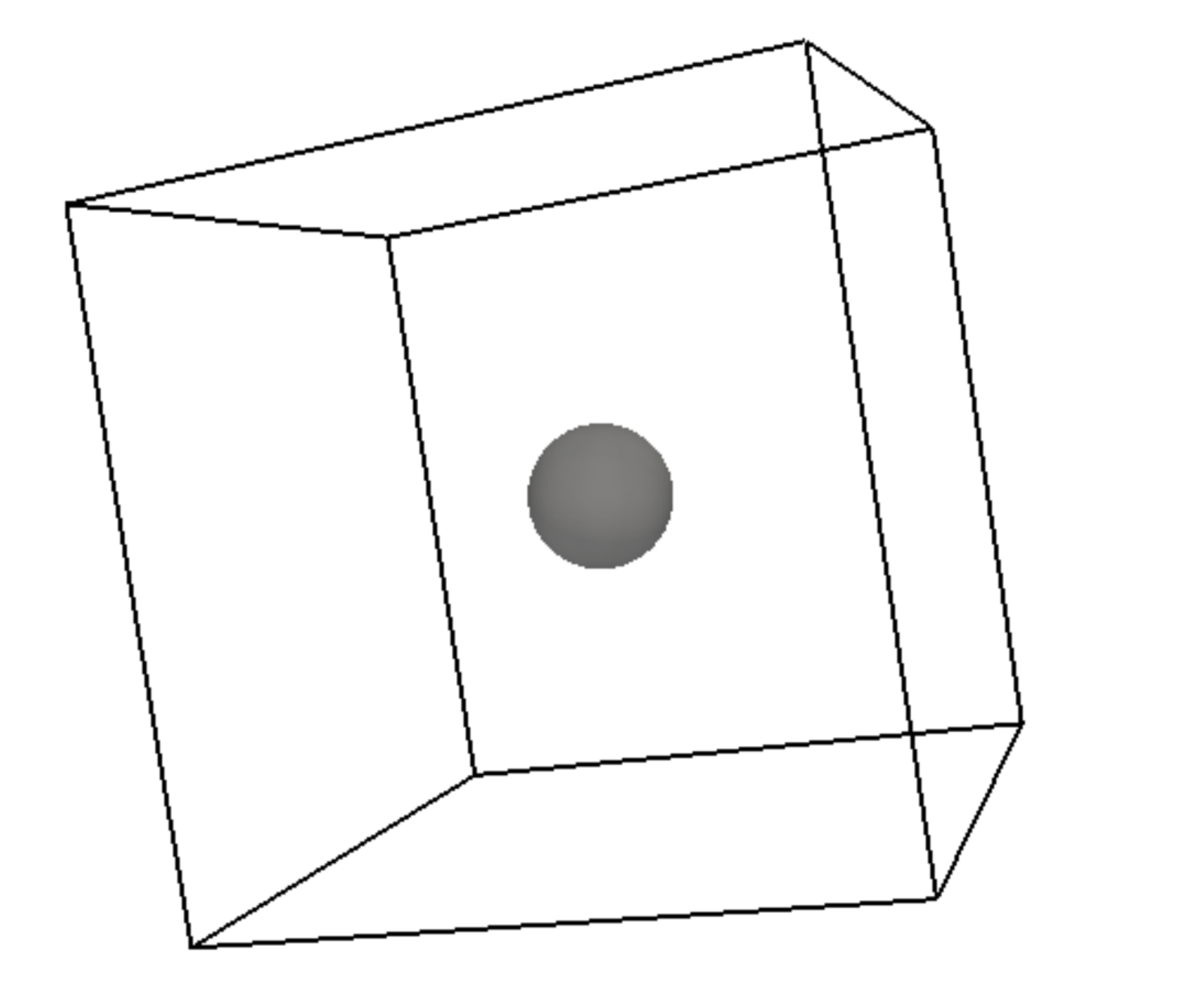}\tabularnewline
\end{tabular}
\end{figure}

\begin{figure}[H]
\caption{\label{fig:Bubble-long-time}Long time behavior of solution with a
bubble. }
\centering{}%
\begin{tabular}{cc}
(a) Solution along a line at $t=0$  & (b) Solution along a line at $t=0.22$\tabularnewline
\includegraphics[width=0.48\columnwidth]{bubble-lineplot-000} & \includegraphics[width=0.48\columnwidth]{bubble-lineplot-022}\tabularnewline
 & \tabularnewline
(c) Solution along a line at $t=0.5$ & (d) Solution along a line at $t=2$\tabularnewline
\includegraphics[width=0.48\columnwidth]{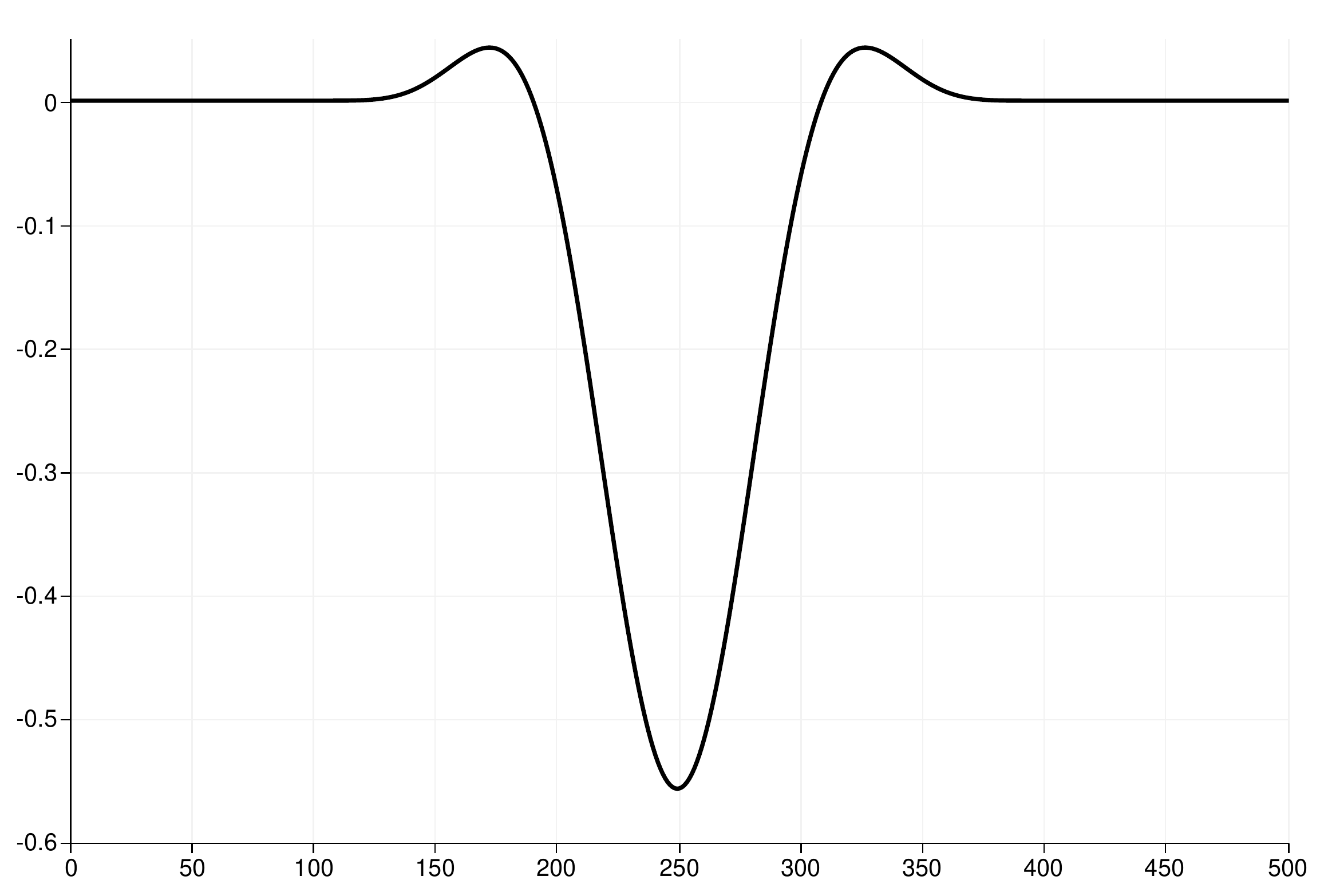} & \includegraphics[width=0.48\columnwidth]{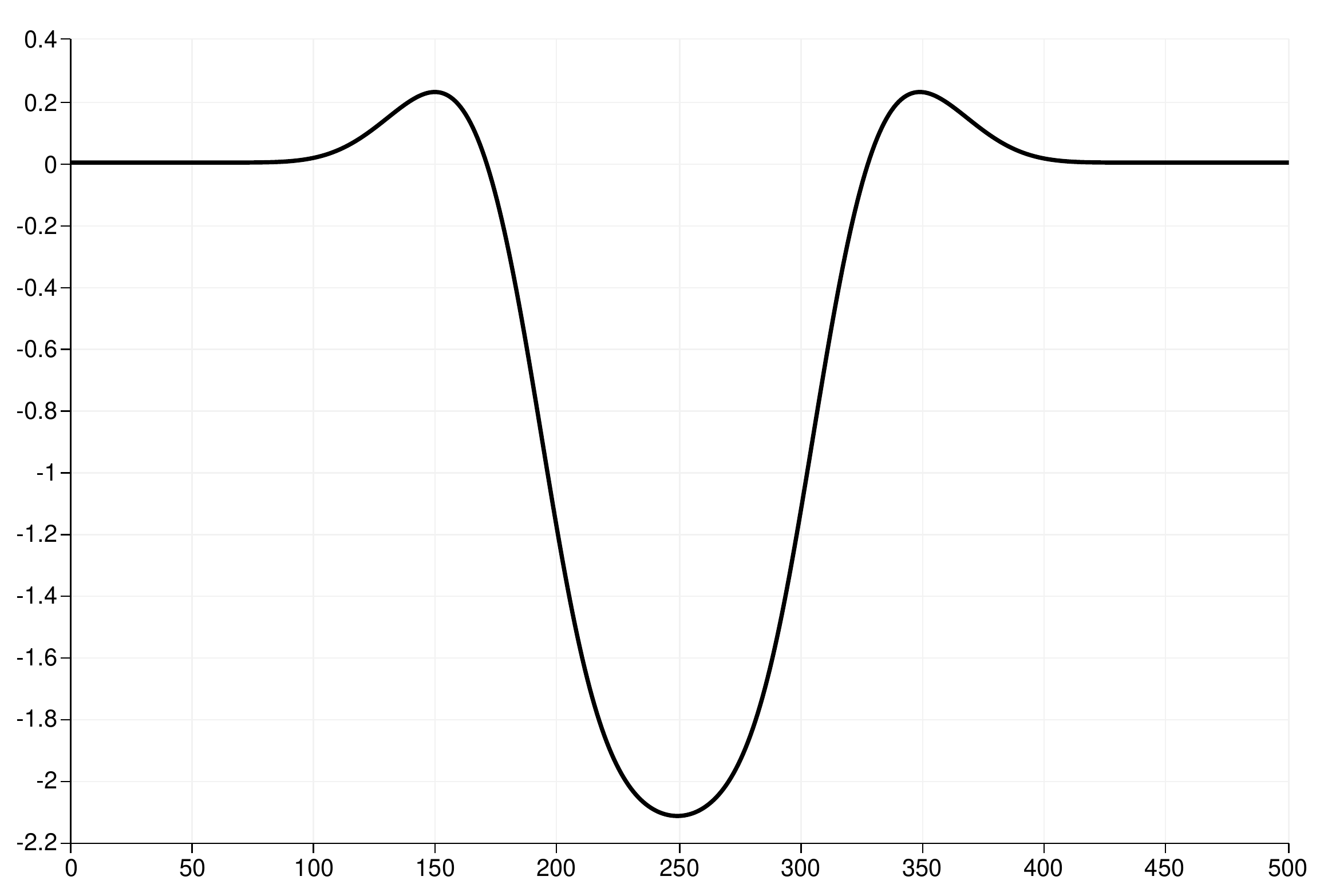}\tabularnewline
 & \tabularnewline
(e) Solution along a line at $t=5$ & (f) Solution along a line at $t=7$\tabularnewline
\includegraphics[width=0.48\columnwidth]{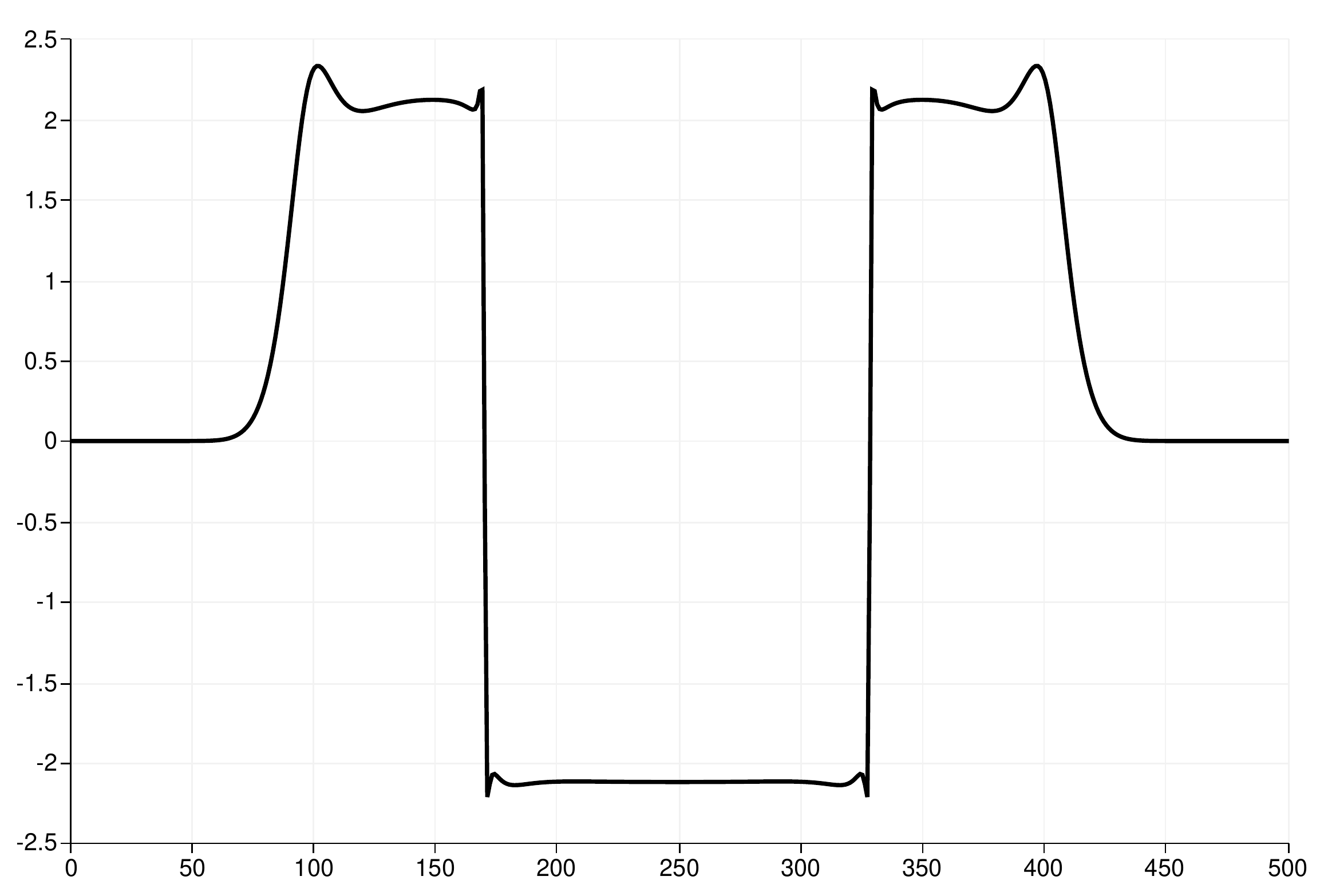} & \includegraphics[width=0.48\columnwidth]{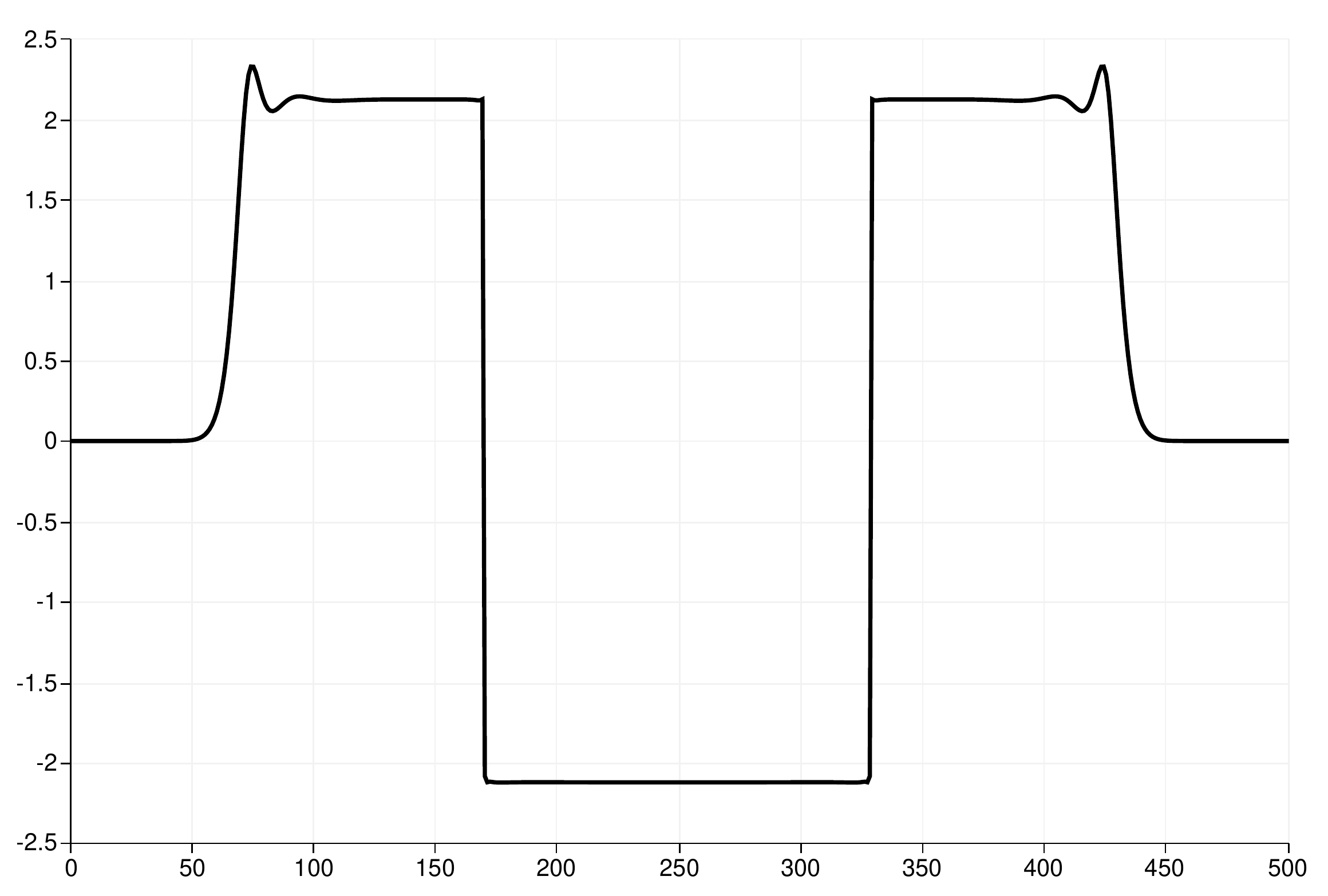}\tabularnewline
\end{tabular}
\end{figure}

\begin{figure}[H]
\caption{\label{fig:Laplacian}$\mathcal{P}\left(t\right)\equiv{\displaystyle \frac{1}{L^{2}}e^{-2t}\max_{\vec{x}\in\Omega}\left|\Delta\phi\left(\vec{x},t\right)\right|}$}
\centering{}\includegraphics[width=0.5\columnwidth]{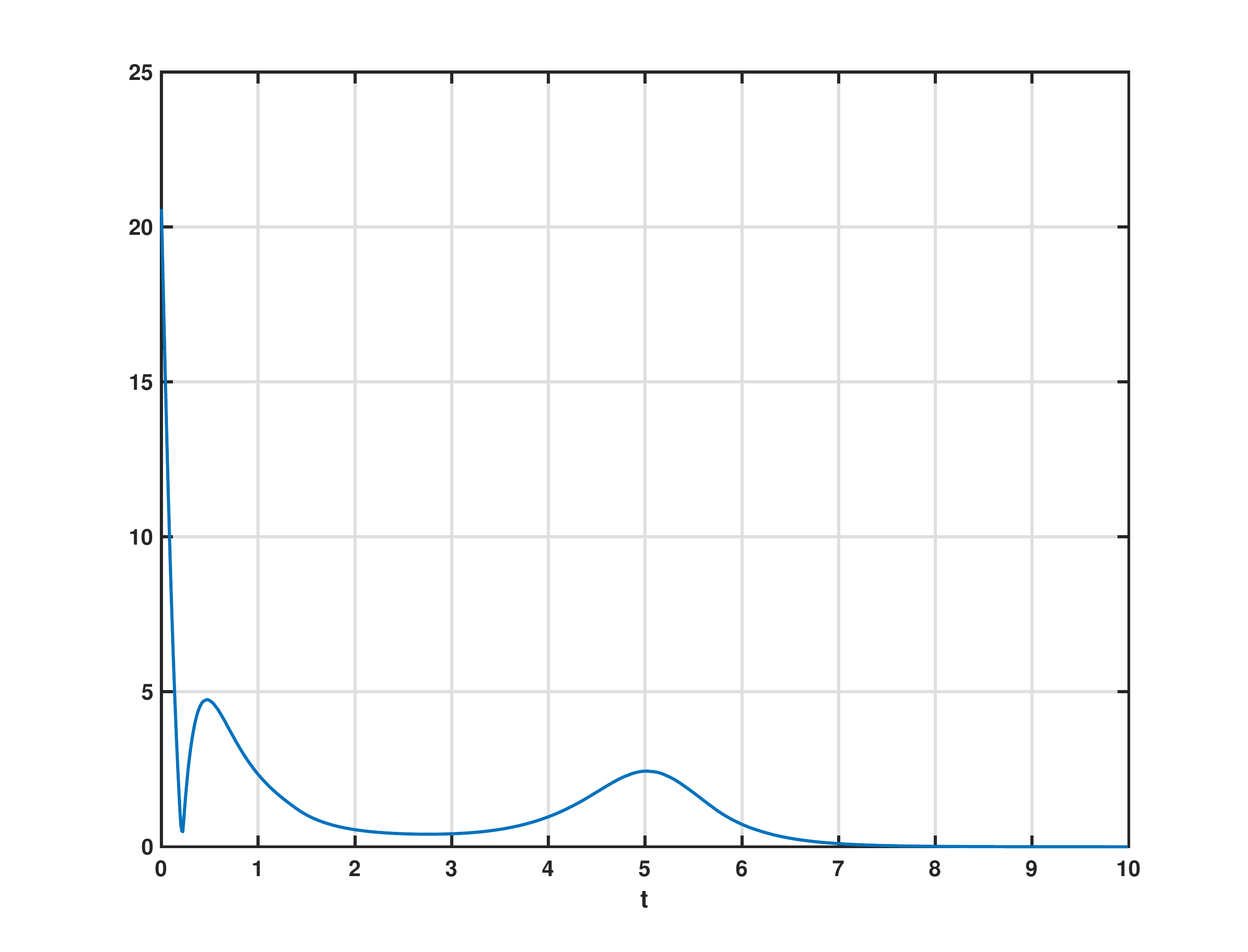}
\end{figure}

\begin{figure}[H]
\caption{\label{fig:Duffing}Phase portrait of the unforced, damped Duffing
equation.}
\centering{}\includegraphics[height=8cm]{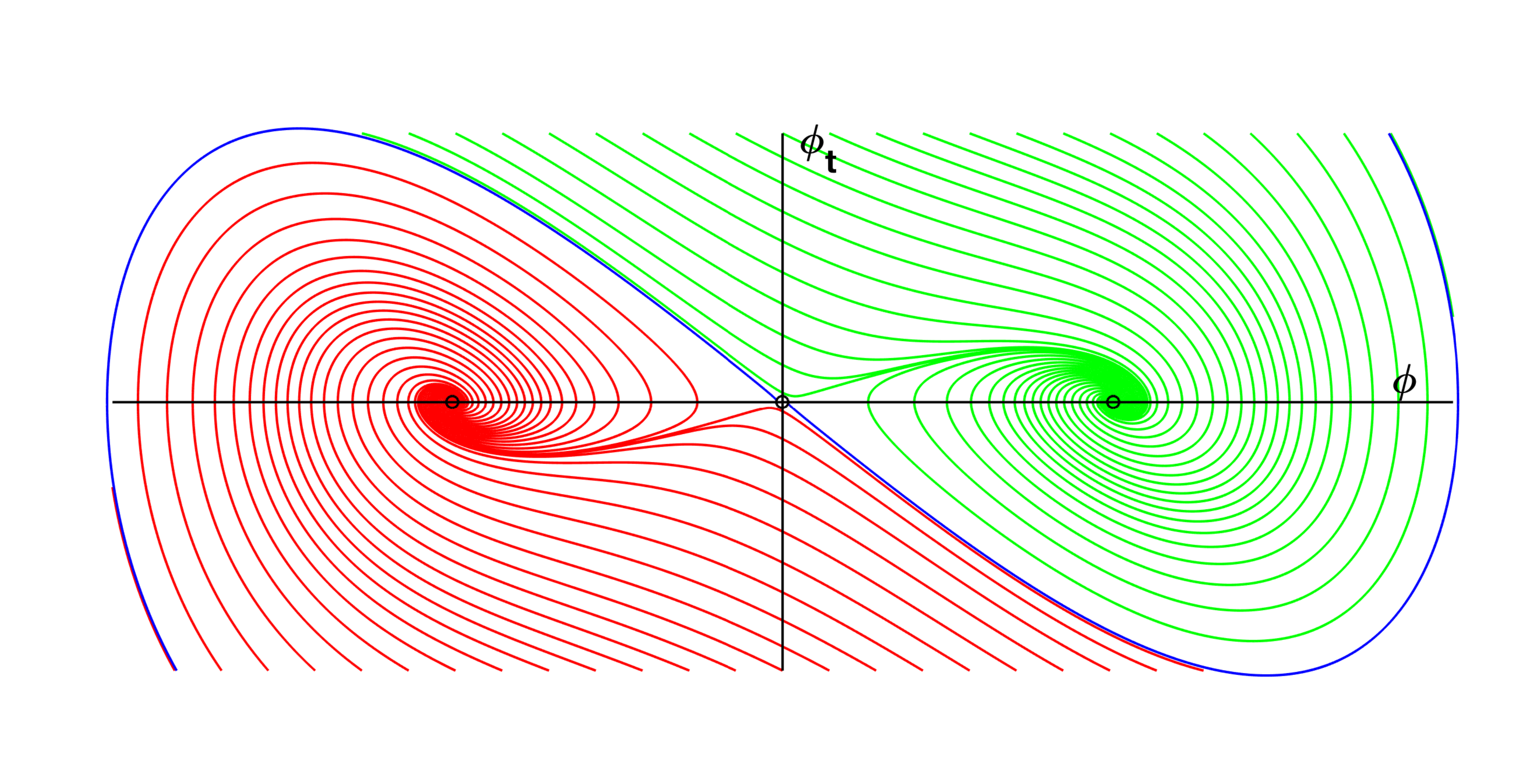}
\end{figure}

\begin{figure}[H]
\caption{\label{fig:nobubble}Lack of bubble formation.}
\centering{}%
\begin{tabular}{cc}
(a) Solution along a line at $t=0$  & (b) Solution along a line at $t=1$\tabularnewline
\includegraphics[width=0.45\columnwidth]{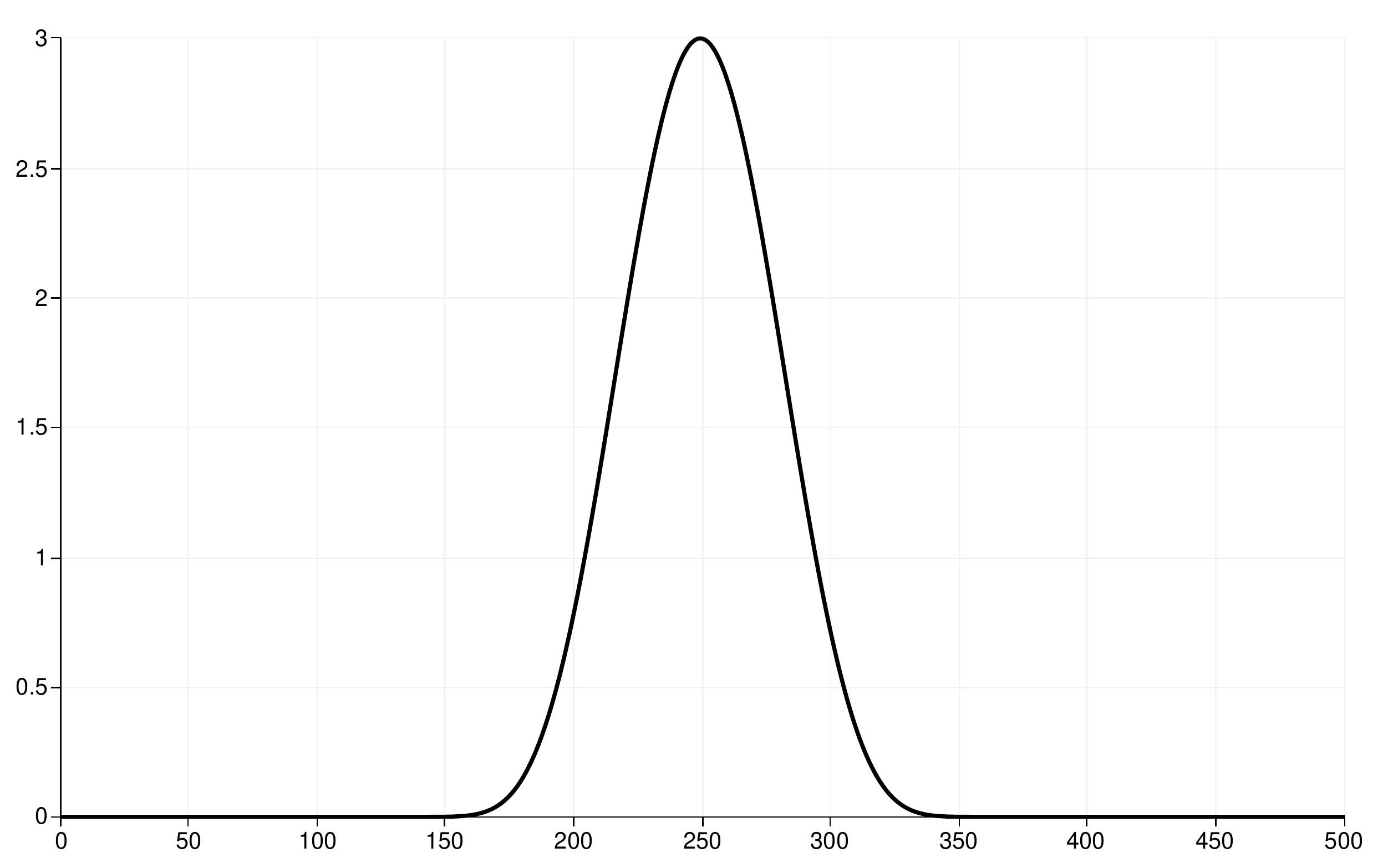} & \includegraphics[width=0.45\columnwidth]{nobubble-010}\tabularnewline
 & \tabularnewline
(c) Solution along a line at $t=2$ & (d) Solution along a line at $t=3$\tabularnewline
\includegraphics[width=0.45\columnwidth]{nobubble-020} & \includegraphics[width=0.45\columnwidth]{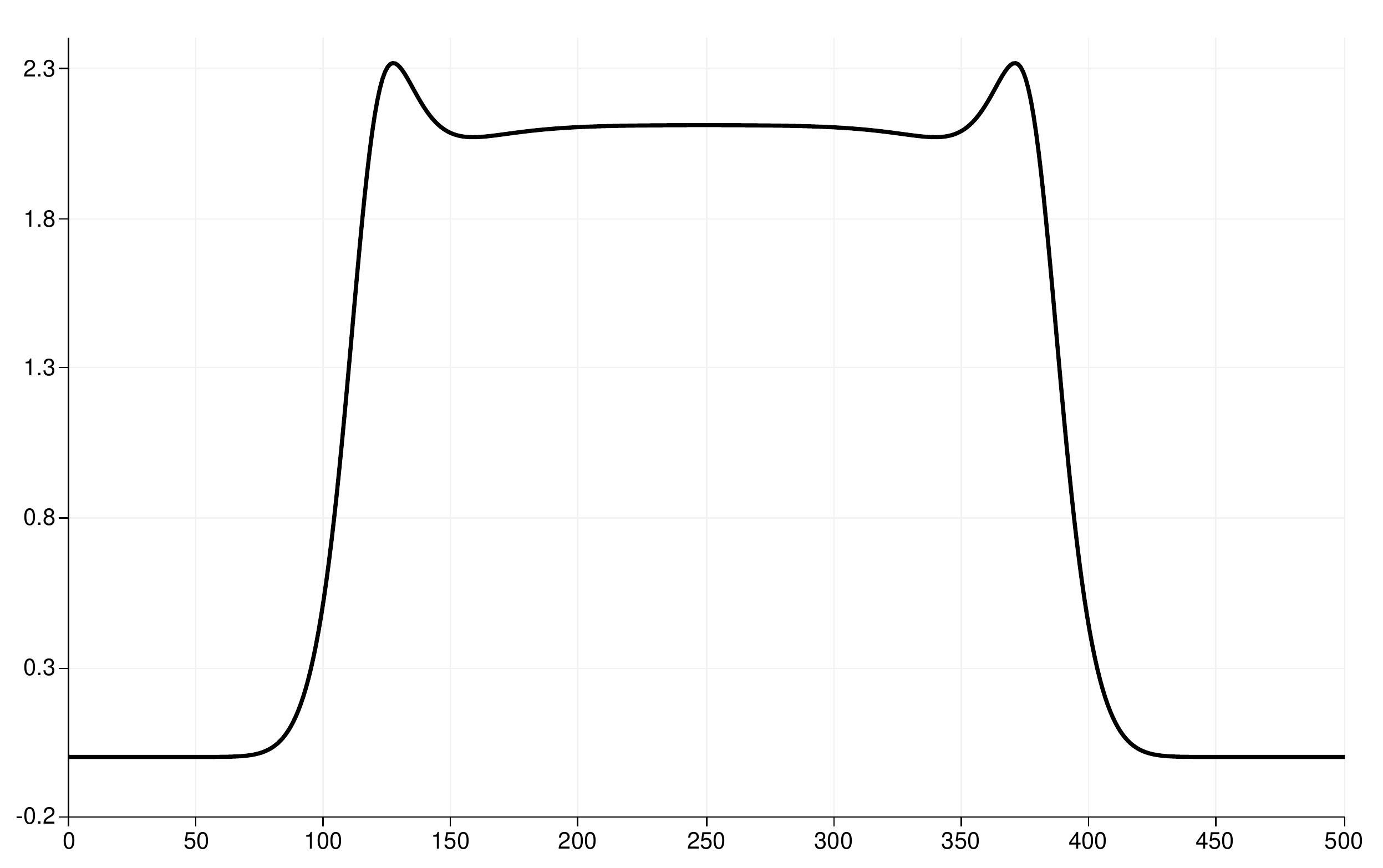}\tabularnewline
 & \tabularnewline
(e) Solution along a line at $t=4$ & (f) Solution along a line at $t=5$\tabularnewline
\includegraphics[width=0.45\columnwidth]{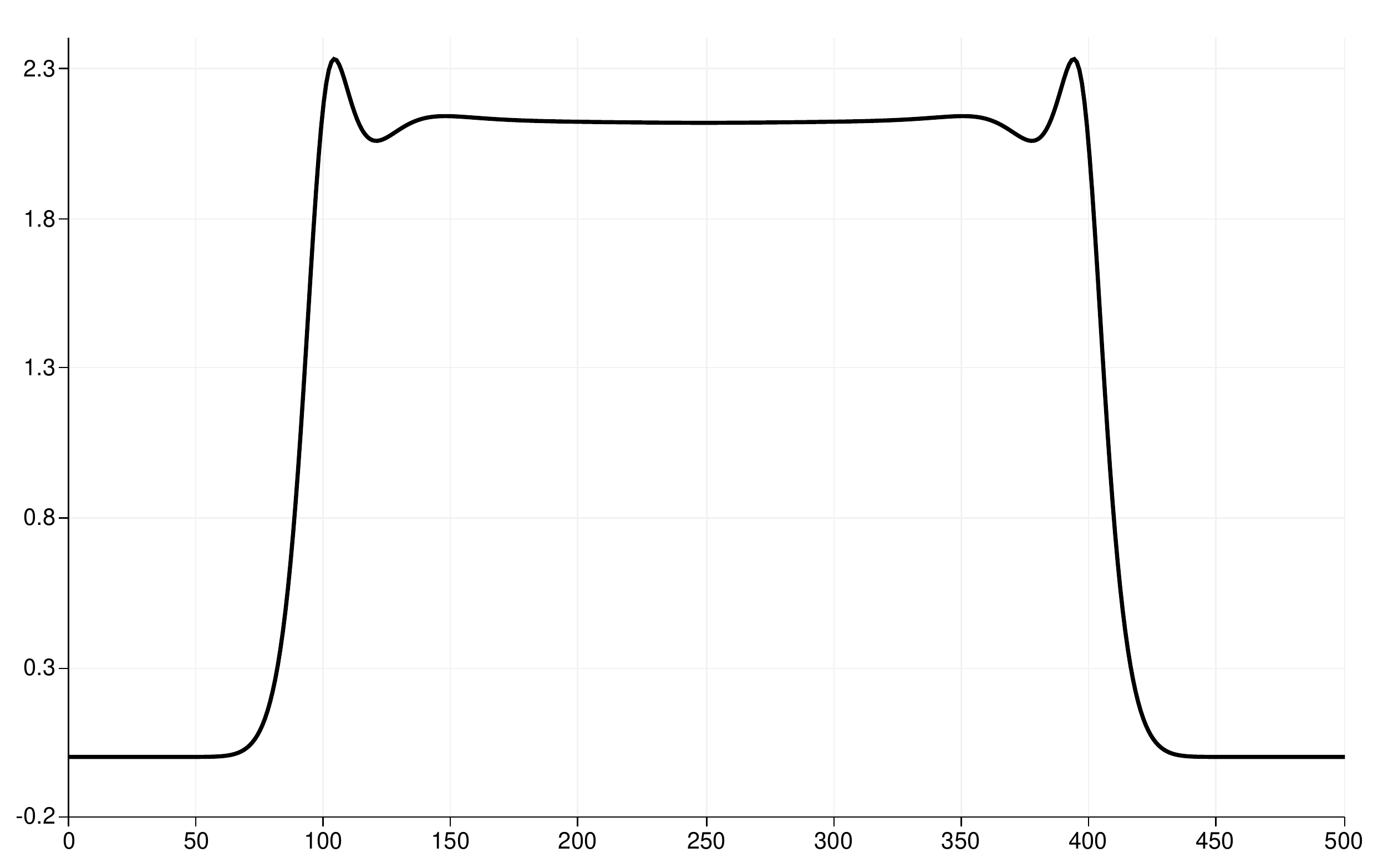} & \includegraphics[width=0.45\columnwidth]{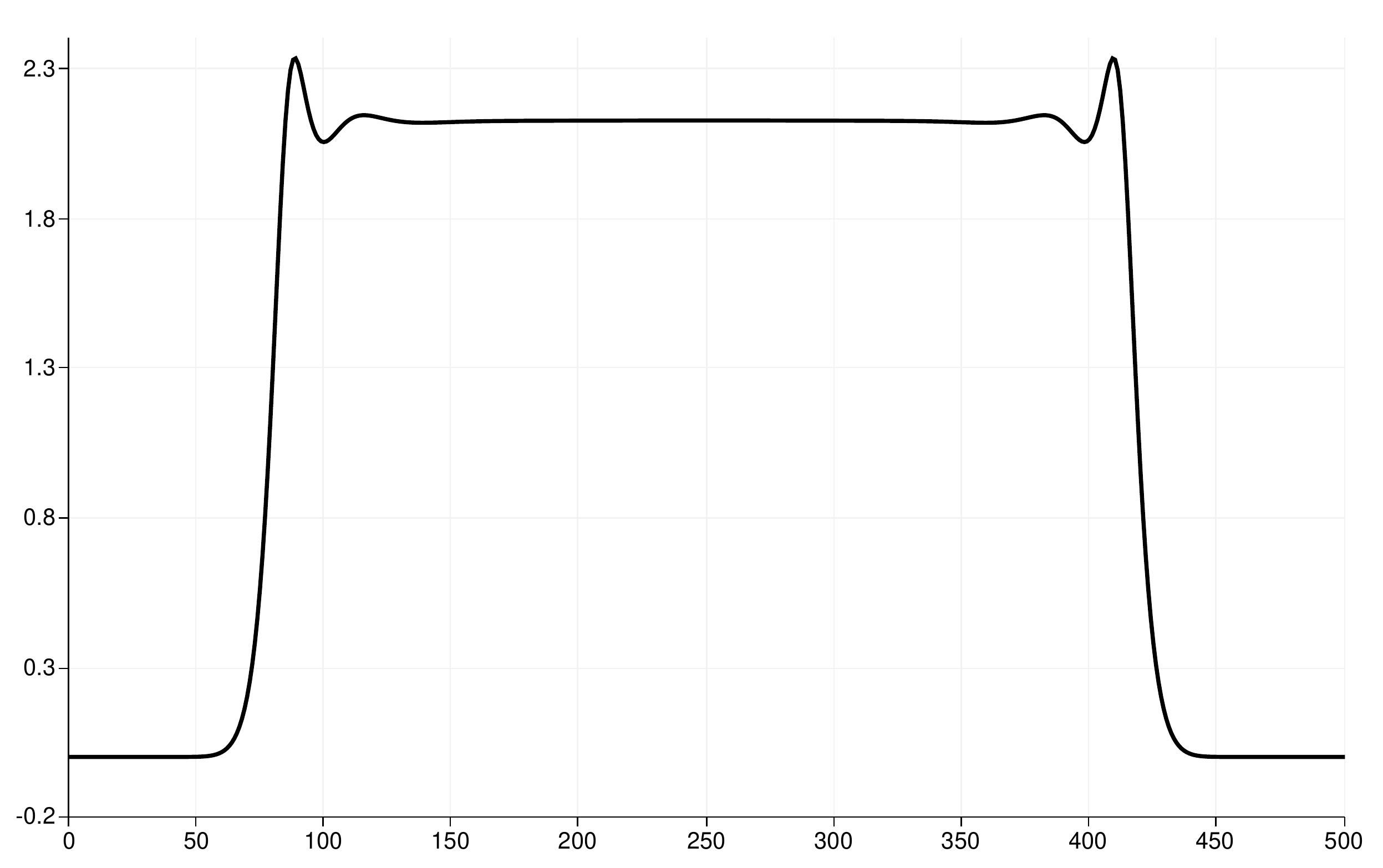}\tabularnewline
 & \tabularnewline
(g) Solution along a line at $t=6$ & (h) Solution along a line at $t=7$\tabularnewline
\includegraphics[width=0.45\columnwidth]{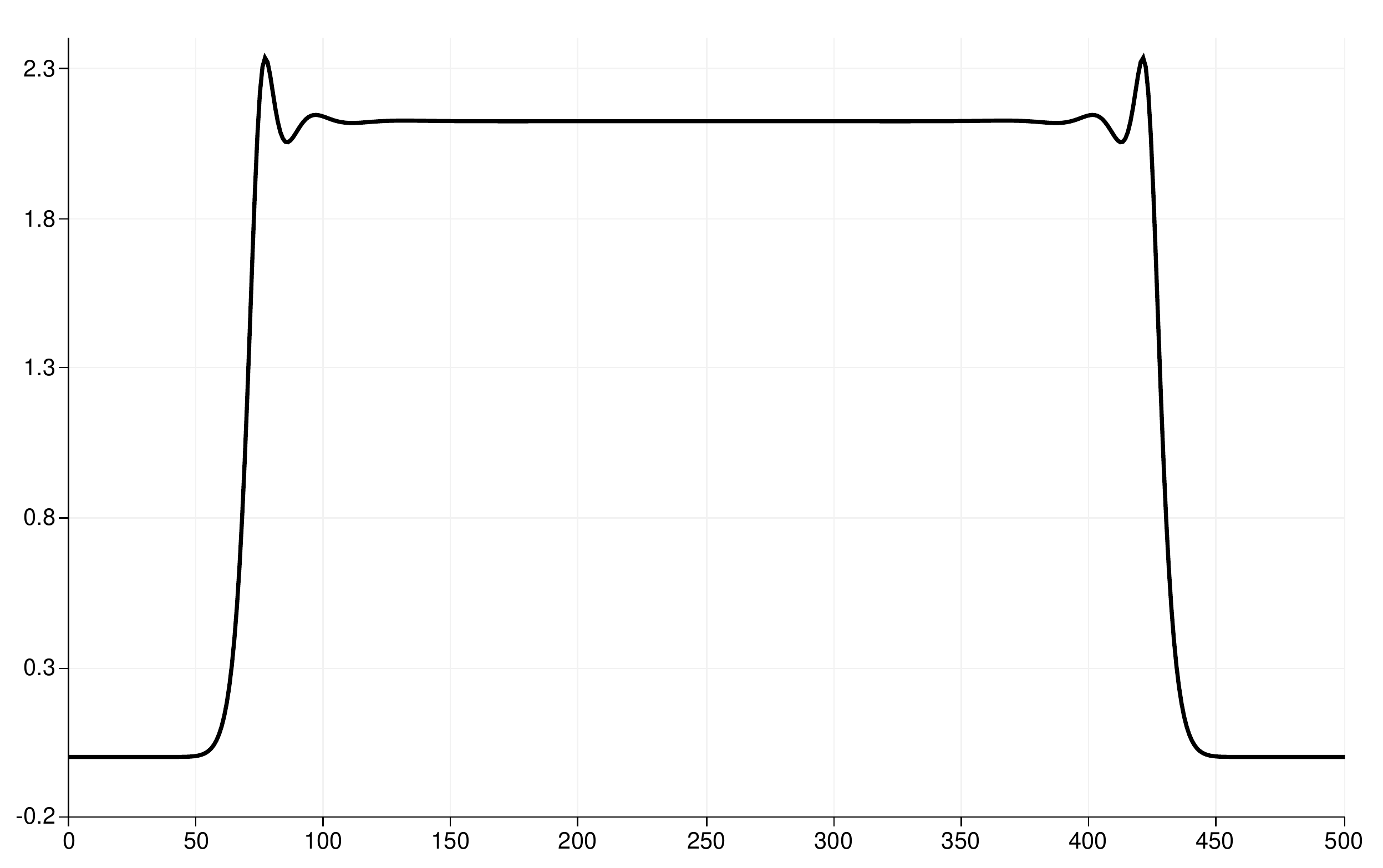} & \includegraphics[width=0.45\columnwidth]{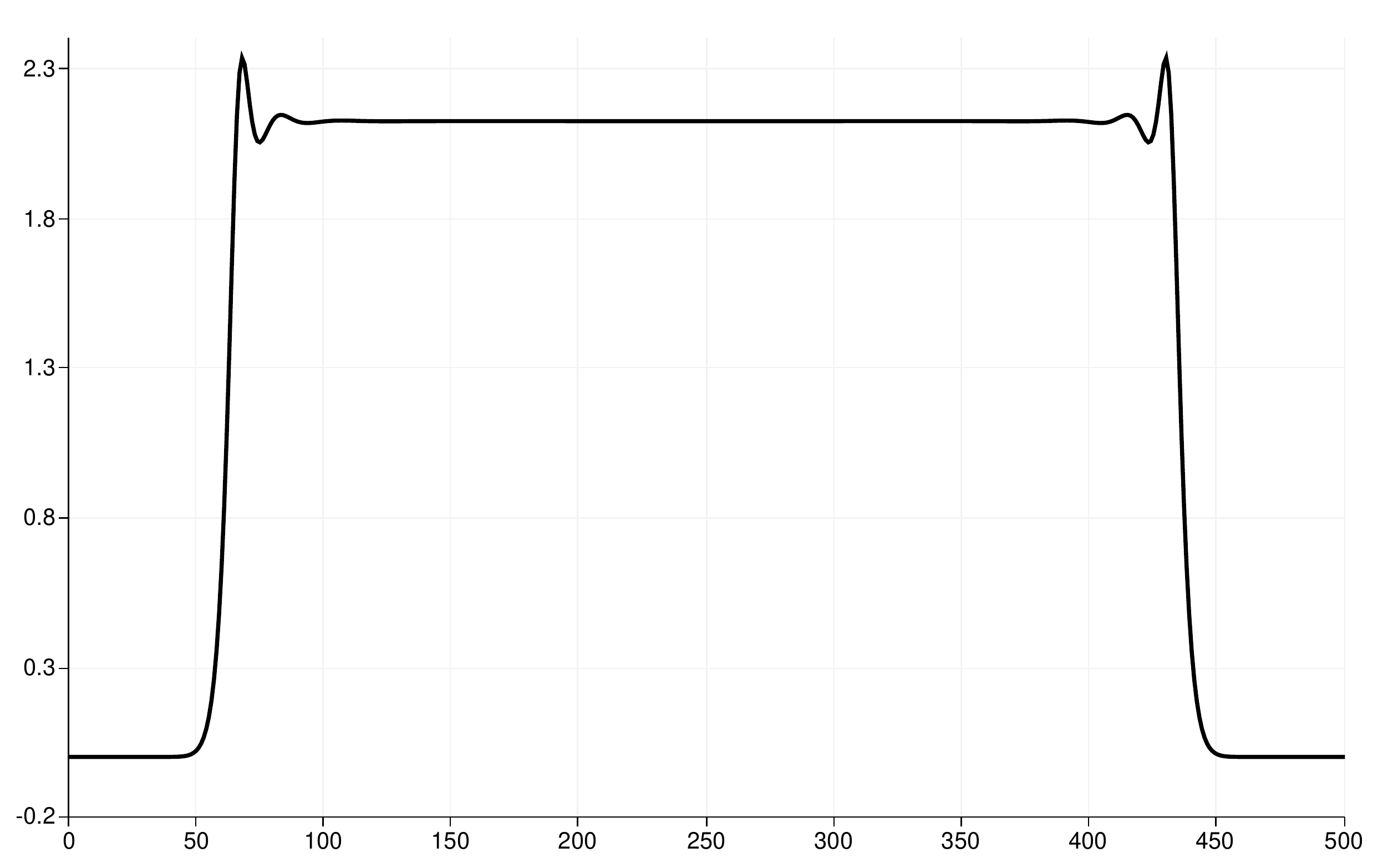}\tabularnewline
\end{tabular}
\end{figure}

\begin{figure}[H]
\caption{\label{fig:duffing_complex}Convergence to the unforced, damped Duffing
equation.}
\centering{}%
\begin{tabular}{cc}
\multicolumn{2}{c}{(a) Initial function $\varphi_{0}$ in the phase portrait of the Duffing
equation.}\tabularnewline
\multicolumn{2}{c}{\includegraphics[width=0.7\columnwidth]{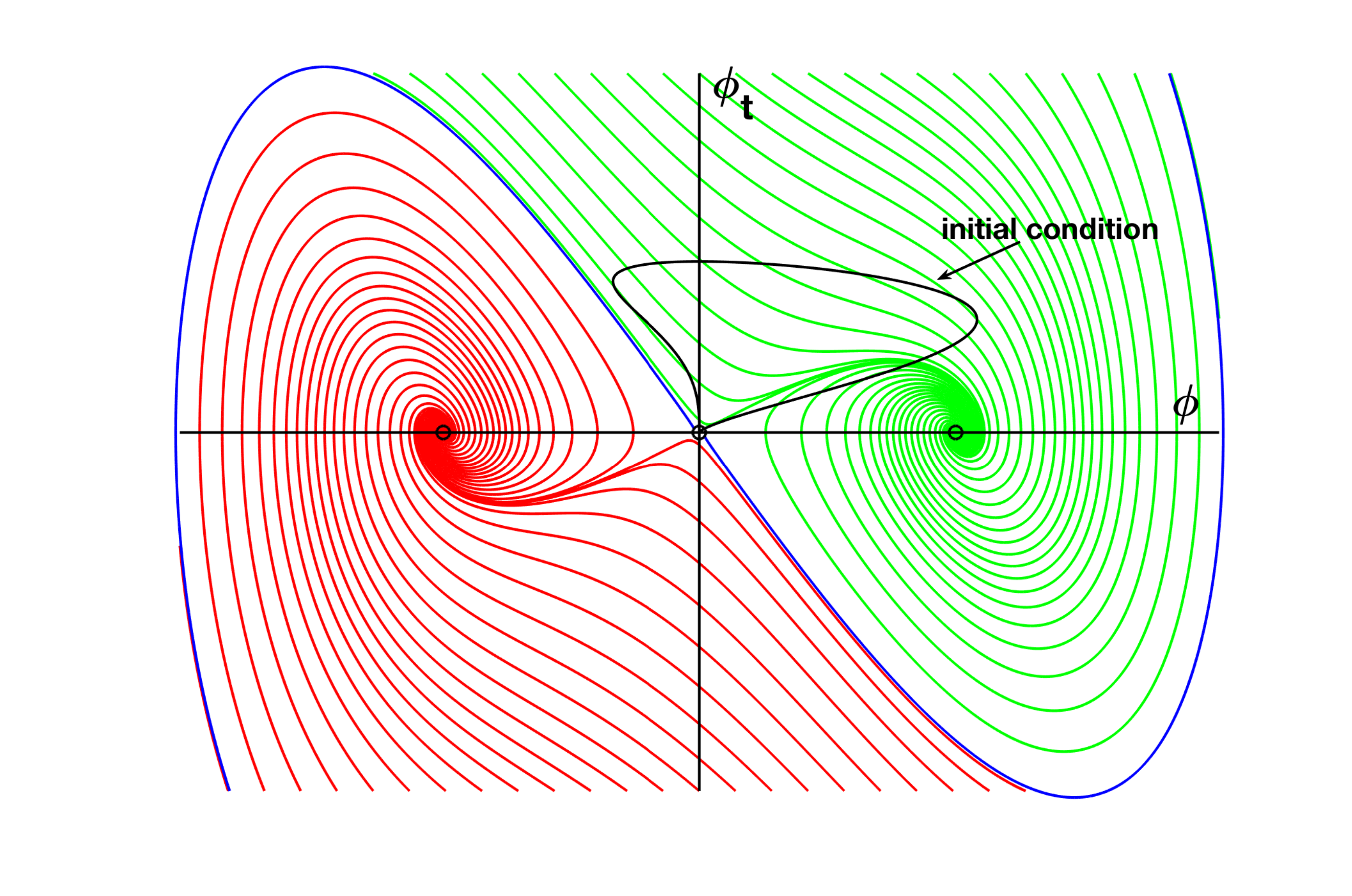}}\tabularnewline
(b) Line plot of initial function $\varphi_{0}$. & (c) Line plot of solution at time $t=4$.\tabularnewline
\includegraphics[width=0.48\columnwidth]{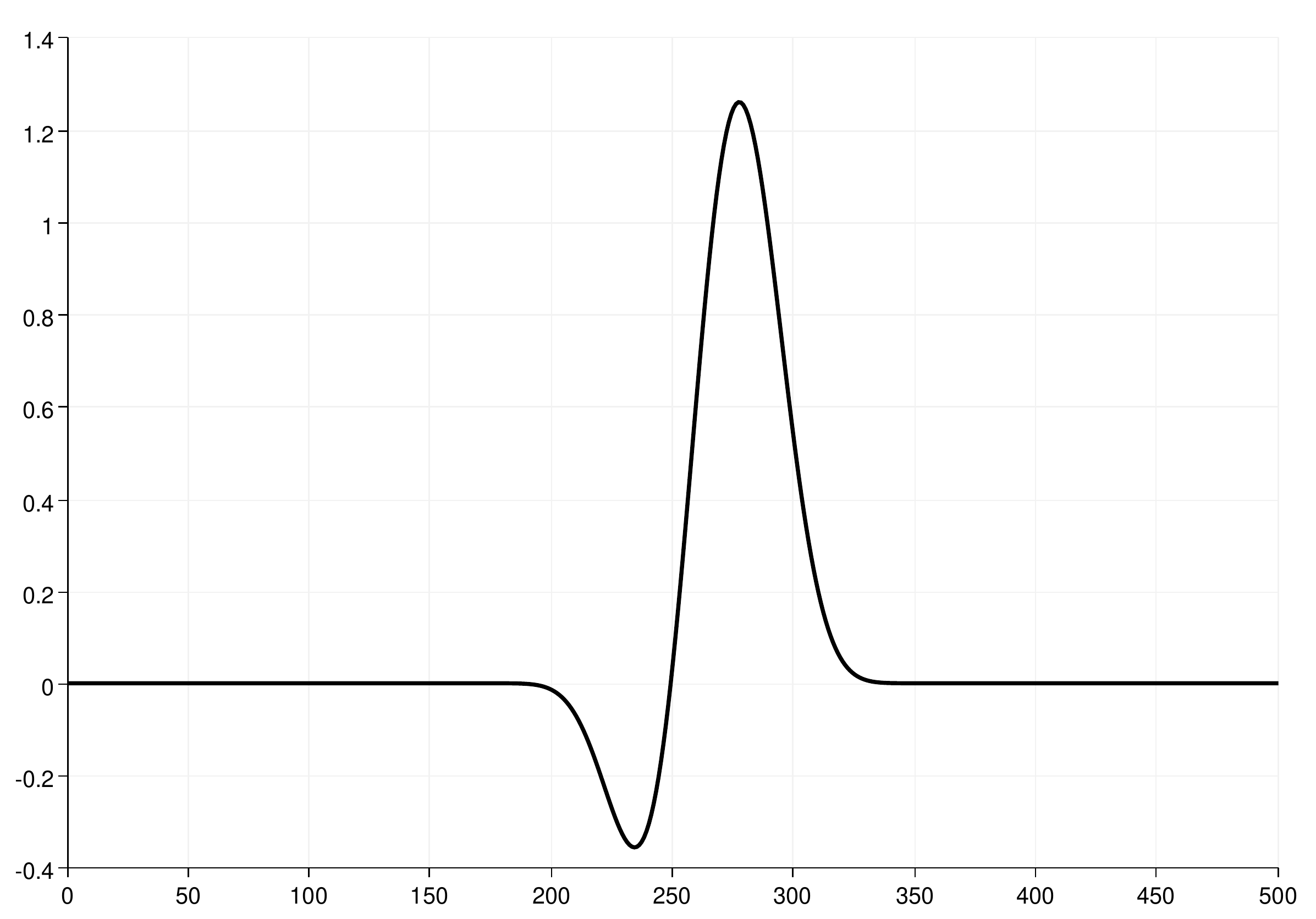} & \includegraphics[width=0.48\columnwidth]{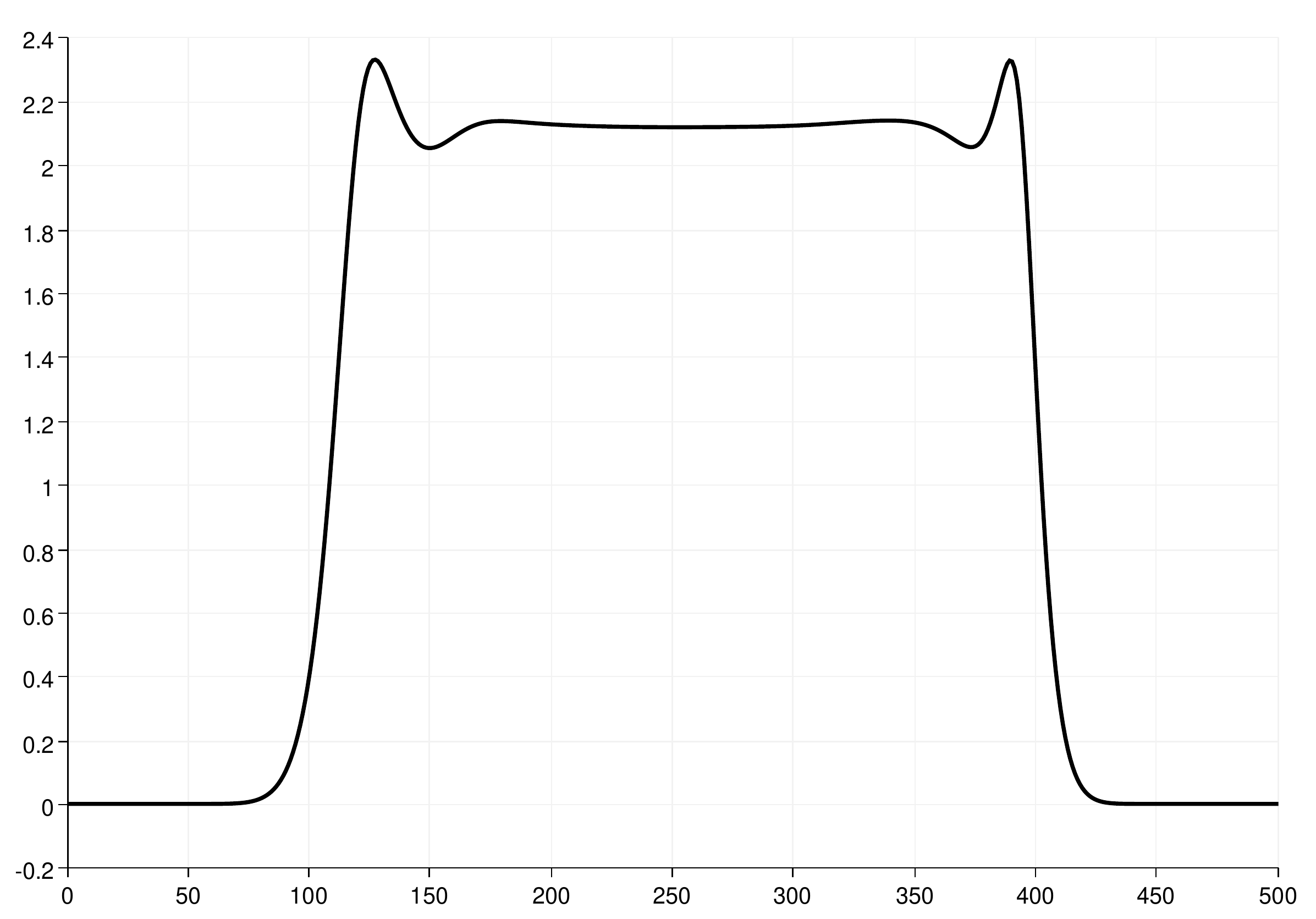}\tabularnewline
\end{tabular}
\end{figure}

\begin{figure}[H]
\caption{\label{fig:two-bubbles-init}Initial value along the diagonal for
interaction of two bubbles.}
\centering{}%
\begin{tabular}{cc}
 & \includegraphics[width=0.5\columnwidth]{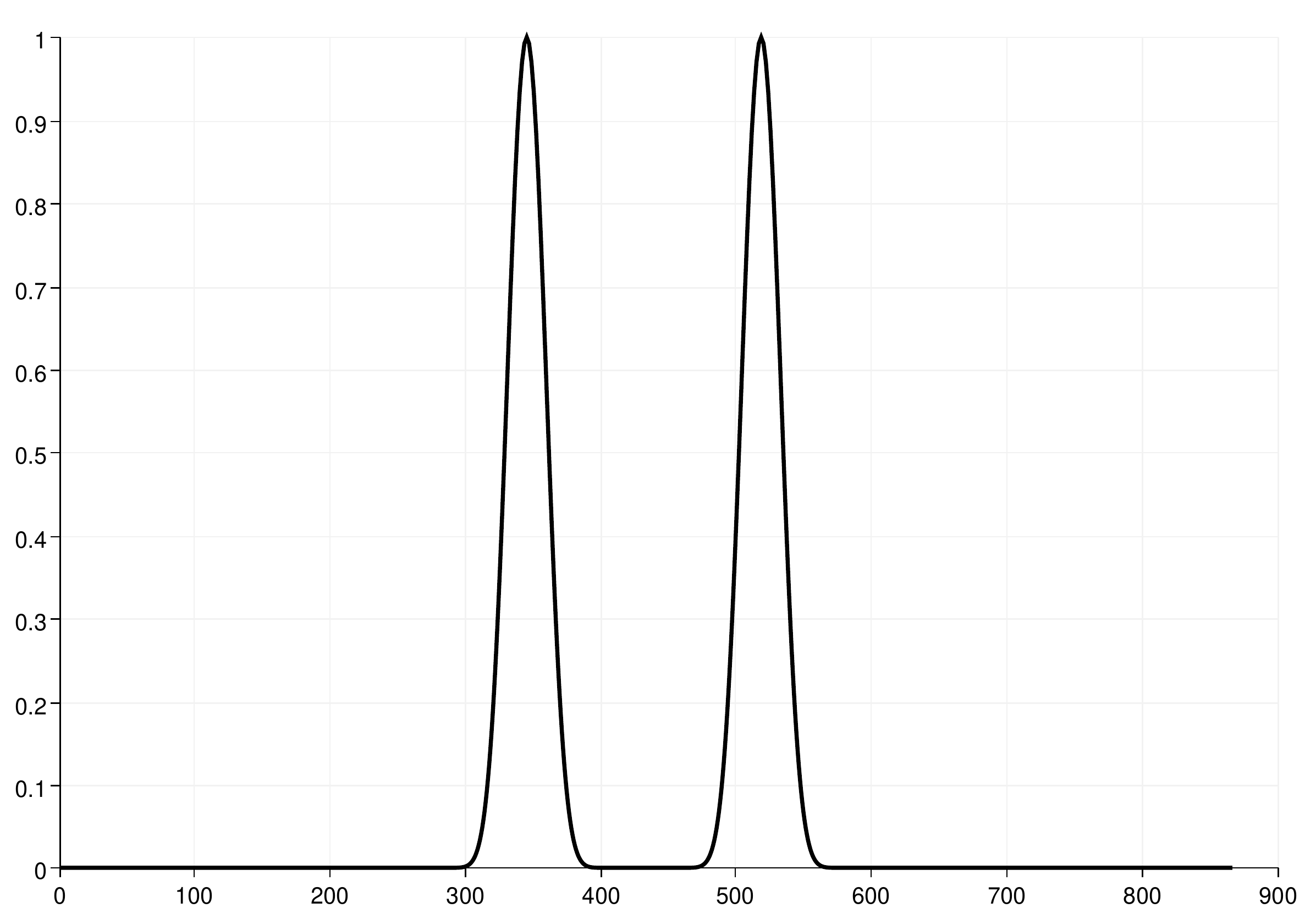}\tabularnewline
\end{tabular}
\end{figure}

\begin{figure}[H]
\caption{\label{fig:two-bubbles-2}Interaction of two bubbles I.}
\centering{}%
\begin{tabular}{cc}
(a.1) Solution along a line at $t=0.08$ & (a.2) 3D bubbles at $t=0.08$\tabularnewline
\includegraphics[width=0.48\columnwidth]{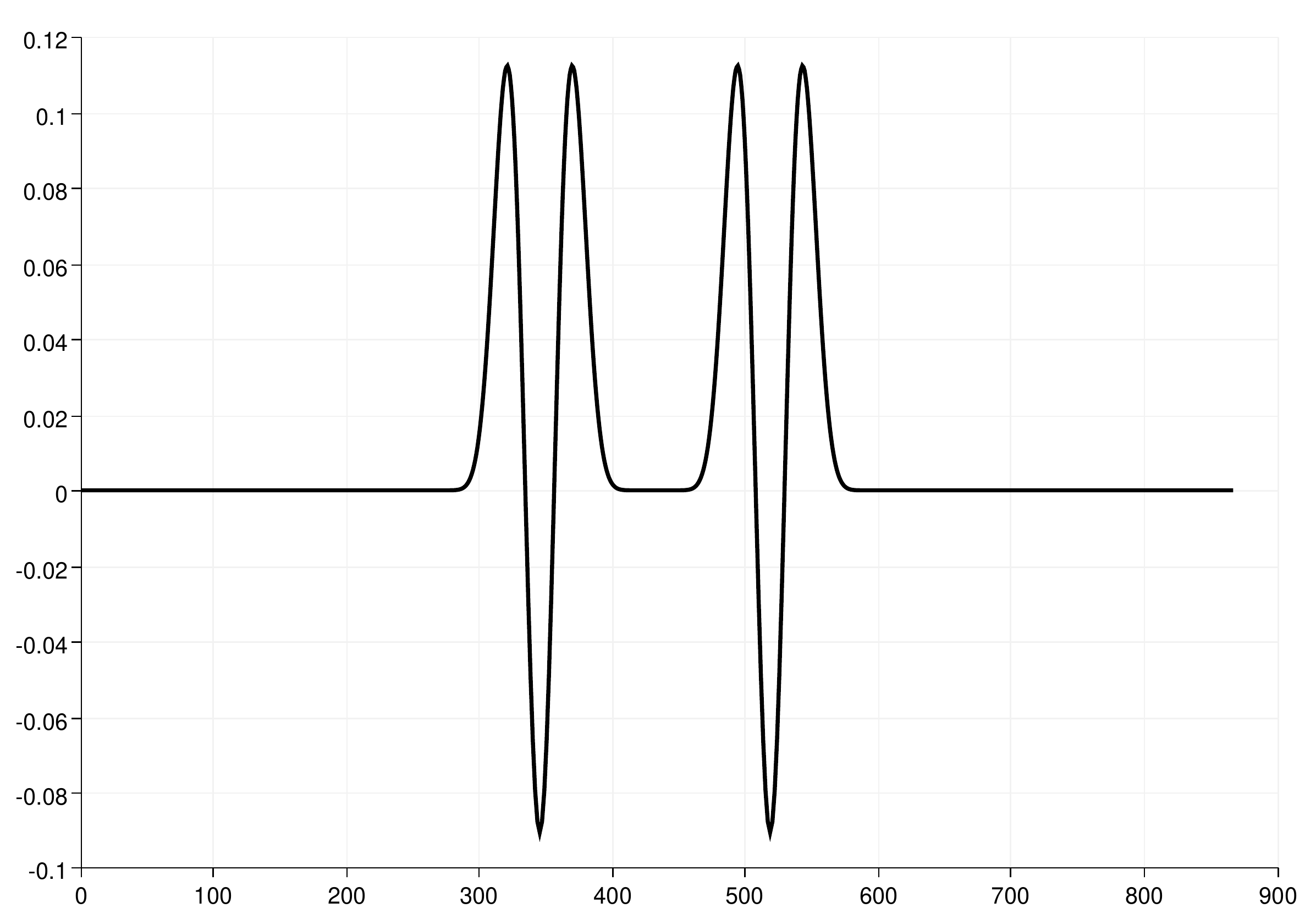} & \includegraphics[width=0.48\columnwidth]{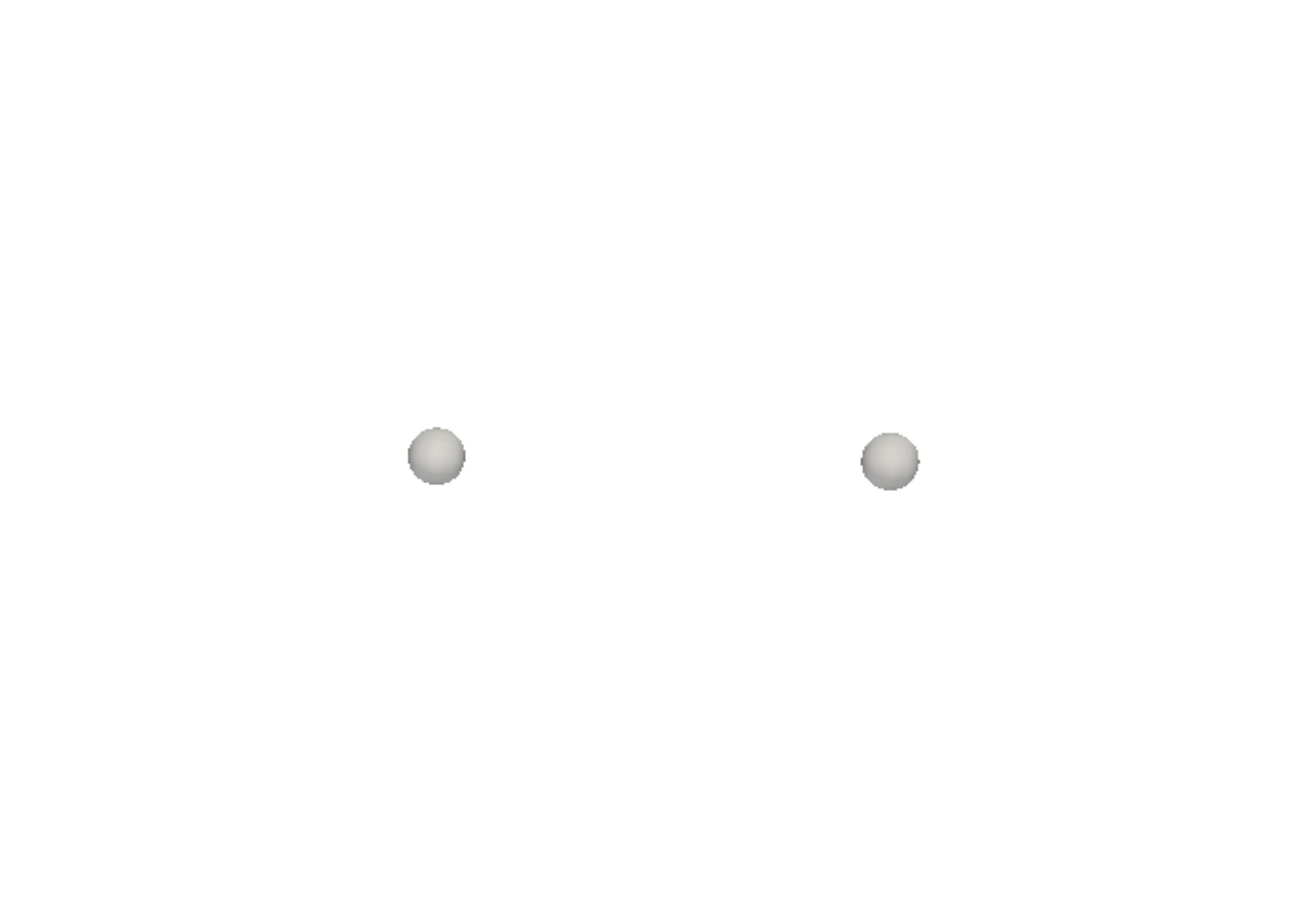}\tabularnewline
 & \tabularnewline
(b.1) Solution along a line at $t=0.2$ & (b.2) 3D bubbles at $t=0.2$\tabularnewline
\includegraphics[width=0.48\columnwidth]{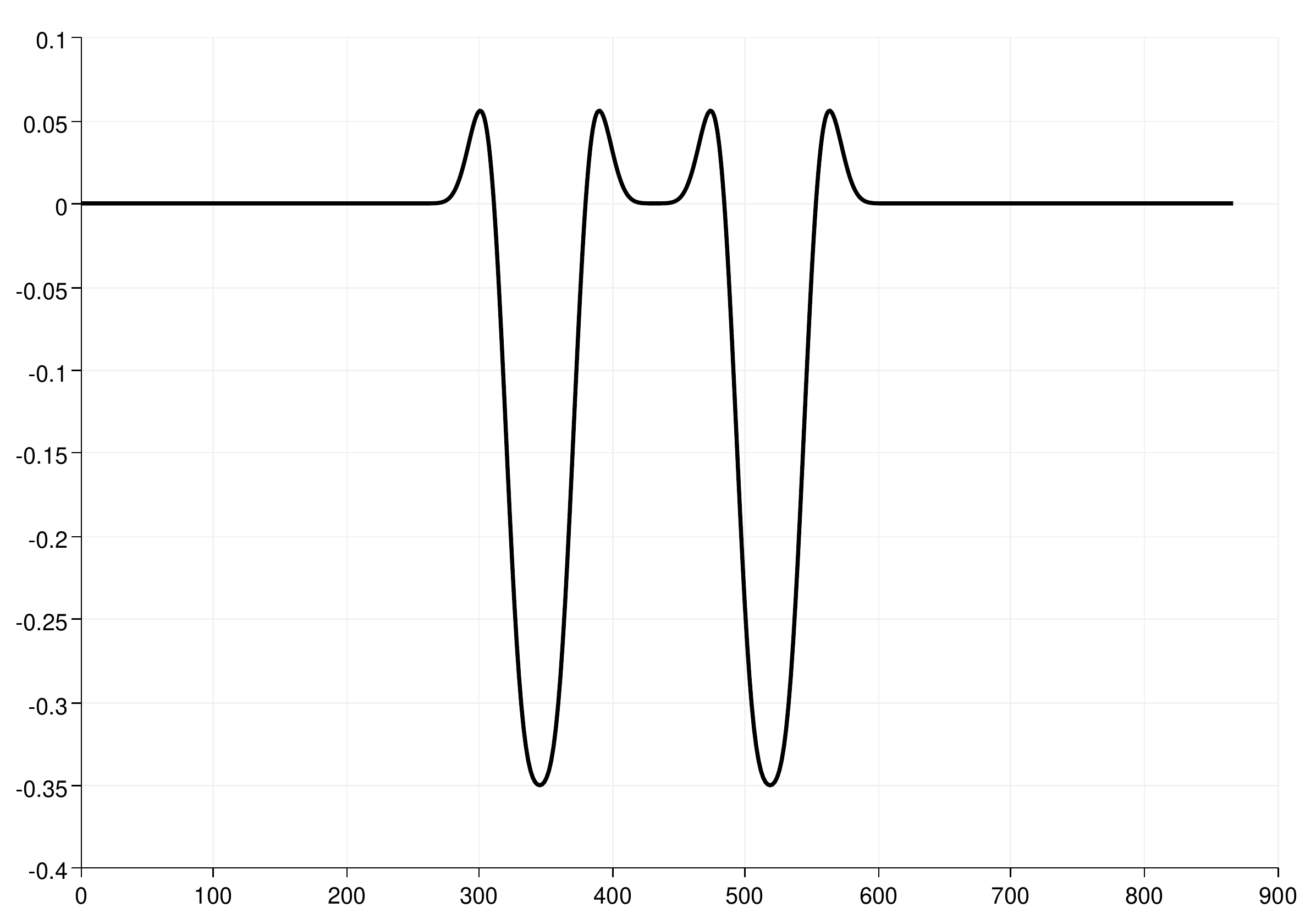} & \includegraphics[width=0.48\columnwidth]{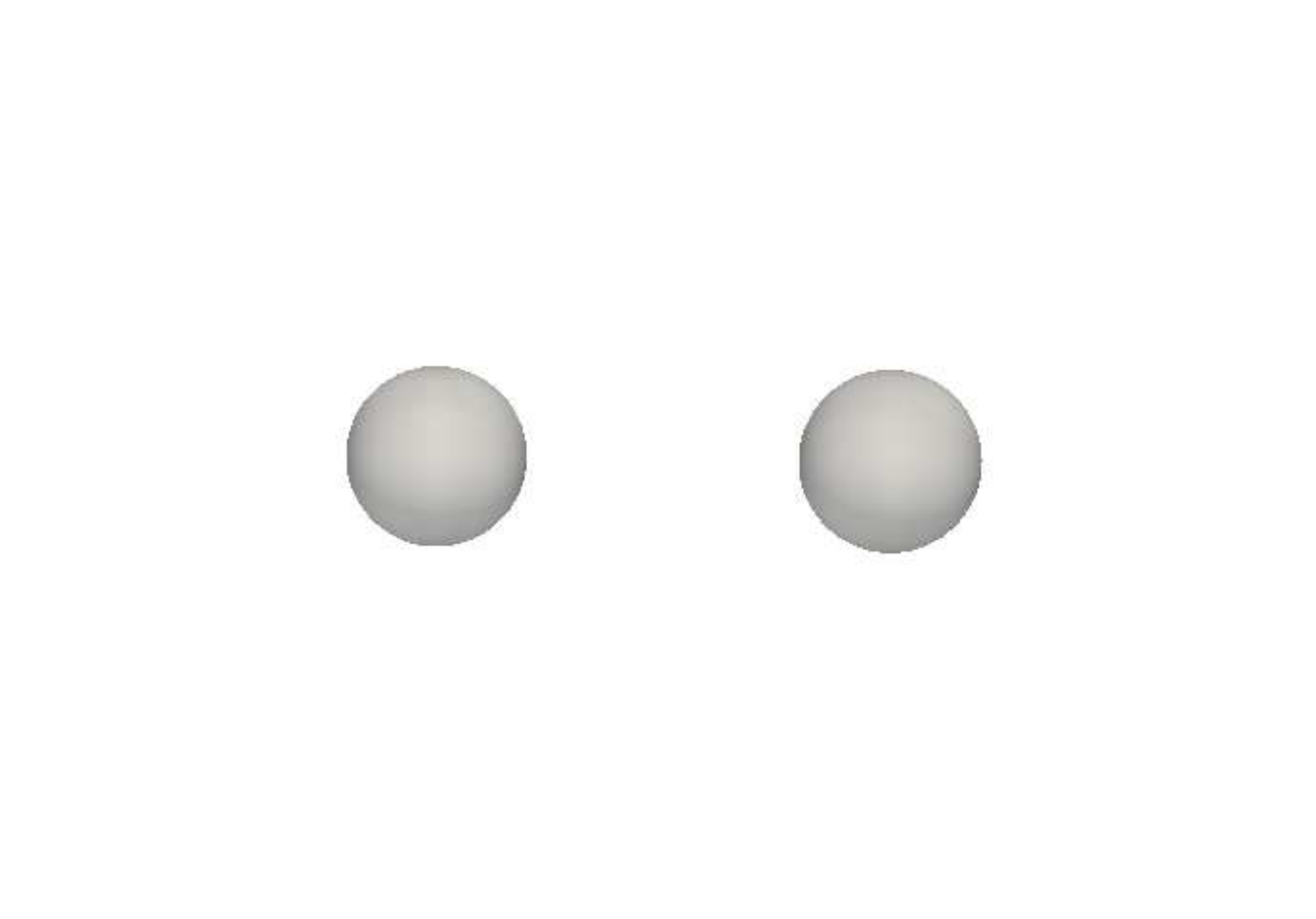}\tabularnewline
 & \tabularnewline
(c.1) Solution along a line at $t=0.69$ & (c.2) 3D bubbles at $t=0.69$\tabularnewline
\includegraphics[width=0.48\columnwidth]{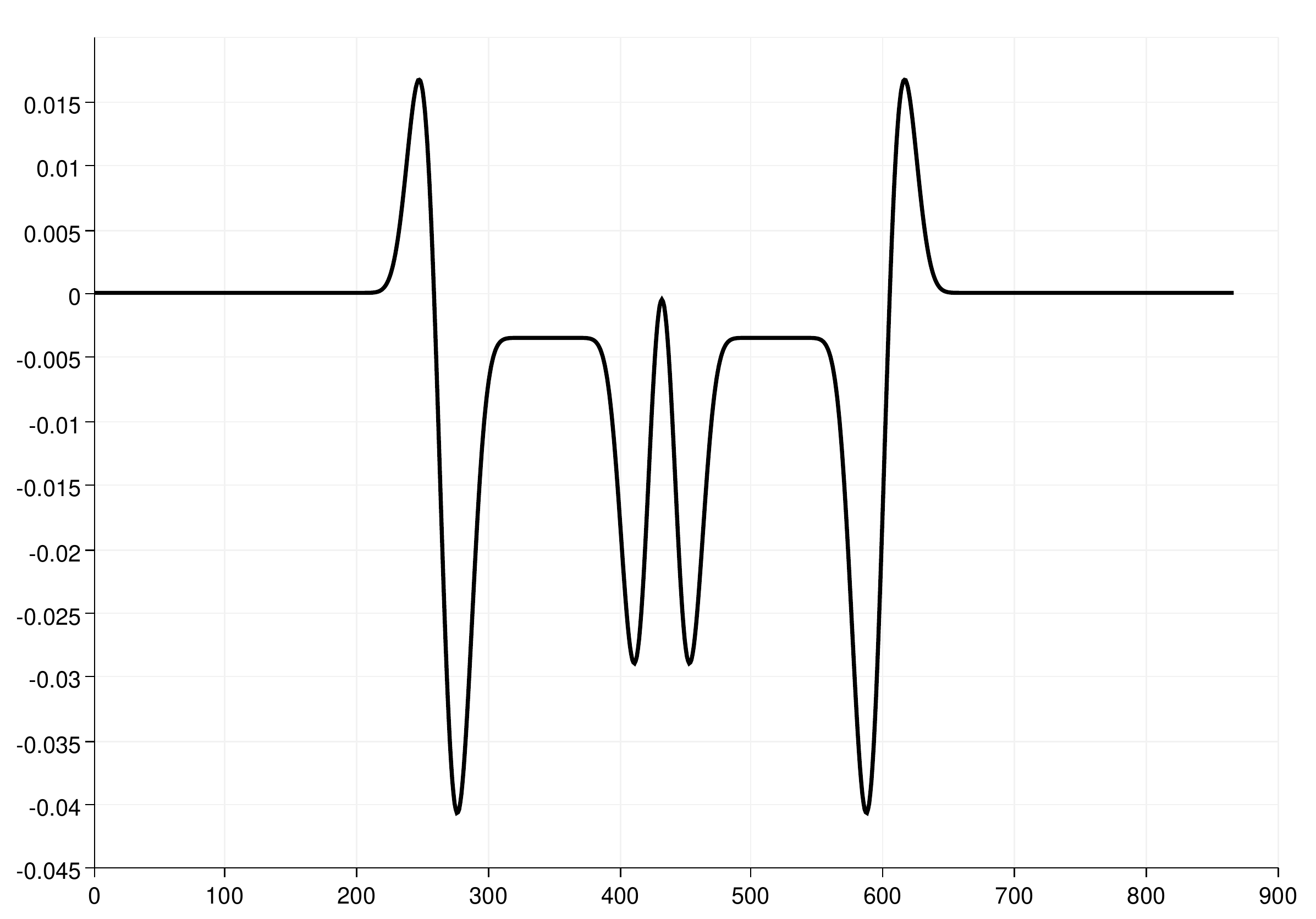} & \includegraphics[width=0.48\columnwidth]{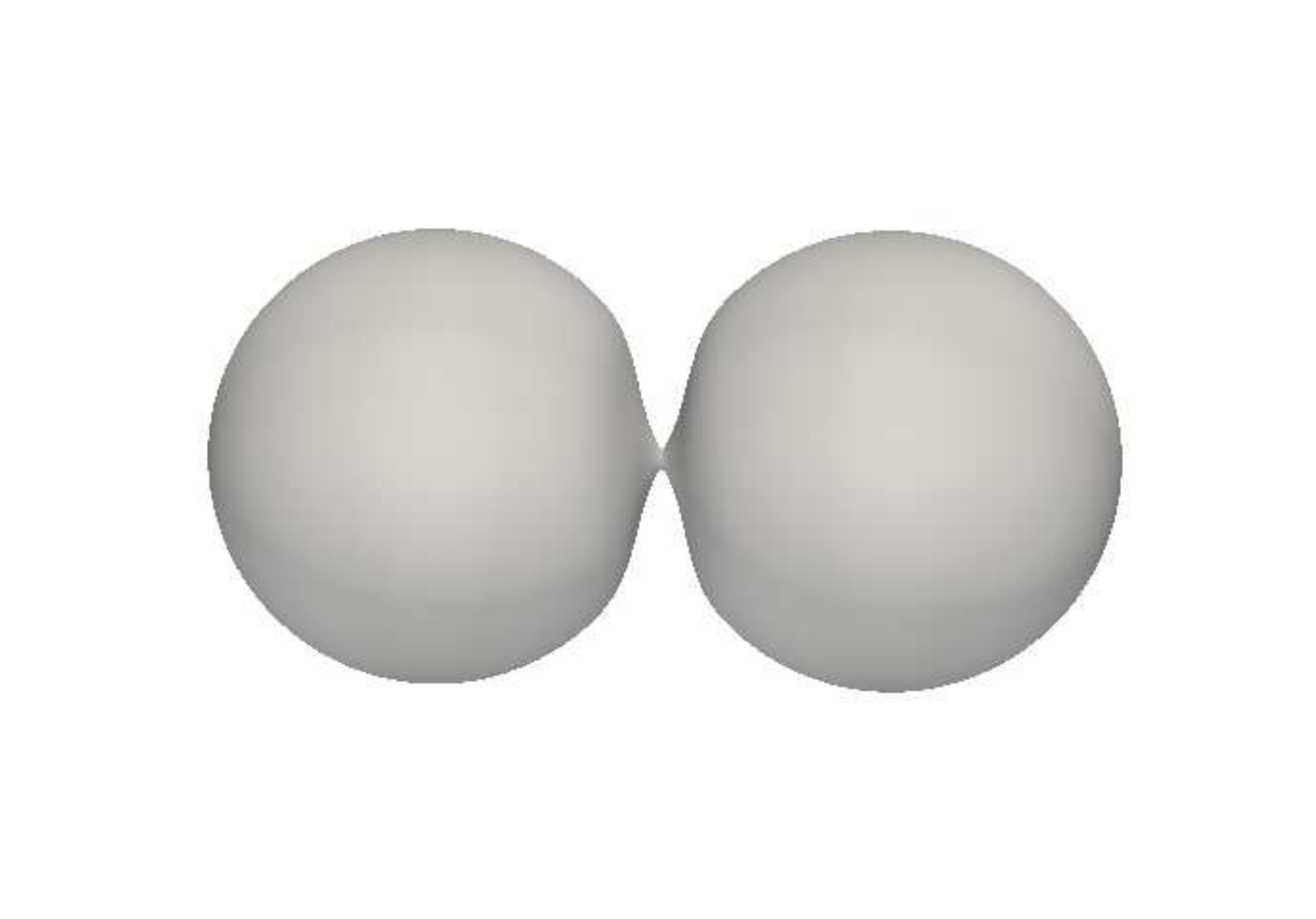}\tabularnewline
\end{tabular}
\end{figure}

\begin{figure}[H]
\caption{\label{fig:two-bubbles-3}Interaction of two bubbles II.}
\centering{}%
\begin{tabular}{cc}
(a.1) Solution along a line at $t=0.8$ & (a.2) 3D bubbles at $t=0.8$\tabularnewline
\includegraphics[width=0.48\columnwidth]{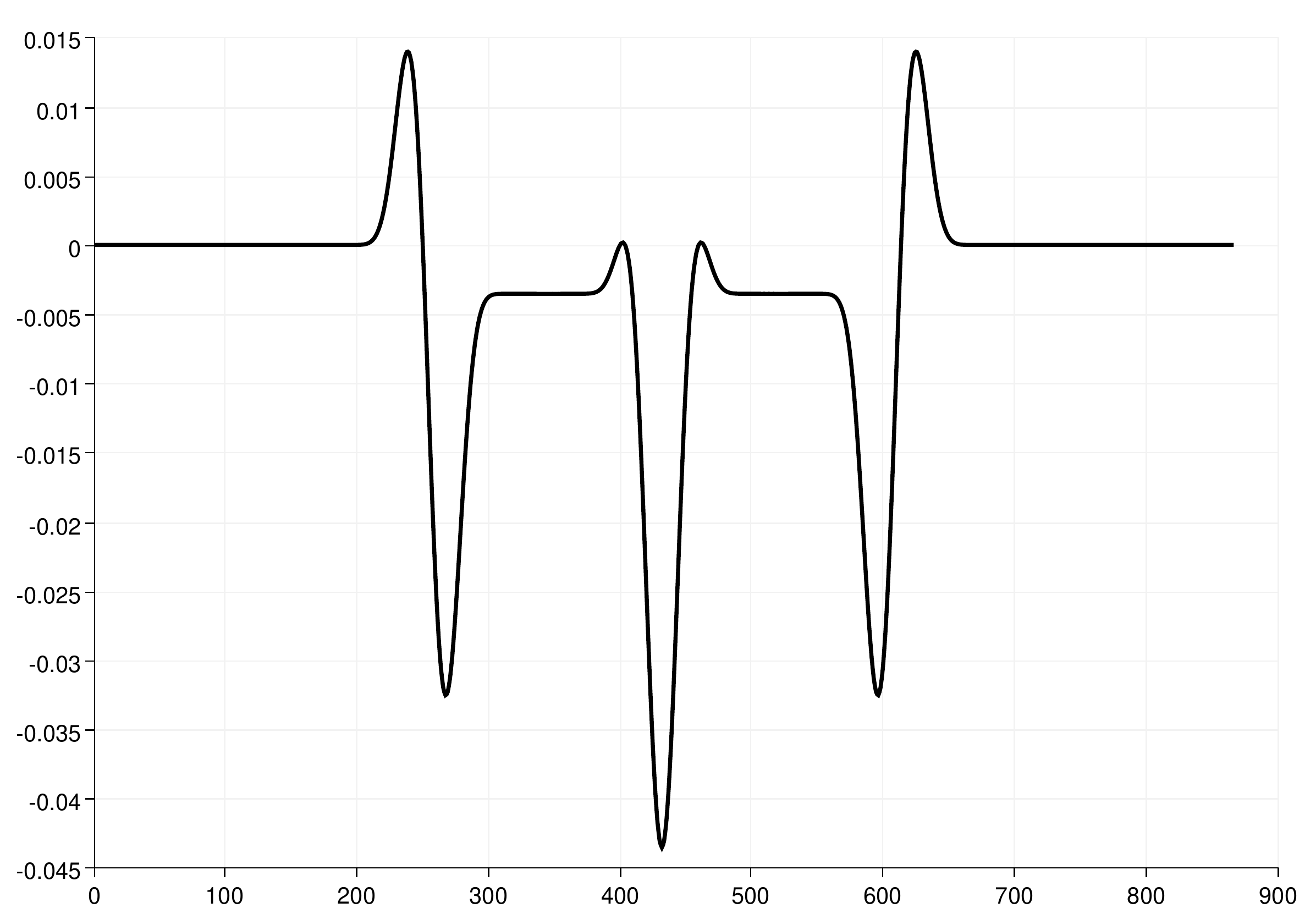} & \includegraphics[width=0.48\columnwidth]{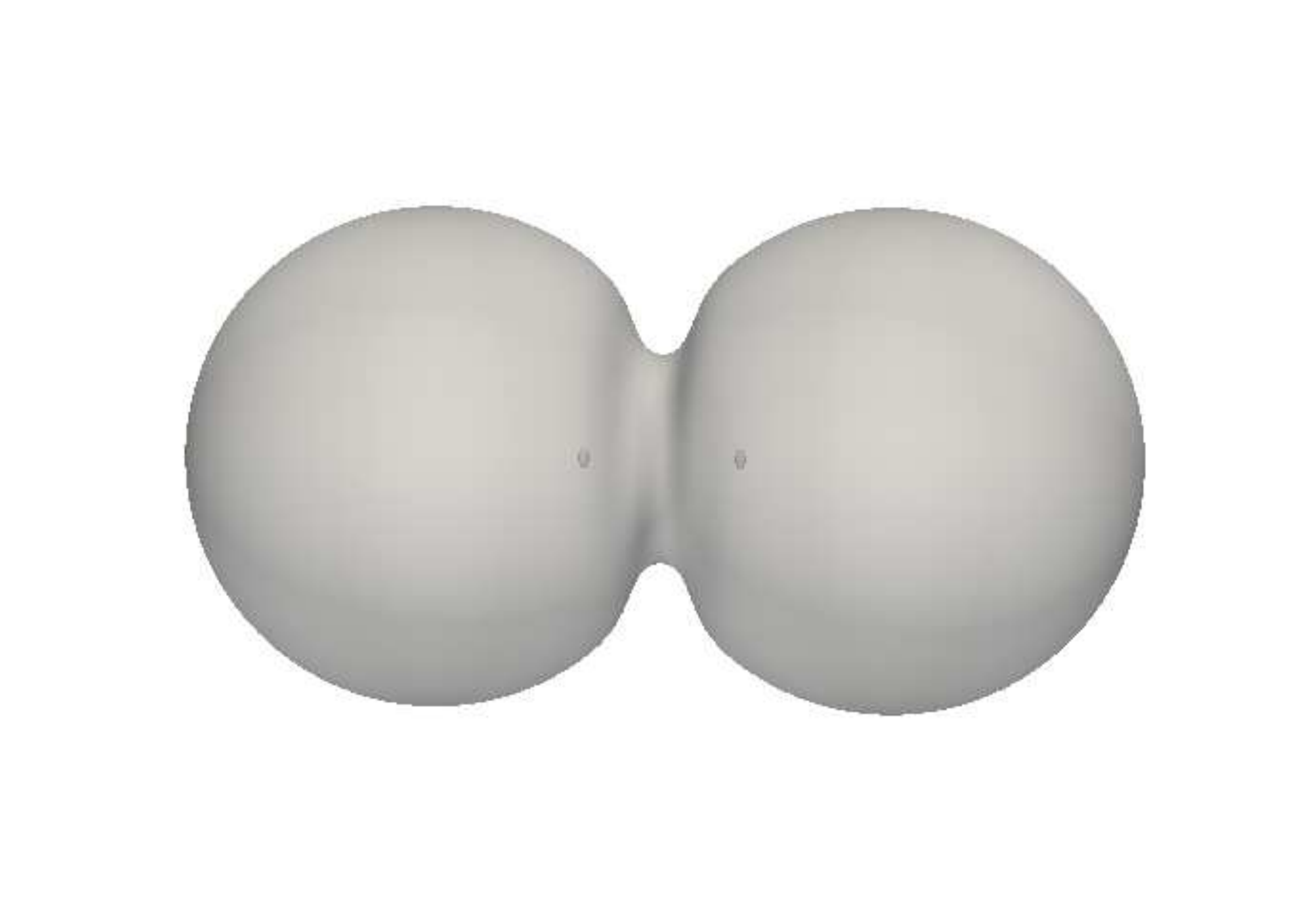}\tabularnewline
 & \tabularnewline
(b.1) Solution along a line at $t=1$ & (b.2) 3D bubbles at $t=1$\tabularnewline
\includegraphics[width=0.48\columnwidth]{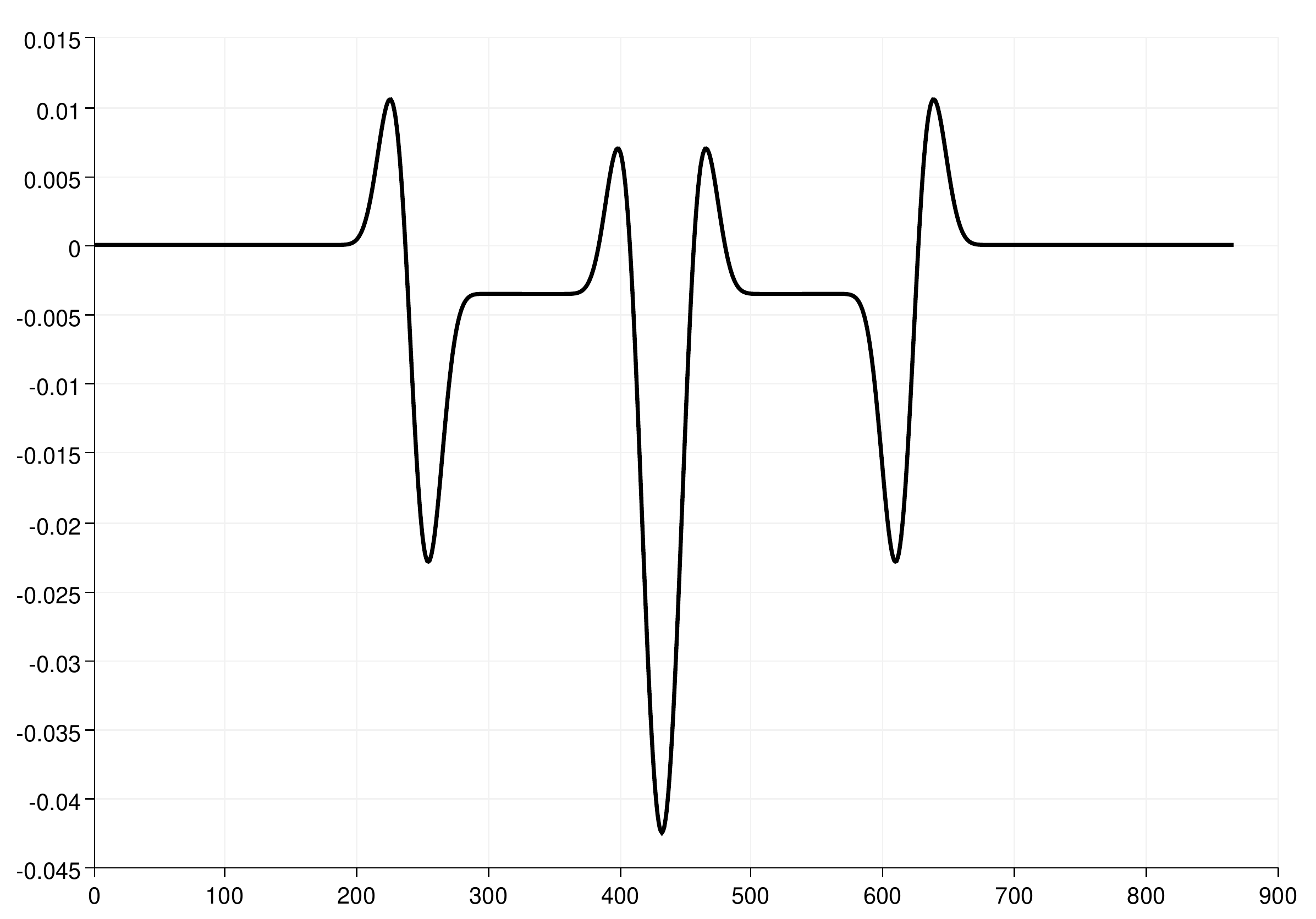} & \includegraphics[width=0.48\columnwidth]{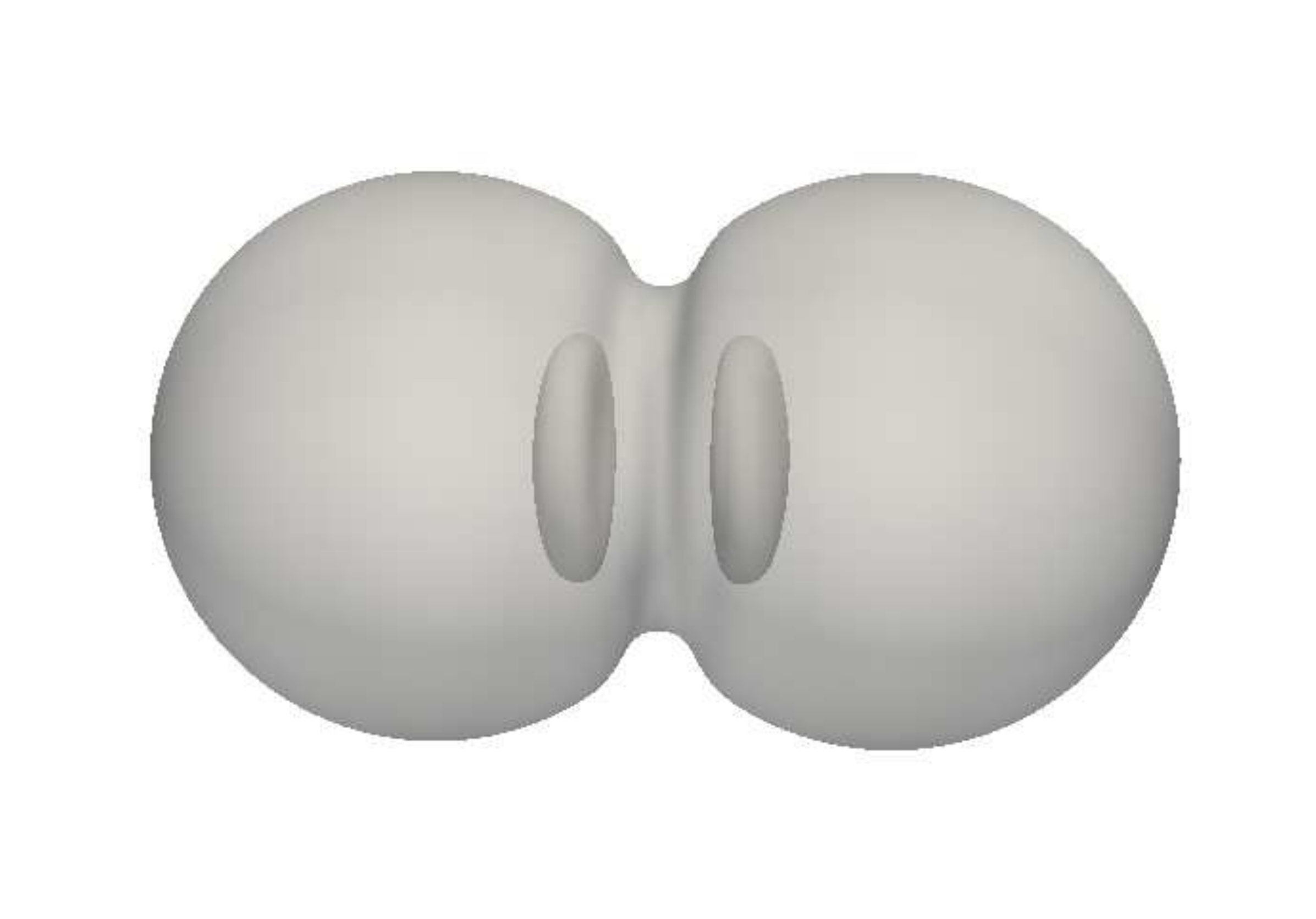}\tabularnewline
 & \tabularnewline
(c.1) Solution along a line at $t=2$ & (c.2) 3D bubbles at $t=2$\tabularnewline
\includegraphics[width=0.48\columnwidth]{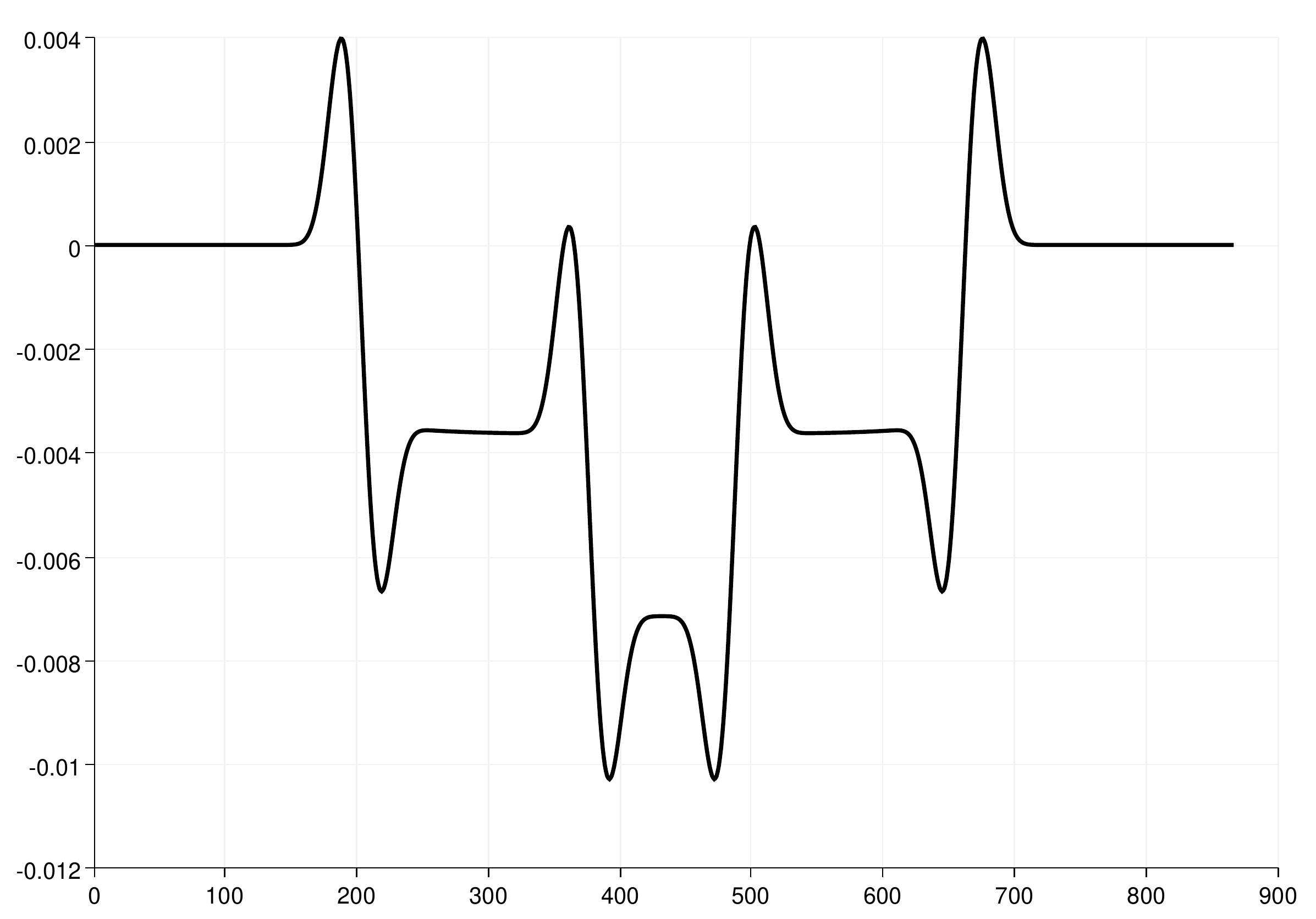} & \includegraphics[width=0.48\columnwidth]{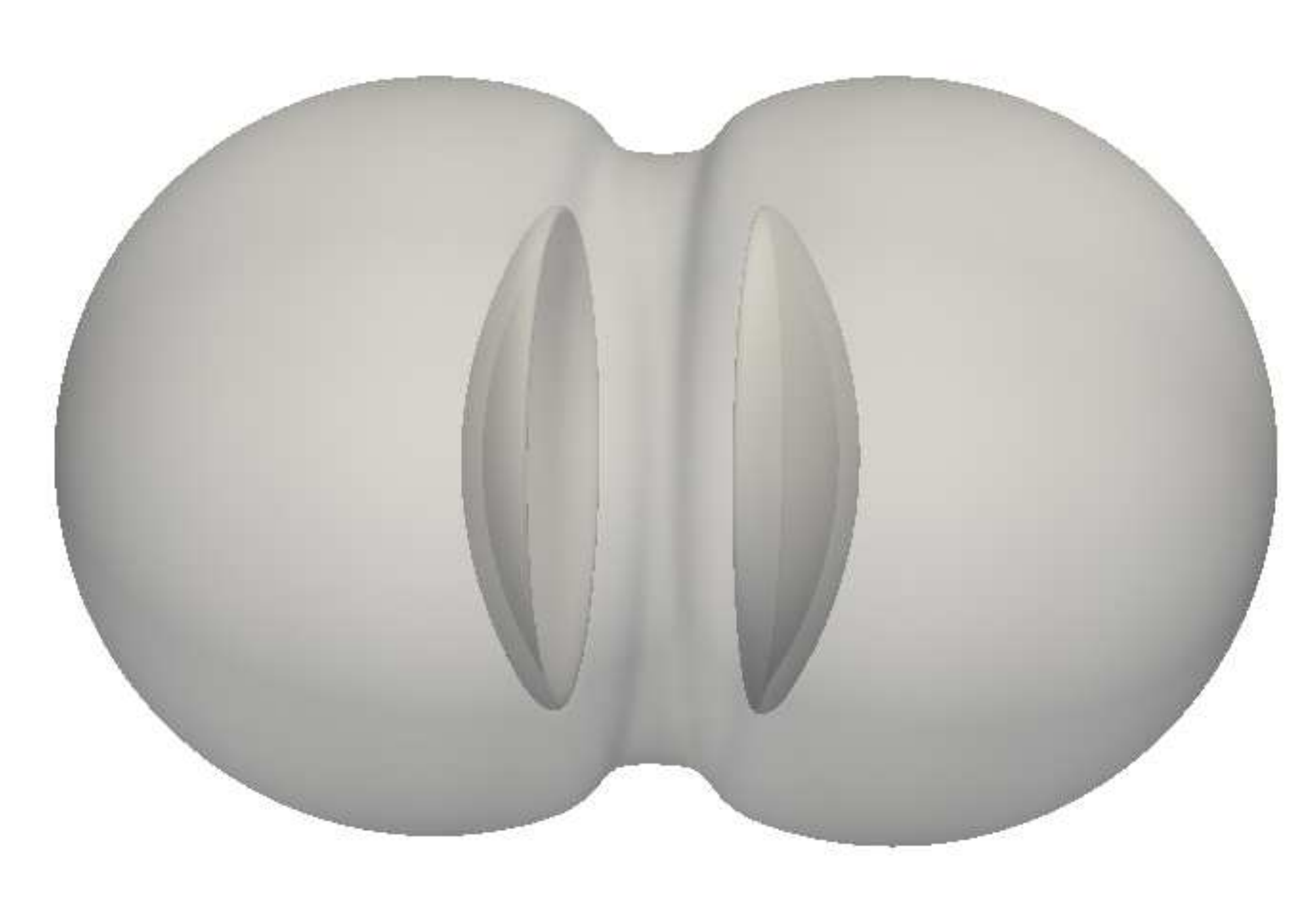}\tabularnewline
\end{tabular}
\end{figure}

\begin{figure}[H]
\caption{\label{fig:two-bubbles-4}Interaction of two bubbles III.}
\centering{}%
\begin{tabular}{cc}
(a.1) Solution along a line at $t=2.15$ & (a.2) 3D bubbles at $t=2.15$\tabularnewline
\includegraphics[width=0.48\columnwidth]{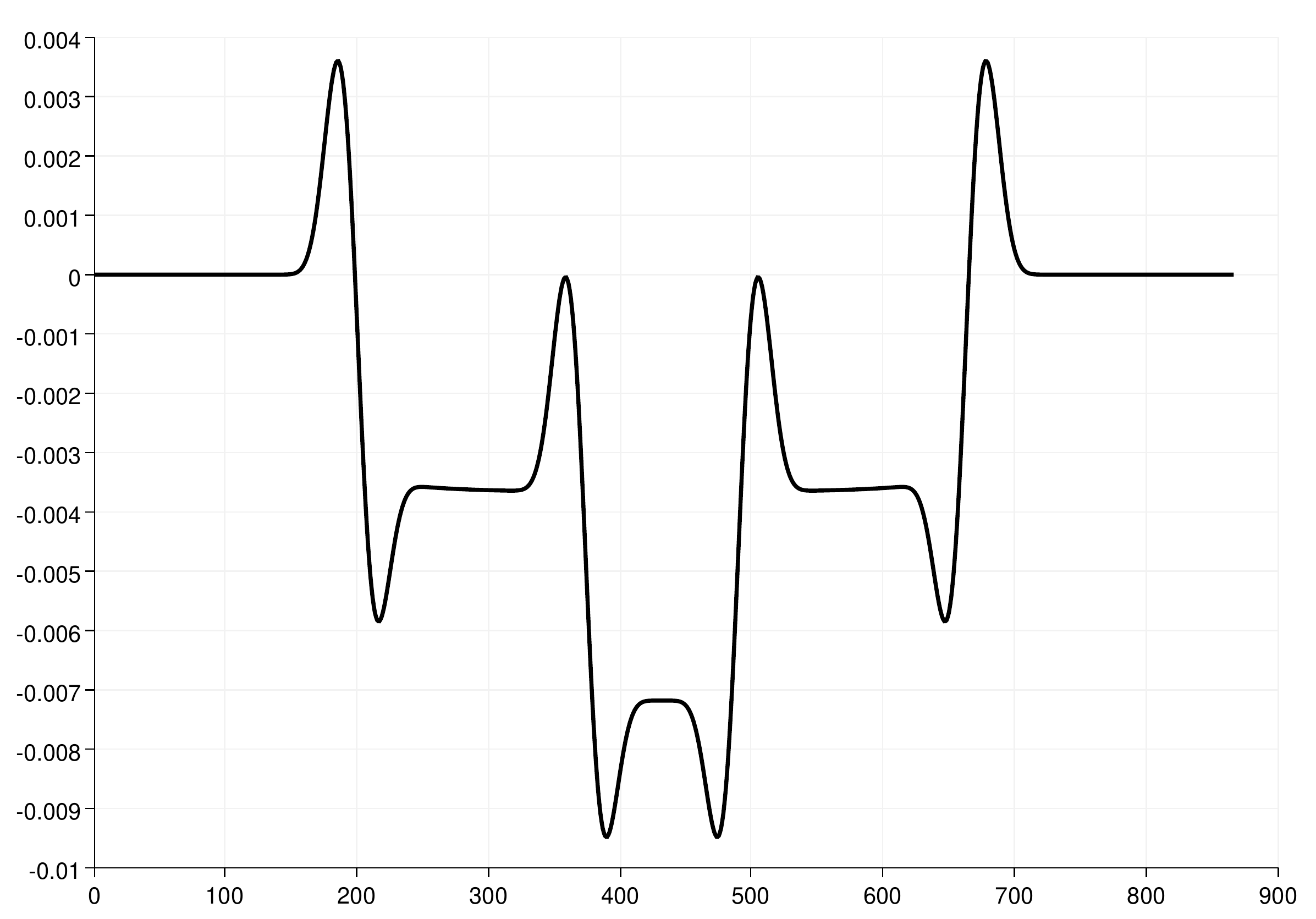} & \includegraphics[width=0.48\columnwidth]{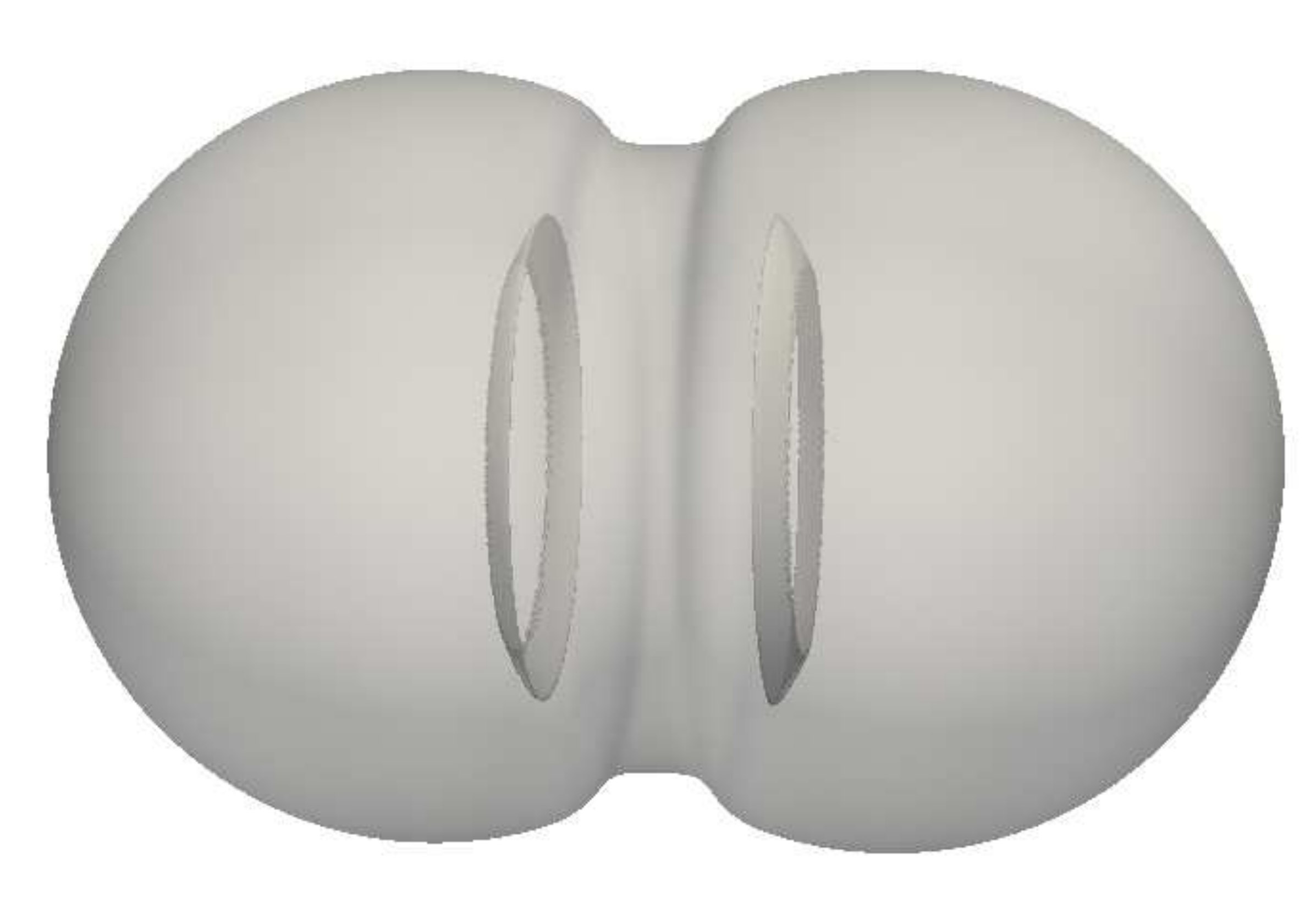}\tabularnewline
 & \tabularnewline
(b.1) Solution along a line at $t=3$ & (b.2) 3D bubbles at $t=3$\tabularnewline
\includegraphics[width=0.48\columnwidth]{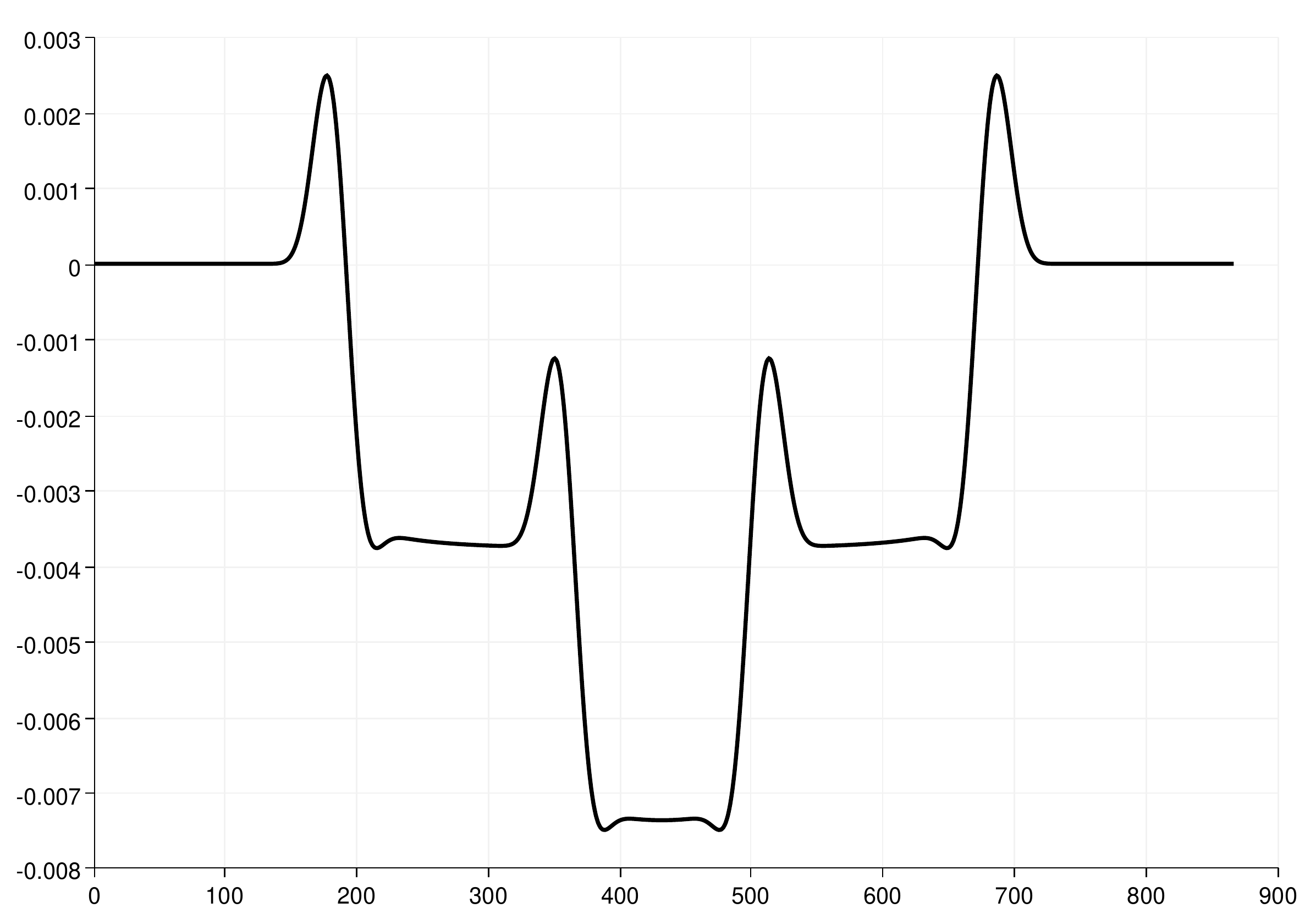} & \includegraphics[width=0.48\columnwidth]{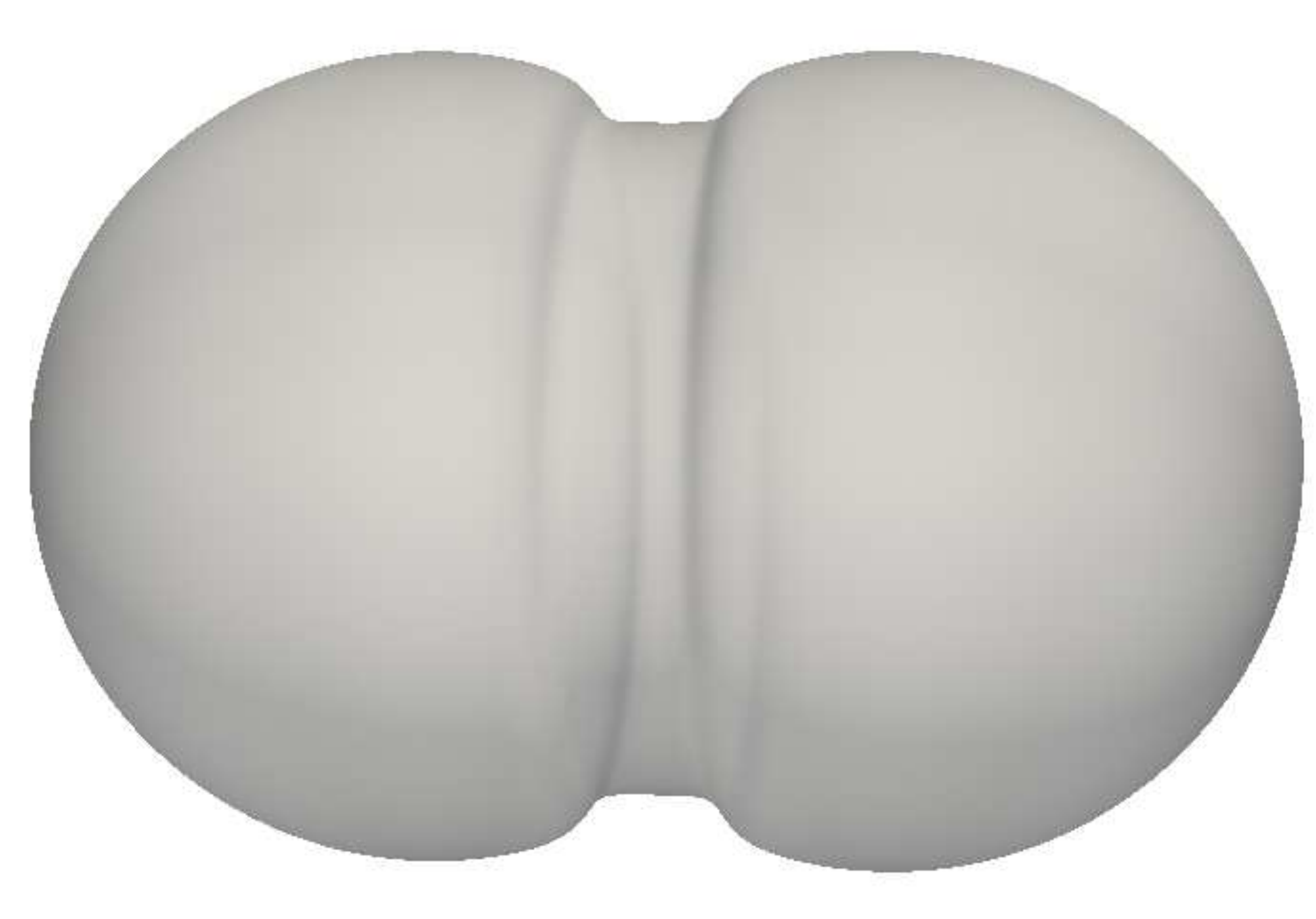}\tabularnewline
 & \tabularnewline
(c.1) Solution along a line at $t=50$ & (c.2) 3D bubbles at $t=50$\tabularnewline
\includegraphics[width=0.48\columnwidth]{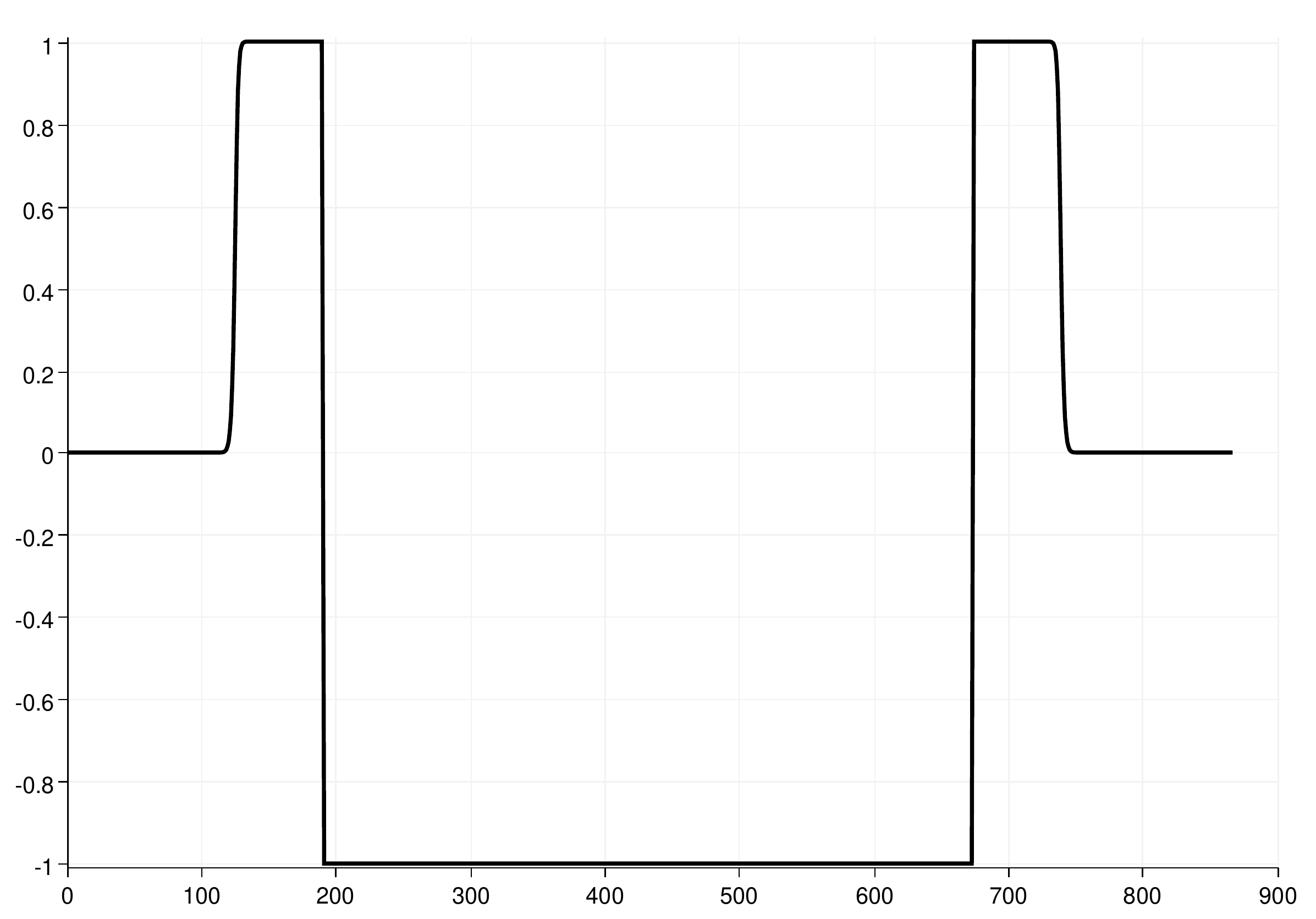} & \includegraphics[width=0.48\columnwidth]{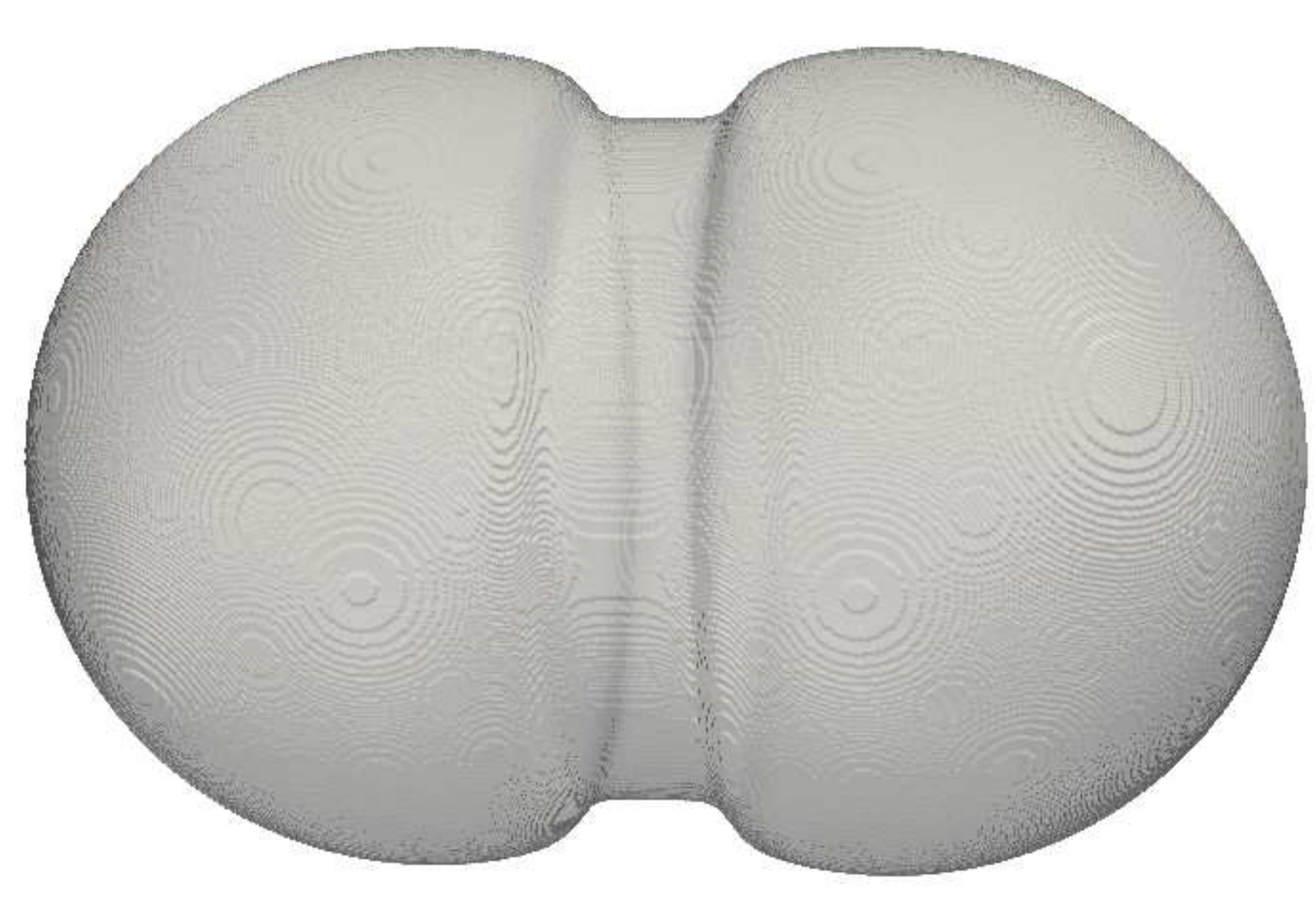}\tabularnewline
\end{tabular}
\end{figure}

In this section we present several examples for different parameter
values of $\mu$, $\lambda$ and for different initial conditions
$\varphi_{0}$ and $\varphi_{1}$. Unless otherwise noted in the examples
the grid size in space is $n\times n\times n=500\times500\times500$,
resulting in uniform grid spacing $\delta x=\delta x=\delta x=2\times10^{-3}$.
With $\delta t=10^{-4}$ the Courant\textendash Friedrichs-{}-Lewy
(CFL) condition ${\displaystyle \left|\phi\right|<\frac{\delta x}{\sqrt{3}\delta t}\approx11.54}$
for stability (see, e.g., \cite{Strang-Comp_Sci}) has been satisfied
for all times and for each of our simulation examples for the Higgs
boson equation. For the scaling factor the value $L=5$ has been used.
For visualization of the solution $\phi\left(\vec{x},t\right)$ we
use line plots along either the line segment connecting the midpoints
of the computational cube's faces parallel to the $yz$\textendash plane
(mid-line parallel to the $x$-axis) or the diagonal line segment
connecting the corners $\left(0,0,0\right)$ and $\left(1,1,1\right)$
of the unit cube (see Figure \ref{fig:Line-segments}). The horizontal
axis shows ranges $\left[0,500\right]$ and $\left[0,900\right]$
respectively on these line plots due to grid the size $n=500$ and
due to the diagonal line segment's length being $500\sqrt{3}\approx866\approx900$.
For initial conditions we mostly use variations of the compactly supported,
infinitely smooth bump function (often used as test functions in the
theory of generalized functions (see, e.g., \cite{FRY2002143}))
\begin{equation}
{\displaystyle B\left(\vec{x};C,R\right)=\begin{cases}
\exp\left(\frac{1}{R^{2}}-\frac{1}{R^{2}-\left|\vec{x}-C\right|^{2}}\right) & \text{if }\left|\vec{x}-C\right|<R,\\
0 & \text{if }\left|\vec{x}-C\right|\geq R
\end{cases}}\label{eq:bump}
\end{equation}
with center $C=\left(C_{1},C_{2},C_{3}\right)\in\Omega$ and with
radius $0<R\ll1$. Here $\left|\vec{x}-C\right|$ denotes the euclidean
distance between points $\vec{x}$ and $C$. Note that these basic
bump functions $B\left(\vec{x};C,R\right)$ are nonnegative with maximum
value $1$ (see Figure \ref{fig:Bump-Function}). 
\begin{example}
\label{exa:Convergence} In order to examine the convergence properties
of our numerical method we look at an example of Equation \eqref{eq:higgs-boson-rescalled}
with parameter values $\mu^{2}=9$, $\lambda=2$ and with the initial
conditions. 
\begin{align}
{\displaystyle \varphi_{0}\left(\vec{x}\right)} & =3B\left(x;\left(0.5,0.5,0.5\right),0.3\right),\qquad\forall\vec{x}\in\Omega\label{eq:init0-conv}\\
\varphi_{1}\left(\vec{x}\right) & =0,\qquad\forall\vec{x}\in\Omega.
\end{align}
This case will be discussed in more details in Example \ref{exa:no-bubble},
right now we only look at how solution changes with respect to different
grid sizes and different precisions. We vary the grid size in space
and compare the results for values $n=200$, $300$, and $400$ to
the result with $n=500$ at times $t=1$ and $t=2$. Since we are
working with compactly supported radial solution starting from a bump
function, the solution is nonzero only around the center of the computational
domain, and we only consider line plots of the pointwise difference
instead of $L^{2}$, or other global error norms. Figure \ref{fig:Convergence}
shows stability as the difference between the numerical solutions
decreases with increasing values of $n$. The difference is the largest
in the regions of high-slope, which moves outward from the center,
suggesting that the difference in dispersion plays a role, as suggested
for example in \cite{cohen-fourth-order} for the acoustic wave equation.
Figure \ref{fig:Precision} shows the difference between running our
code in single and double precision at times $t=1$ and $t=2$. While
the difference is insignificant, and does not seem to increase with
time, the double precision code required about $80\%$ more run-time
than the single precision code. 
\end{example}
\begin{example}
\label{exa:blowup}As a second computational example we demonstrate
that our numerical code can predict the blow up of the sign-preserving
solution for not the Higgs boson equation but for the related Klein-Gordon
equation with an imaginary mass (see Theorem 1.1 in \cite{yagdjian_DCDS}).
For this purpose we set the parameter values as $\mu^{2}=1$ and $\lambda=-1$.
As initial conditions we use the bump functions 
\begin{align}
\varphi_{0}\left(\vec{x}\right) & =2B\left(\vec{x};\left(0.5,0.5,0.5\right),0.2\right),\qquad\forall\vec{x}\in\Omega,\\
\varphi_{1}\left(\vec{x}\right) & =10B\left(\vec{x};\left(0.5,0.5,0.5\right),0.2\right),\qquad\forall\vec{x}\in\Omega.
\end{align}
Figure \ref{fig:Blowup}, parts (a)-(e) show the line plot of the
solution $\phi\left(x,t\right)$ for various times. In particularly,
we may notice from parts (a) and (b) of Figure \ref{fig:Blowup},
that from time $t=0$ to time $t=2$ the magnitude of the solution
is around $2$ and it does not increase, but the compact support of
the solution becomes larger instead. After time $t=2$ the magnitude
of solution suddenly starts to increase and it reaches the value of
approximately $130$ by time $t=2.97$ (part (e) of Figure \ref{fig:Blowup}).
The blow up of the integral ${\displaystyle \int_{\Omega}\phi\left(\vec{x},t\right)d\vec{x}}$
shortly after time $t=2.97$ can be observed on part (f) of Figure
\ref{fig:Blowup}, demonstrating the theoretical results of Theorem
1.1 in \cite{yagdjian_DCDS}. 
\end{example}
\begin{example}
\label{example:short-bubble} In this example we demonstrate that
under certain conditions the solution forms bubbles. For this purpose
we reformulate the theoretical results of Corollary 1.4 in \cite[page 453]{CPDE2012}:
\begin{cor}
\cite[Corollary 1.4]{CPDE2012} Bubble forms if the initial data satisfy
\begin{align}
\left(\frac{n}{2}+\sqrt{\frac{n^{2}}{4}+\mu^{2}}\right)\varphi_{0}\left(\vec{x}\right)+\varphi_{1}\left(\vec{x}\right) & <0\quad\forall\vec{x}\in\Omega,\label{eq:bubble-cond1}\\
\int_{\mathbb{R}^{3}}\phi^{3}\left(\vec{x},t\right)d\vec{x} & \geq0\quad\forall t\geq0.\label{eq:bubble-cond2}
\end{align}
\end{cor}
In particular, we consider as initial data the bump functions 
\begin{align}
{\displaystyle \varphi_{0}\left(\vec{x}\right)} & =B\left(\vec{x};\left(0.5,0.5,0.5\right),0.3\right),\qquad\forall\vec{x}\in\Omega,\label{eq:mollifier}\\
\varphi_{1}\left(\vec{x}\right) & =-5B\left(\vec{x};\left(0.5,0.5,0.5\right),0.3\right),\qquad\forall\vec{x}\in\Omega.
\end{align}
The line plot of the initial function \eqref{eq:mollifier} is depicted
in part (a) of Figure \ref{fig:Bubble-Line}. Note that the shift
$\vec{x}_{1}=\vec{x}-C$ makes the initial data and hence the solution
radial. The parameter values are $\mu^{2}=9$, $\lambda=2$, and hence
condition \eqref{eq:bubble-cond1} is satisfied inside the domain
of support. Condition \eqref{eq:bubble-cond2} is also satisfied for
$t=0$ and, by continuity, at least for small $t\geq0$. Part (b)
of Figure \ref{fig:Bubble-Line} shows that at time $t=0.21$ the
values of solution are still all positive inside the domain of support.
In parts (c), (d), and (e) of Figure \ref{fig:Bubble-Line} at time
instances $t=0.22$, $t=0.23$, and $t=0.4$ we can see two points,
where the solution has zero values with sign change inside the domain
of support. These zero places move from the center towards the border
of the domain as time goes by. In thee dimensions the isosurface corresponding
to these zero values has the shape of a ball (hence the name bubble)
centered at the point $C=\left(0.5,0.5,0.5\right)$, as shown in part
(f) of Figure \ref{fig:Bubble-Line}. The radius of the ball increases
in time, with the rate of increase decreasing exponentially following
the finite speed of propagation ${\displaystyle \frac{e^{-2t}}{L^{2}}}$. 
\end{example}
\begin{example}
\label{example:long-bubble} We go back now to Example \ref{example:short-bubble}
and continue it for larger time in order to investigate the long-time
behavior of the solution. Figure \ref{fig:Bubble-long-time} shows
line plots of the solution for various times. The solution converges
to a piecewise constant (step) function with values $0$, $2.12$,
and $-2.12$. The step function's shape with sharp corners raises
the possibility that the solution looses smoothness. In order to examine
the question of differentiability we look at how the quantity $\mathcal{P}\left(t\right)\equiv{\displaystyle \frac{1}{L^{2}}e^{-2t}\max_{\vec{x}\in\Omega}\left|\Delta\phi\left(\vec{x},t\right)\right|}$
(related to the second term from equation \eqref{eq:higgs-boson-rescalled})
changes in time. Figure \ref{fig:Laplacian} shows that the laplacian
does not blow up in finite time and hence the solution remains smooth
at least up to the second order derivative. Moreover, since $\mathcal{P}\left(t\right)\equiv{\displaystyle \frac{1}{L^{2}}e^{-2t}\max_{\vec{x}\in\Omega}\left|\Delta\phi\left(\vec{x},t\right)\right|}$
converges fast to zero as the time $t$ increases, we obtain that
equation \eqref{eq:higgs-boson-rescalled} converges to the unforced,
damped Duffing equation 
\begin{align}
\phi_{tt}+3\phi_{t} & =\mu^{2}\phi-\lambda\phi{}^{3},\quad\vec{x}\in\Omega,\quad t\in\left(0,T\right],\label{eq:duffing}\\
\phi\left(\vec{x},0\right) & =\varphi_{0}\left(\vec{x}\right),\quad\vec{x}\in\Omega,\\
\phi_{t}\left(\vec{x},0\right) & =\varphi_{1}\left(\vec{x}\right),\quad\vec{x}\in\Omega.
\end{align}
For this ordinary differential equation the two stable equilibrium
points are ${\displaystyle \pm\sqrt{\frac{\mu^{2}}{\lambda}}=\pm\sqrt{\frac{9}{2}}\approx\pm2.12}$,
and the zero is an unstable equilibrium point (see Figure \ref{fig:Duffing}
for a phase portrait of equation \eqref{eq:duffing}). Comparing these
values to parts (e)-(f) of Figure \ref{fig:Bubble-long-time} we can
observe that after some initial time during which the dissipative
term is not negligible, eventually the solution of the Higgs boson
equation \eqref{eq:higgs-boson-rescalled} converges to a step function
corresponding to the positive and negative equilibrium points of the
Duffing equation \eqref{eq:duffing}. It is our future plan to investigate
and mitigate the numerical difficulty to correctly simulate the places
with the sudden changes in the solution. Remarkably, our numerical
method did not break down for the time intervals we present here.
We will discuss the connection to Duffing equations further in the
next examples. 
\end{example}
\begin{example}
\label{exa:no-bubble}In this example we demonstrate that under some
conditions bubbles do not form. This is a new conjecture supported
by our computational result, with currently no theoretical proof existing
for it in the literature. The parameter values are $\mu^{2}=9$, $\lambda=2$
as before, and the initial conditions are
\begin{align}
{\displaystyle \varphi_{0}\left(\vec{x}\right)} & =3B\left(x;\left(0.5,0.5,0.5\right),0.3\right),\qquad\forall\vec{x}\in\Omega\label{eq:init0-nobubble}\\
\varphi_{1}\left(\vec{x}\right) & =0,\qquad\forall\vec{x}\in\Omega.
\end{align}
Figure \ref{fig:nobubble} shows the line plot of the solution for
various times. There are no zeros inside the domain of support, which
means that no bubbles formed. Note that this initial data lies on
the nonnegative horizontal axis of the phase portrait shown in Figure
\ref{fig:Duffing}, hence, other than the unstable zero equilibrium,
this initial data is in the positive equilibrium's basin of attraction
for all $\vec{x}$ in the domain of support. Hence, it should be expected
that the solution converges to a step function with values zero and
${\displaystyle \sqrt{\frac{\mu^{2}}{\lambda}}\approx2.12}$, as suggested
by the dynamical properties of the corresponding damped, unforced
Duffing equation. The main effect the initially non-negligible dissipative
term in equation \eqref{eq:higgs-boson-rescalled} seems to have is
to propagate the values in space while they converge to the step function's
values. 
\end{example}
\begin{example}
\label{exa:duffing}This example shows a slightly more complicated
combination of initial conditions in order to demonstrate the connection
between the asymptotic behavior of solutions to the Higgs boson equation
\eqref{eq:higgs-boson-rescalled} and the associated damped, unforced
Duffing equation \eqref{eq:duffing}. The parameter values are $\mu^{2}=9$,
$\lambda=2$ as before, and the initial conditions are
\begin{align}
\varphi_{0}\left(\vec{x}\right) & =-10B\left(\vec{x};\left(0.5,0.5,0.5\right),0.3\right)B\left(\vec{x};\left(0.55,0.55,0.55\right),0.3\right)\sin\left(2\pi x\right),\qquad\forall\vec{x}\in\Omega,\label{eq:B1timesB2}\\
\varphi_{1}\left(\vec{x}\right) & =5B\left(\vec{x};\left(0.5,0.5,0.5\right),0.3\right),\qquad\forall\vec{x}\in\Omega.
\end{align}
Note that in the function $\sin\left(2\pi x\right)$ of initial condition
\eqref{eq:B1timesB2} the variable $x$ is the first component of
the thee-dimensional space variable $\vec{x}$. Part (b) of Figure
\ref{fig:duffing_complex} shows the line plot of the initial function
$\varphi_{0}$. Part (a) of Figure \ref{fig:duffing_complex} depicts
the curve $\left(\varphi_{0}\left(x,0.5,0.5\right),\varphi_{1}\left(x,0.5,0.5\right)\right)$
for $x\in\left[0,1\right]$ in the phase portrait of the Duffing equation.
Note also that this initial data cannot be made radial by a simple
shift in the spatial variable. On the other hand, this initial data
lies in the damped, unforced Duffing equation's stable positive equilibrium's
basin of attraction (and the unstable zero equilibrium point), even
though $\varphi_{0}$ changes sign. Part (a) of Figure \ref{fig:duffing_complex}
also suggests that the negative part of the initial function $\varphi_{0}$
has to be small relative to the positive part of $\varphi_{1}$ in
order for the initial data to stay in the positive equilibrium's basin
of attraction. We can observe on the line plot of solution for time
$t=4$ in part (c) of Figure \ref{fig:duffing_complex} that the nonzero
part of the solution becomes and stays positive for larger time and
indeed it converges to the positive equilibrium point ${\displaystyle \sqrt{\frac{\mu^{2}}{\lambda}}=\sqrt{\frac{9}{2}}\approx2.12}$
of the damped, unforced Duffing equation, as expected.
\end{example}
\begin{example}
In our final example we look at a more complicated scenario with initial
conditions that are the sums of two bump functions: 
\begin{align}
\varphi_{0}\left(\vec{x}\right) & =B\left(\vec{x};\left(0.4,0.4,0.4\right),0.2\right)+B\left(\vec{x};\left(0.6,0.6,0.6\right),0.2\right),\qquad\forall\vec{x}\in\Omega,\label{eq:B1plusB2}\\
\varphi_{1}\left(\vec{x}\right) & =-5\varphi_{0}\left(\vec{x}\right).\qquad\forall\vec{x}\in\Omega.
\end{align}
The parameter values are $\lambda=\mu^{2}=0.1$. Figure \ref{fig:two-bubbles-init}
shows a line plot of the first initial data $\varphi_{0}\left(\vec{x}\right)$,
which cannot be made radial by a simple shift in the spatial variable.
Initially the dissipative term is large compared to the speed at which
the solution converges to the Duffing equation's step function. As
a result two bubbles form and they interact with each other. For line
plots we used the main diagonal line segments of the computational
cube, since the centers of the initial conditions's bump functions
are on that diagonal. Figures \ref{fig:two-bubbles-2}-\ref{fig:two-bubbles-4}
show the formation and interactions of bubbles. Initially there is
no bubble present. After the two bubbles form at around $t=0.08$,
their size grows continuously in time. Around time $t=0.69$ the two
bubbles touch, and from that time on they are attached to each other.
At time $t=0.8$ (shown on parts (a.1) and (a.2) of Figure \ref{fig:two-bubbles-3}\foreignlanguage{american}{)}
an additional, tiny bubble forms inside each of the now merged bubbles.
These additional bubbles grow (parts (b.1) and (b.2) of Figure \ref{fig:two-bubbles-3}
at time $t=1$); then they flatten and become concave (part (c.1)
and (c.2) of Figure \ref{fig:two-bubbles-3} at time $t=2$). Later
hole forms in these inside bubbles, and they become toroidal (parts
(a.1) and (a.2) of Figure \ref{fig:two-bubbles-4} at time $t=2.15$).
Finally they disappear (parts (b.1) and (b.2) of Figure \ref{fig:two-bubbles-4}
at time $t=3$). The growth of the larger outer bubble exponentially
slows down and it does not seem to change its shape after time $t=3$.
Parts (c.1) and (c.2) of Figure \ref{fig:two-bubbles-4} show the
bubble in a quasi-steady state due to the dissipative term in Equation
\eqref{eq:higgs-boson-rescalled} being very close to zero and the
solution very close to the step function with constant values $0$
and ${\displaystyle \pm\sqrt{\frac{\mu^{2}}{\lambda}}=\pm1}$ of the
corresponding damped, unforced Duffing equation's equilibrium points. 
\end{example}

\section{\label{sec:Concluding-Remarks}Conclusion}

In this paper, we obtained numerical solutions of the Higgs boson
equation in the de Sitter spacetime, which is a nonlinear hyperbolic
partial differential equation on three spatial dimensions and one
time dimension. Our approach was based on a fourth order finite difference
method in space and an explicit fourth order Runge-Kutta method in
time for the discretization. High performance computations performed
using NVIDIA Tesla K40c GPU Accelerator demonstrated existing theoretical
results and suggested several previously not known properties of the
solution for which theoretical proofs do not exist yet. We conjecture
that for smooth initial conditions a smooth solution exists globally
for all time; under some initial conditions no bubbles form; and solutions
converge to step functions determined by the initial conditions and
related to unforced, damped Duffing equations. The theoretical proofs
of these conjectures is a future plan along with the numerical treatment
of the sharp corners developing in the solutions.

\section*{Acknowledgments}

The authors acknowledge the Texas Advanced Computing Center (TACC)
at The University of Texas at Austin for providing high performance
computing and visualization resources that have contributed to the
research results reported within this paper. URL: \url{http://www.tacc.utexas.edu}.
We also gratefully acknowledge the support of NVIDIA Corporation with
the donation of the Tesla K40 GPU used for this research.

\newpage{}

\bibliographystyle{amsplain}
\addcontentsline{toc}{section}{\refname}\bibliography{higgs}

\end{document}